\documentclass[11pt]{article}
 \usepackage{indentfirst, latexsym,bm}
 \usepackage{amsmath}
 \usepackage{pifont}
 \usepackage{amsfonts}
 \usepackage{mathrsfs}
 \usepackage{array}
 \usepackage{multirow} 
 \usepackage{graphicx}
 \usepackage{subfigure}
 \usepackage{picinpar}
 \RequirePackage[colorlinks,citecolor=blue,urlcolor=blue,linkcolor=blue]{hyperref} 
 
 \usepackage{authblk}
 
 \RequirePackage[numbers]{natbib}
 \usepackage{enumitem} 

 \usepackage{booktabs} 
 \usepackage{threeparttable}
 \usepackage{mathrsfs,amsfonts,amsmath,amssymb} 
 \usepackage{color}

 \makeatletter
 \def\namedlabel#1#2{\begingroup
  #2%
  \def\@currentlabel{#2}%
  \phantomsection\label{#1}\endgroup
  }
 \makeatother

 \newcommand\email[1]{\href{mailto:#1}{ \nolinkurl{#1}}}

 \renewcommand{\theequation}{\arabic{section}.\arabic{equation}}

 \newtheorem{theorem}{Theorem}[section]
 \newtheorem{definition}[theorem]{Definition}
 \newtheorem{lemma}[theorem]{Lemma}
 \newtheorem{corollary}[theorem]{Corollary}
 \newtheorem{proposition}[theorem]{Proposition}
 \newtheorem{remark}[theorem]{Remark}
 \newtheorem{condition}[theorem]{Condition}
 \newtheorem{example}{Example}[section]
 
 \newtheorem{assumption}[theorem]{Assumption}

 \def\blemma{\begin{lemma}}
 	\def\elemma{\end{lemma}}
 \def\bproposition{\begin{proposition}}
 	\def\eproposition{\end{proposition}}
 \def\ttheorem{\begin{theorem}}
 	\def\etheorem{\end{theorem}}
 \def\bcorollary{\begin{corollary}}
 	\def\ecorollary{\end{corollary}}
 \def\bremark{\begin{remark}}
 	\def\eremark{\end{remark}}
 \def\bcondition{\begin{condition}}
 	\def\econdition{\end{condition}}

 \def\benumerate{\begin{enumerate}}
    \def\eenumerate{\end{enumerate}}
 \def\bitemize{\begin{itemize}}
    \def\eitemize{\end{itemize}}

 \def\beqlb{\begin{eqnarray}}
 	\def\eeqlb{\end{eqnarray}}
 \def\beqnn{\begin{eqnarray*}}
 	\def\eeqnn{\end{eqnarray*}}
 \def\ar{\!\!\!&}

 \def\proof{\noindent{\it Proof.~~}}
    \def\qed{\hfill$\Box$\medskip}

\begin{document}

 \title{\bf\Large  Scaling Limit Theorems for Multivariate Hawkes Processes and Stochastic Volterra Equations with Measure Kernel}
 \author{Wei Xu\footnote{School of Mathematics and Statistics, Beijing Institute of Technology, 100081 Beijing, China. Email: xuwei.math@gmail.com}}

 \maketitle

 \begin{abstract}
 	
 This paper is devoted to establishing the full scaling limit theorems for multivariate Hawkes processes.
 Under some mild conditions on the exciting kernels,   we develop a new way to prove that after a suitable time-spatial scaling, the asymptotically critical multivariate Hawkes processes converge weakly to the unique solution of a multidimensional stochastic Volterra equation with convolution kernel being the potential measure associated to a matrix-valued extended  Bernstein function.    
 Also, based on the observation of their affine property and generalized branching property, we provide an exponential-affine representation of the Fourier-Laplace functional of scaling limits in terms of the unique solutions of multidimensional Riccati-Volterra equations with measure kernel. 
 The regularity of limit processes and their alternate representations are also investigated by using the potential theory of L\'evy subordinators.
 
 \bigskip \medskip
 
 \noindent {\it \textbf{MSC 2020 subject classifications:}} Primary 60G55, 60F17; secondary 60H20, 45D05.
  
 \smallskip
  
 \noindent {\it \textbf{Keywords and phrases:}} Hawkes process, scaling limit, (stochastic) Volterra equation, measure kernel, extended Bernstein function.

 \end{abstract}
  	

 \section{Introduction} 
 \label{Sec.Introduction}
 \setcounter{equation}{0}
 
 For an integer $d\geq 1$ and $\mathtt{D}:= \{1,2,\cdots, d\}$, 
 a $d$-variate Hawkes process $\boldsymbol{N}:= \{(N_i(t))_{i\in\mathtt{D}} :t\geq 0\}$ is a $d$-dimensional random point process on $\mathbb{R}_+$ that models self- and mutually exciting arrivals of $d$ kinds of distinguishable random events, which are usually called type-$1,2,\cdots, d$.  
 Its intensity  at time $t$, denoted by $\boldsymbol{\varLambda}(t):= (\varLambda_i(t))_{i\in\mathtt{D}}$, is of the form
 \beqlb\label{Eqn.MultiHPdensity}
 \varLambda_i(t)
 := \mu_i(t) + \sum_{j=1}^d \sum_{k=1}^{N_j(t)} \phi_{ij}(t-\tau_{j,k})
 =\mu_i(t) + \sum_{j=1}^d  \int_{[0,t)} \phi_{ij}(t-s) dN_j(s), 
 \quad i \in \mathtt{D},
 \eeqlb
 for some \textsl{exogenous density} $\boldsymbol{\mu}:= (\mu_i)_{i\in\mathtt{D}}\in L^1_{\rm loc} (\mathbb{R}_+;\mathbb{R}_+^{d}) $ that represents the \textsl{external excitation} on the arrival of 
 future events of various types and  \textsl{kernel} $\boldsymbol{\phi}:=(\phi_{ij})_{i,j\in \mathtt{D}}\in L^1 (\mathbb{R}_+;\mathbb{R}_+^{d\times d})$ that captures the \textsl{self-} or \textsl{mutually exciting} impact of past events of various types on the arrivals of future events of various types.   
 It can be reconstructed as a Poisson cluster process associated to an age-dependent multi-type branching process (also well-known as \textsl{Crump-Mode-Jagers branching process with immigration}) with  mean offspring matrix
 \beqnn
 \|\boldsymbol{\phi} \|_{L^1} 
 =\big( \big\| \phi_{ij}\big\|_{L^1} \big)_{i,j\in\mathtt{D}} = \Big( \int_0^\infty \phi_{ij}(t)\, dt \Big)_{i,j\in\mathtt{D}}
 \eeqnn
 being identically the average number of child
 events triggered by a mother event of various types. 
 Depending on the \textsl{spectral radius} of $\|\boldsymbol{\phi}\|_{L^1}$, denoted by $\rho (\|\boldsymbol{\phi}\|_{L^1})$, three phases of $(\boldsymbol{N},\boldsymbol{\varLambda})$ arise from  the criticality of branching processes, i.e., it is said to be \textsl{subcritical}, \textsl{critical} or \textsl{supercritical} if $\rho (\|\boldsymbol{\phi}\|_{L^1})<1$, $=1$ or $>1$ respectively. 
 We refer to Section~2.1 in \cite{BacryMastromatteoMuzy2015,Xu2021} for details. 
 This paper is concerned with the asymptotic behavior of $(\boldsymbol{N},\boldsymbol{\varLambda})$ at a large time scale as it tends to be critical.
 
 \subsection{Literature review} 
 
 Since Hawkes processes were firstly introduced by  A.G. Hawkes \cite{Hawkes1971a,Hawkes1971b} to understand the cross-dependencies among earthquakes and their aftershocks, they have been widely used to model the cascade phenomenon and clustering effect that have been widely observed in  biology and neuroscience \cite{HorstXu2021,Johnson1996,Reynaud-BouretSophieSchbath2010}, sociology and criminology \cite{Crane2008,Mohler2011}, financial contagion  \cite{AitSahaliaCacho-DiazLaeven2015,JorionZhang2009},   limit order books \cite{HorstXu2019,Large2007} and stochastic volatility  \cite{ElEuchFukasawaRosenbaum2018,HorstXu2022}.  
 We refer to  \cite{BacryMastromatteoMuzy2015,HorstXu2023} for reviews on Hawkes processes and their applications. 
 In recent years, motivated by modeling the tangled and complicated self-/mutually excitation in various stochastic systems, Hawkes processes have been generalized in many directions, e.g. marked Hawkes processes \cite{Ogata1988}, nonlinear Hawkes processes \cite{BremaudMassoulie2001}, infinite-dimensional Hawkes processes \cite{HorstXu2019},  Hawkes processes with simultaneous occurrence \cite{BieleckiJakubowskiNieweglowski2022} and so on.

 Along with their growing applications in various fields,  \textsl{microscopic} and \textsl{macroscopic} properties of Hawkes processes and their generalized models have now been studied in a rapidly increasing theoretical literature. 
 For instance, many well-known microscopic properties, including statistical characterization, cluster representation, probability generating function, genealogy and event cascades, were established in \cite{BacryMuzy2014,GaoZhu2018c,Hawkes1971a,Hawkes1971b,HawkesOakes1974,HorstXu2022} for subcritical (generalized) Hawkes processes. 
 Their macroscopic properties are usually investigated  via three different kinds of limit theorems: \textsl{law of large numbers {\rm (LLN)} and central limit theorem} (CLT), \textsl{large deviation principle} (LDP) and \textsl{scaling limit theorem} (SLT). 
 Under some light-tailed conditions on the kernel, (functional) LLNs and CLTs were established in \cite{BacryDelattreHoffmannMuzy2013,HawkesOakes1974} for Hawkes processes and \cite{HorstXu2021, KarabashZhu2015} for marked Hawkes processes. 
 A nearly full functional LLN and CLT were established in \cite{HorstXu2023} recently for subcritical and critical Hawkes processes under either light-tailed or heavy-tailed conditions. 
 LPDs have been established in \cite{GaoZhu2018b} for Markovian Hawkes processes with large exogenous density, \cite{Zhu2014,Zhu2015} for nonlinear Hawkes processes and \cite{GaoGaoZhu2023,GaoZhu2023} for their mean-field limit. 
 
 SLTs for nearly unstable uni-variate Hawkes processes were firstly established by Jaisson and Rosenbaum \cite{JaissonRosenbaum2015,JaissonRosenbaum2016} in the study of asymptotic behavior of Hawkes-based price models in the context of high-frequency trading. 
 Under a light-tailed condition on the kernel, they proved that the rescaled intensity and point process converge weakly to a Feller diffusion and the relevant integrated process. 
 Under a heavy-tailed condition, they also established the weak convergence of rescaled point process to the integral of a rough fractional diffusion. 
 The weak convergence of rescaled intensities to the  rough fractional diffusion was proved in a recent work \cite{HorstXuZhang2023a}. 
 Analogous results in the multivariate case were provided in \cite{ElEuchFukasawaRosenbaum2018,RosenbaumTomas2021}.
 For marked Hawkes processes, SLTs were established in \cite{HorstXu2022} and \cite{Xu2021}  respectively for uni-variate case and multivariate case under a light-tailed condition on the kernel.  
 To generalize the second Ray-Knight theorem for spectrally postive L\'evy processes, Xu \cite{Xu2025,Xu2024b} established a SLT for a class of heavy-tailed marked Hawkes processes and proved the weak convergence of their rescaled intensities to the unique solution of a stochastic Volterra equation driven by a Gaussian white noise and a Poisson random measure. 
 Recently, Horst et al. \cite{HorstXuZhang2023b} considered a bench market model with market order and limit order modeled respectively by a heavy-tailed Hawkes process and a long-range dependent marked Hawkes process. They proved that after a suitable time-spatial scaling the volatility process converges to a fractional Heston model driven by an additional Poisson random measure. 
 
 In the existing literature, limiting processes of Hawkes-based models are usually either (drifted) Brownian motions and Markovian affine-like processes in the light-tailed case or Gaussian processes of Riemann-Liouville type and rough fractional affine-like processes in the heavy-tailed case. 
 They usually can be characterized by the following multi-dimensional stochastic Volterra equation (SVE) of convolution type 
 \beqlb\label{General.SVE}
 \boldsymbol{X}(t)= \boldsymbol{G}(t) 
 + \int_0^t \boldsymbol{F}(t-s) \boldsymbol{\sigma} \big(s,\boldsymbol{X}(s)\big)\, d\boldsymbol{Z}(s),\quad t\geq 0,
 \eeqlb
 for some given function $\boldsymbol{G}$,  convolution kernel $\boldsymbol{F}$, coefficient $\boldsymbol{\sigma}$ and semimartingale $\boldsymbol{Z}$. 
 Non-singular Lipschitz SVEs driven by Brownian noise were studied firstly at least in \cite{BergerMizel1980a,BergerMizel1980b,RaoTsokos1975} and later generalized to the case of general semimartingales in \cite{Protter1985} and anticipating coefficients in \cite{PardouxProtter1990}. 
 The non-Lipschitz SVEs driven by Brownian noise and Poisson random measure were investigated in two recent works \cite{AlfonsiSzulda2024,PromelScheffels2023}. 
 The study of SVEs with singular kernels and (near-)Lipschitz coefficients can be found in \cite{CoutinDecreusefond2001,Decreusefond2002,Wang2008,Zhang2010}. 
 The existence and uniqueness of solutions to \eqref{General.SVE} with (near-)Lipschitz coefficients were usually  proved by using the standard Picard iteration in the aforementioned literature.

 By contrast, the well-posedness of  \eqref{General.SVE} is considerably challenging in the case of non-Lipschitz coefficients and singular kernels, due to the failure of Picard iteration as well as the loss of martingality and Markovianity.  
 To the best of our knowledge, the weak well-posedness of one-dimensional continuous SVE \eqref{General.SVE} with $(1/2+\varepsilon)$-H\"older coefficient and fractional convolution kernel was firstly considered by Mytnik and Salisbury in the pioneering work \cite{MytnikSalisbury2015}. 
 Nowadays, singular SVEs with non-Lipschitz coefficients have attracted a considerable renewed interests in probability theory and mathematical finance because of their striking success in modeling stochastic rough volatility and large network; see e.g.  \cite{DelattreFournierHoffmann2016,ElEuchFukasawaRosenbaum2018}.   
 Abi Jaber et al. \cite{Jaber2021,JaberLarssonPulido2019} established the weak well-posedness of \eqref{General.SVE} with affine coefficients as well as the kernel being  continuous, non-negative, non-increasing and admitting a non-negative, non-increasing resolvent of the first kind.  
 Later, the weak existence of $L^p$-solutions to \eqref{General.SVE} was proved by them  \cite{JaberCuchieroLarssonPulido2021} under some linear growth conditions on the coefficients and some  regularity conditions on the kernel. 
 Their methods were further extended in \cite{PromelScheffels2023b,PromelScheffels2024} to deal with the general cases. 
 It is worth mentioning that the strong well-posedness of \eqref{General.SVE} with singular kernel and $1/2$-H\"older coefficient is still open up to now. 
 
 \subsection{Our contributions and methodologies}
 
 In this paper, we establish the full SLTs for multivariate Hawkes processes that are asymptotically critical and characterize their scaling limits in various ways. 
 Consider a sequence of $d$-variate Hawkes processes $\{(\boldsymbol{N}_n,\boldsymbol{\varLambda}_n)\}_{n\geq 1}$ with parameters $\{ (\boldsymbol{\mu}_n,\boldsymbol{\phi}_n) \}_{n\geq 1}$  satisfying that $\rho\big(\|\boldsymbol{\phi}_n \|_{L^1}\big)\to 1$, which particularly reduces to $\|\phi_n  \|_{L^1} \to 1$ in the uni-variate case. 
 However, the infinitely many possible limit matrices with unit spectral radius and unaccountable number of ways to tend towards criticality in the multivariate case give rise to the first challenge, that is that the asymptotic criticality is much more complicated than that in the uni-variate case.  
 So necessarily, associated with a sequence $\{\theta_n\}_{n\geq 1}$ increasing to infinity and a non-negative matrix $\boldsymbol{K}$ with  $\rho (\boldsymbol{K} ) =1 $ we specify our \textsl{asymptotical criticality condition} as follows
 \beqnn
  \big\|\boldsymbol{\phi}_n \big\|_{L^1} \to \boldsymbol{K} 
  \quad \mbox{and} \quad
  \sqrt{n\cdot \theta_n} \cdot  \big( \boldsymbol{K}-  \mathcal{L}_{\boldsymbol{\phi}_n} (\lambda/n) \big) \to \boldsymbol{\varPhi}(\lambda),\quad \lambda \geq 0,
 \eeqnn
 where  $\mathcal{L}_{\boldsymbol{\phi}_n}$ represents the Laplace transform of $\boldsymbol{\phi}_n$ and  $\boldsymbol{\varPhi}$ is a continuous function on $\mathbb{R}_+$ that is proved to be an extended Bernstein function with L\'evy triplet $(\boldsymbol{b}^{\boldsymbol{\varPhi}},\boldsymbol{\sigma}^{\boldsymbol{\varPhi}},\boldsymbol{\nu}^{\boldsymbol{\varPhi}})$. 
 Additionally, we also assume the limit couple $(\boldsymbol{K},\boldsymbol{\varPhi})$ satisfying \textsl{the admissibility condition} (see Definition~\ref{Def.Admissible}) to guarantee that the scaling limits are non-degenerated. 
 In particular, the admissibility condition in uni-variate case reduces to be $\varPhi(\infty)>0$, that is identified to be the necessary and sufficient condition for the convergence of rescaled Hawkes processes. 
  
 Our first main result in this paper establishes the joint convergence of Hawkes point processes and their compensated point processes after a suitable time and spatial scaling,  i.e., 
 \beqlb\label{eqn.1002}
 \big(\boldsymbol{N}^{(n)}(t) , \boldsymbol{\widetilde{N}}^{(n)}(t) \big)
 :=\bigg( \frac{\boldsymbol{N}_n(nt)}{n\cdot \theta_n} , \frac{\boldsymbol{\widetilde{N}}_n(nt)}{\sqrt{n\cdot \theta_n}} \bigg) 
 \to \big(   \boldsymbol{\varXi},  \boldsymbol{M}  \big)
 \eeqlb
 weakly in $D(\mathbb{R}_+;\mathbb{R}_+^d\times \mathbb{R}^d)$ with the limit process $\boldsymbol{\varXi}$ being non-decreasing and $ \boldsymbol{M}$ being a martingale. 
 Usually, the standard machinery to establish the weak convergence of stochastic processes  consists in \textsl{proving tightness} and \textsl{characterizing limit points} that are remarkable challenges in the non-Markovnian or non-martingale framework, 
 because very few instruments provided by modern probability theory and stochastic analysis (see \cite{Billingsley1999,EthierKurtz1986,JacodShiryaev2003}) are available and effective. 
 Even worse, due to the long-range dependence among events occurring at different time, the mutual excitation among events of various types and the lack of high-order moment estimates, 
 it turns to be more difficult and challenging to establishing the weak convergence for multivariate Hawkes processes. 
 
 Instead of assuming some additional stronger moment or regularity conditions on the pre-limit models as in the aforementioned references,
 we develop a new way to overcome the preceding difficulties:
 \begin{enumerate}
 	\item[$\bullet$] To prove the tightness, we first endow the function space $D(\mathbb{R}_+;\mathbb{R}_+^d\times \mathbb{R}^d)$ with the \textsl{$S$-topology} that was introduced by Jakubowski \cite{Jakubowski1997} and is a little weaker than the Skorokhod topology but much stronger than the pseudo-path topology and also the topology generated by the finite-dimensional convergence. 
 	Then, we prove the relative compactness of the sequence 
 	$ \big\{ (\boldsymbol{N}^{(n)},\boldsymbol{\widetilde{N}}^{(n)} ) \big\}_{n\geq 1}$ in the sense of $S$-topology by identifying the stochastic boundedness of the numbers of their oscillations. 
 	
 	\item[$\bullet$]  To achieve the uniqueness in law of accumulation points, we first use the duality method developed in \cite{JaberLarssonPulido2019,HorstXu2023} to establish an exact exponential-affine representation of the Fourier-Laplace functional 
 	\beqnn
 	\mathbf{E}\Big[  \exp\Big\{ \int_{[0,t]}    \boldsymbol{f}(t-s)  \boldsymbol{N}^{(n)}(ds) + \int_{[0,t]}    \boldsymbol{h}(t-s)  \boldsymbol{\widetilde{N}}^{(n)}(ds) \Big\} \Big]
 	\eeqnn
 	in terms of the unique solution of a nonlinear Volterra integral equation. 
 	In collaboration with a priori uniform estimate for solutions of these nonlinear Volterra integral equations, we prove their uniform convergence to the following Riccati-Volterra equation with measure kernel 
 	\beqlb\label{eqn.1001}
 	\boldsymbol{\mathcal{V}}(t)  =   \int_{[0,t]}\boldsymbol{\mathcal{W}}(t-s)\, \boldsymbol{\varPi}(ds)  
 	\quad \mbox{with} \quad
 	\boldsymbol{\mathcal{W}}(t):=  \boldsymbol{f}(t)  + \frac{1}{2} \cdot \big(\boldsymbol{\mathcal{V}}(t)+\boldsymbol{h}(t)\big)^2.
 	\eeqlb
 	where $\boldsymbol{\varPi}(ds)$ is the potential measure associated to $(\boldsymbol{K},\boldsymbol{\varPhi})$. 
 	Finally, we obtain the uniqueness in law of  accumulation points by proving that the uniqueness holds for \eqref{eqn.1001} and the Fourier-Laplace functional of every accumulation point $ (\boldsymbol{\varXi},  \boldsymbol{M}  )$ admits the following common representation 
 	\beqlb\label{eqn.1006}
 	\mathbf{E}\Big[\exp \Big\{ \int_{[0,t]}\boldsymbol{f} (t-s)\, d\boldsymbol{\varXi} (t) + \int_{[0,t]} \boldsymbol{h}(t-s)\, d\boldsymbol{M} (s)   \Big\} \Big]
 	= \exp \Big\{ \int_{[0,t]}\boldsymbol{\mathcal{W}}(t-s)\, d\boldsymbol{\varUpsilon} (s) \Big\} 
 	\eeqlb
 	with $\boldsymbol{\varUpsilon}$ being a non-decreasing function that captures the accumulated impact of exogenous events on the limit system. 
 	
 	\item[$\bullet$]  Based on the preceding results, we finally obtain the weak convergence \eqref{eqn.1002} in the sense of the $S$-topology, and, further, improve it to the sense of the Skorokhod topology, provided that either $\boldsymbol{\varXi}$ or $\boldsymbol{M}$ is continuous. 
 \end{enumerate}
 It is worth mentioning that due to the possible atoms in the measure $\boldsymbol{\varPi}(ds)$, the well-posedness of \eqref{eqn.1001} seems impossible to be established by using the standard Picard iteration and Banach's fixed point theorem. 
 For all we know, it seems the first time to obtain the existence and uniqueness of solutions to Riccati-Volterra equations with measure kernel.
 
 As the second main result in the work, we provide another characterization for the limit process  $(\boldsymbol{\varXi},\boldsymbol{M})$ in terms of the unique solution to the following Volterra integration equation 
 \beqlb\label{eqn.1003}
 \boldsymbol{\varXi}(t)= \boldsymbol{\varUpsilon}(t)+ \int_{[0,t]} \boldsymbol{\varPi}(ds)\,\boldsymbol{M}(t-s).
 \eeqlb
 Provided the continuity of $\boldsymbol{\varXi}$ or $\boldsymbol{M}$, the martingale $\boldsymbol{M}$ has quadratic variation $\boldsymbol{\varXi}$ and can be represented as $\boldsymbol{M} = \boldsymbol{B} \circ \boldsymbol{\varXi} $ for some standard $d$-dimensional Brownian motion $\boldsymbol{B}$.
 Additionally, when $\boldsymbol{\varUpsilon}$ is differentiable and the potential measure $\boldsymbol{\varPi}(ds)$ has a locally square integrable density $\boldsymbol{\pi}$ (no other conditions are necessary), the limit process $\boldsymbol{\varXi}$ has a predictable derivative that uniquely solves the following SVE
 \beqnn
 	\boldsymbol{\xi} (t) =  \boldsymbol{\varUpsilon}'(t)
 	+ \int_0^t \boldsymbol{\pi}(t-s)\, {\rm diag } \big(\sqrt{\boldsymbol{\xi} (s)}\big)d\boldsymbol{B}(s)   . 
 \eeqnn 
 The challenge to establish the preceding equations mainly comes from the lake of a priori estimates on the high-order moments and oscillations of pre-limit models. 
 The problem is overcome via establishing the weak convergence results for stochastic Volterra integrals by applying the theory of weak convergence of It\^o's stochastic integrals developed in \cite{Jakubowski1996} and then proving the weak convergence of SVEs solved by rescaled Hawkes point processes to \eqref{eqn.1003}.  
 
 The last contribution we make in this work is to 
 provide the following alternate stochastic equation for $(\boldsymbol{\varXi}, \boldsymbol{M})$ that separates the impact of drift parameter $\boldsymbol{b}^{\boldsymbol{\varPhi}} $ from those of  other two parameters $ \boldsymbol{\sigma}^{\boldsymbol{\varPhi}}$ and $\boldsymbol{\nu}^{\boldsymbol{\varPhi}} $ 
 \beqlb \label{eqn.1004}
 \boldsymbol{\varXi}(t) \ar=\ar  \int_{[0,t]} \boldsymbol{\varPi}_0 (ds)  \boldsymbol{\varGamma}(t-s) - \int_{[0,t]} \boldsymbol{\varPi}_0 (ds)   \big(\boldsymbol{b}^{\boldsymbol{\varPhi}}\cdot \boldsymbol{\varXi}(t-s) \big) + \int_{[0,t]} \boldsymbol{\varPi}_0 (ds)   \boldsymbol{M}(t-s) 
 \eeqlb
 with $\boldsymbol{\varPi}_0(ds)$ being the potential measure associated to the extended Bernstein function with L\'evy triplet $(\mathbf{0}, \boldsymbol{\sigma}^{\boldsymbol{\varPhi}},\boldsymbol{\nu}^{\boldsymbol{\varPhi}})$. 
 This representation also provides a criticality criterion for the process $\boldsymbol{\varXi}$ according as the spectral radius of $\boldsymbol{b}^{\boldsymbol{\varPhi}}$. 
 Our method is mainly inspired by the potential theory of Markov processes. 
 In detail, we first establish a resolvent equation analogous to that of L\'evy subordinators to connect the two measures $\boldsymbol{\varPi}(ds)$ and $\boldsymbol{\varPi}_0(ds)$, and then identify the mutual equivalence between \eqref{eqn.1003} and \eqref{eqn.1004} by solving the relevant Wiener-Hopf type equations.

 \subsection{Conclusion} 
  
 Each non-degenerated scaling limit of multivariate Hawkes processes uniquely solves the SVE  \eqref{eqn.1003} that are fully determined by the potential measure $\boldsymbol{\varPi}(ds)$ and hence the matrix $\boldsymbol{K}$ and the matrix-valued extended Bernstein function $\boldsymbol{\varPhi}$. 
 Conversely, for any given admissible couple $(\boldsymbol{K},\boldsymbol{\varPhi})$, we also prove that there exists a sequence of scaling parameters and a sequence of multivariate Hawkes processes that are asymptotically critical such that the weak convergence \eqref {eqn.1002} holds with limit process uniquely solving the SVE  \eqref{eqn.1003}. 
 In conclusion, there exists an one-to-one correspondence among scaling limits of multivariate Hawkes processes and the SVEs \eqref{eqn.1003} as well as the admissible couples $(\boldsymbol{K},\boldsymbol{\varPhi})$. 
 Drawing from the theory of Volterra equations, we may say the limit process $(\boldsymbol{\varXi},\boldsymbol{M})$ is a \textsl{$(\boldsymbol{K},\boldsymbol{\varPhi})$-Hawkes Volterra process} and refer to \eqref{eqn.1003} as  \textsl{$(\boldsymbol{K},\boldsymbol{\varPhi})$-Hawkes Volterra equation}. 
 
 As we mentioned before, multivariate Hawkes processes are multi-type general branching processes with immigration and enjoy the \textsl{``branching property"}, which naturally should be inherited by the limit process $(\boldsymbol{\varXi},\boldsymbol{M})$. Indeed, this can be identified by comparing the Fourier-Laplace functional \eqref{eqn.1006} with those of continuous-state branching processes with immigration; see Theorem~1.1 in \cite{KawazuWatanabe1971}. 
 Hence all scaling limits of multivariate Hawkes processes constitute a family of \textsl{general (non-Markovian) continuous-state branching processes with immigration} whose \textsl{branching mechanisms} are determined by the relevant couple $(\boldsymbol{K},\boldsymbol{\varPhi})$. 
 Furthermore, the comparison between the two Fourier-Laplace functional \eqref{eqn.1006} and (2.2) in \cite{DuffieFilipovicSchachermayer2003} as well as the two equations \eqref{eqn.1001} and (2.14)-(2.16) in \cite{DuffieFilipovicSchachermayer2003} allows us to claim that the limit process $(\boldsymbol{\varXi},\boldsymbol{M})$ shares the \textsl{``affine property"} analogous to that of affine processes. 
 Hence all scaling limits of multivariate Hawkes processes also can be called \textsl{affine Volterra processes} as in \cite{JaberLarssonPulido2019}.  
 It is necessary to remind that our SVEs cannot be covered by models studied in the existing literature, e.g. \cite{Jaber2021,JaberCuchieroLarssonPulido2021,JaberLarssonPulido2019,PromelScheffels2023b,PromelScheffels2023,PromelScheffels2024}.


 \medskip
 
 \textit{\textbf{Organization of this paper.}} In Section~\ref{Sec.HPSL}, we first recall some elementary properties of multivariate Hawkes processes and then formulate our main results as well as some specific examples of processes obtained through our SLTs. Section~\ref{Sec.AsymResolvent} is devoted to provide some asymptotic results for the rescaled resolvent that will play an important role in the proofs of our main results. 
 The explicit exponential-affine representations of Fourier-Laplace functionals of Hawkes processes and their scaling limits are established in Section~\ref{Sec.ConvergenceFLF}. 
 All detailed proofs of our main theorems and corollaries are given in Section~\ref{ProofMainThm}. 
 The basic theories of non-negative matrix,  Volterra equation and $S$-topology are recalled respectively in Appendix~\ref{Sec.Matrices}, \ref{Appendix-Volterra} and \ref{Sec.TopologiesD}. 
 Some supporting results about extended Bernstein functions are presented in Appendix~\ref{Appendix--EBF}.

 \medskip
 
 \textit{\textbf{Notation.}} Every matrix/vector is denoted by a bold symbol and its elements are represented by the corresponding non-bold symbol with subscript. In the one-dimensional case, all bold symbols reduce to the corresponding non-bold ones without subscript.  
 Without specifying the dimension whenever there is no ambiguity, we write $\mathbf{1}$ for the all-ones row vector, $\mathbf{0}$ for the zero matrix and $\mathbf{Id}$ for the identity matrix.  
 Let $\boldsymbol{A}^{\rm T}$, $\boldsymbol{A}^{-1}$, $\rho (\boldsymbol{A})$ and ${\rm det}(\boldsymbol{A})$ denote the transpose, inverse, spectral radius and determinant of matrix $\boldsymbol{A}$ respectively whenever they exist. 
 
 For two integers $m,n \in \mathbb{Z}_+$, $p\in(0,\infty]$ and an interval $\mathcal{T}\subset[0,\infty)$, we define several function spaces and measure spaces as follows:
 
 \begin{enumerate}
 	\item[$\bullet$] $D(\mathcal{T},\mathbb{C}^{m\times n})$: the space of all c\`adl\`ag $\mathbb{C}^{m\times n}$-valued functions on $\mathcal{T}$.
 	
 	\item[$\bullet$] $C(\mathcal{T},\mathbb{C}^{m\times n})$: the space of all continuous $\mathbb{C}^{m\times n}$-valued functions on $\mathcal{T}$
 	
 	\item[$\bullet$]  $L^p(\mathcal{T}; \mathbb{C}^{m\times n})$: the space of all $\mathbb{C}^{m\times n}$-valued measurable functions $\boldsymbol{f}$ on $\mathcal{T}$ satisfying   
 	\beqnn
 	\big\|\boldsymbol{f}\big\|_{L^p_\mathcal{T}}^p:=  \Big(\int_\mathcal{T} \big|f_{ij}(t)\big|^p dt \Big)_{i,j\in\mathtt{D}} <\infty.
 	\eeqnn
    Let $L^p_{\rm loc}(\mathbb{R}_+; \mathbb{C}^{m\times n}) :=  \cap_{T\geq 0}L^p([0,T]; \mathbb{C}^{m\times n})$. 	
    We also write $\|\boldsymbol{f}\|_{L^p_T}$ for $\|\boldsymbol{f}\|_{L^p_{[0,T]}}$ and $\|\boldsymbol{f}\|_{L^{^p}}$ for $\|\boldsymbol{f}\|_{L^p_\infty}$. 
 	
 	\item[$\bullet$] $ M(\mathcal{T} ;\mathbb{C}^{m\times n})$: the space of all finite measures on $\mathcal{T}$ with range contains in $\mathbb{C}^{m\times n}$ and total variation norm. 
 	Let $ M_{\rm loc}(\mathbb{R}_+;\mathbb{C}^{m\times n}):= \cap_{T\geq 0} M([0,T] ;\mathbb{C}^{m\times n}) $.
 	
 \end{enumerate}
 
 We make the convention that for $x, y\in \mathbb{R}$ with $y\geq x$,
 \beqnn
 \int_x^y=-\int_y^x= \int_{(x,y]},\quad \int_{x-}^{y-}=  \int_{[x,y)}
 \quad\mbox{and}\quad
 \int_x^\infty = \int_{(x,\infty)}.
 \eeqnn 
 For $\boldsymbol{f} \in L^1_{\rm loc}(\mathbb{R}_+;\mathbb{C}^{l\times m})$ and  $\boldsymbol{g} \in L^1_{\rm loc}(\mathbb{R}_+;\mathbb{C}^{m\times n})$, their convolutions
 $\boldsymbol{f}* \boldsymbol{g}$ and $\boldsymbol{f}* d\boldsymbol{g}$ are defined by
 \beqnn
 \boldsymbol{f}* \boldsymbol{g}(t) 
 \ar:=\ar \bigg( \sum_{k=1}^m \int_{[0,t]} f_{ik}(t-s)g_{kj}(s)ds \bigg)_{i\leq l,j\leq n}
 \quad \mbox{and}\quad 
 \boldsymbol{f}* d\boldsymbol{g}(t) 
 := \bigg( \sum_{k=1}^m  \int_{[0,t]} f_{ik}(t-s)dg_{kj}(s) \bigg)_{i\leq l,j\leq n} ,
 \eeqnn
 for $t\geq 0$. 
 In particular, if $n=m$ and  $\boldsymbol{g}\equiv \mathbf{Id}$, then   $\mathcal{I}_{\boldsymbol{f}}:= \boldsymbol{f}* \mathbf{Id}$  turns to be the integral of $\boldsymbol{f}$. 
 For $\boldsymbol{\nu}\in M_{\rm loc}(\mathbb{R}_+;\mathbb{C}^{m\times n})$, the convolution $\boldsymbol{f}* \boldsymbol{\nu}$ is defined by 
 \beqnn 
 \boldsymbol{f}* \boldsymbol{\nu}(t) 
 \ar:=\ar \bigg( \sum_{k=1}^m  \int_{[0,t]} f_{ik}(t-s)\nu_{kj}(ds) \bigg)_{i\leq l,j\leq n} ,\quad t\geq 0. 
 \eeqnn
 
 Let $\overset{\rm f.d.d.}\longrightarrow$, $\overset{\rm a.e.}\to$,  $\overset{\rm a.s.}\to$, $\overset{\rm d}\to$  and  $\overset{\rm p}\to$ be the convergence in the sense of finite dimensional distributions,  almost everywhere convergence, almost sure convergence,  convergence in distribution and convergence in probability respectively.
 We also use $\overset{\rm a.s.}=$, $\overset{\rm a.e.}=$, $\overset{\rm d}=$ and $\overset{\rm p}=$ to denote almost sure equality, almost everywhere equality, equality in distribution and equality in probability respectively. 
 
 We use $C$ to denote a positive constant whose value might change from line to line.

  \section{Preliminaries and main results}
 \label{Sec.HPSL} 
 \setcounter{equation}{0}

 In this section we first recall some elementary properties of multivariate Hawkes processes and then formulate their scaling limit theorems. 
 Without loss of generality, we assume that all our processes are  defined on a complete probability space $(\Omega,\mathscr{F},\mathbf{P})$ endowed with a filtration $\{\mathscr{F}_t:t\geq 0\}$ that satisfies the usual hypotheses. 
 Without confusion, for a function $g$ on $\mathbb{C}$ and $\boldsymbol{x}\in \mathbb{C}^{m\times n}$ we make the convention that $g(\boldsymbol{x}):= \big(g(x_{ij})\big)_{i\leq m,j\leq n}$, e.g. $(\boldsymbol{x})^2:= (x_{ij}^2)_{i\leq m,j\leq n}$.  
 Also, we always extend each function $f$ on a subset $\mathcal{T}\subset\mathbb{R}$ to the whole real line by setting $f(z)=0$ if $z\in \mathbb{R}\setminus \mathcal{T}$.

 For convention, we set $\boldsymbol{\phi}(0)=\mathbf{0}$, replace the interval of integration $(0,t)$ in (\ref{Eqn.MultiHPdensity}) by $[0,t]$ and rewrite (\ref{Eqn.MultiHPdensity}) into the following matrix form 
 \beqlb
 \boldsymbol{\varLambda}(t) \ar=\ar \boldsymbol{\mu}(t) +  \boldsymbol{\phi} * d\boldsymbol{N} (t) ,\quad t\geq 0.
 \eeqlb
 Associated to the kernel $\boldsymbol{\phi}$ we define its \textsl{resolvent} $\boldsymbol{R}:=\big\{ \big( R_{ij}(t) \big)_{i,j\in \mathtt{D}}: t\geq 0 \big\} \in L^1_{\rm loc}(\mathbb{R}_+;\mathbb{R}_+^{d\times d})$  by the unique solution of the \textsl{resolvent equation} 
 \beqlb\label{Resolvent}
 \boldsymbol{R} (t)
 =\boldsymbol{\phi}(t) +\boldsymbol{R}*\boldsymbol{\phi}(t)  
 =\boldsymbol{\phi}(t) +\boldsymbol{\phi}*\boldsymbol{R}(t), \quad t\geq 0. 
 \eeqlb
 In addition to the resolvent $\boldsymbol{R}$, we will also need the following integrated process and integrated function
 \beqnn
 \mathcal{I}_{\boldsymbol{\varLambda}}(t):= \big\|\boldsymbol{\varLambda}\big\|_{L^1_t} = \int_{[0,t]} \boldsymbol{\varLambda}(s)ds
 \quad\mbox{and}\quad 
 \mathcal{I}_{\boldsymbol{R}}(t):= \big\|\boldsymbol{R} \big\|_{L^1_t} = \int_{[0,t]} \boldsymbol{R}(s)ds , \quad t\geq 0,
 \eeqnn
 to describe the \textsl{cumulative event rate} of the point process $\boldsymbol{N}$ and the \textsl{cumulative impact} of $d$ immigrant events of various types respectively.  
 In particular, the \textsl{total impact} $ \mathcal{I}_{\boldsymbol{R}}(\infty)$ is finite if and only if the spectral radius $\rho(\|\boldsymbol{\phi}\|_{L^1})<1$. In this case, we have 
 \beqnn
 \mathcal{I}_{\boldsymbol{R}}(\infty)=  \big\| \boldsymbol{R} \big\|_{L^1} 
 = \big( \mathrm{I}- \big\|\boldsymbol{\phi} \big\|_{L^1}  \big)^{-1}\cdot \big\|\boldsymbol{\phi} \big\|_{L^1} 
 = \big\|\boldsymbol{\phi} \big\|_{L^1} \cdot \big( \mathrm{I}- \big\|\boldsymbol{\phi} \big\|_{L^1}  \big)^{-1} . 
 \eeqnn 
 
 Note that the random point process $\boldsymbol{N}$ has compensator $\mathcal{I}_{\boldsymbol{\varLambda}}$ and hence the compensated point process 
 \beqlb\label{eqn.CompensatedN}
 \widetilde{\boldsymbol{N}}(t):=\boldsymbol{N}(t)-\mathcal{I}_{\boldsymbol{\varLambda}}(t),\quad t\geq 0,
 \eeqlb
 is a $(\mathscr{F}_t)$-martingale with \textsl{predictable quadratic co-variation} and \textsl{quadratic co-variation} (see \cite[Chapter~I.4]{JacodShiryaev2003})
 \beqlb\label{eqn.QV}
 \big\langle\widetilde{N}_i,\widetilde{N}_j\big\rangle_t = \mathbf{1}_{\{ i=j  \}}\cdot   \mathcal{I}_{\varLambda_i}(t)
 \quad \mbox{and}\quad 
 \big[\widetilde{N}_i,\widetilde{N}_j\big]_t= \mathbf{1}_{\{ i=j  \}}\cdot  N_i(t),
 \quad t\geq 0,\ i,j\in\mathtt{D}.  
 \eeqlb 
 The following useful martingale representation for the intensity process $\boldsymbol{\varLambda}$ comes from \cite{HorstXu2023}. 

 \begin{lemma}[Martingale representation]
 \label{Lemma.MartRep}
 The intensity process $\boldsymbol{\varLambda}$ is the unique solution to
 \beqlb\label{SVR}
 \boldsymbol{\varLambda}(t) \ar=\ar \boldsymbol{H}(t) +  \boldsymbol{R} * d\widetilde{\boldsymbol{N}} (t) 
 \quad \mbox{with}\quad 
 \boldsymbol{H}(t) :=  \boldsymbol{\mu} (t)  +   \boldsymbol{R} * \boldsymbol{\mu}(t) ,\quad t\geq 0.
 \eeqlb 
 \end{lemma}

 \subsection{Scaling limits} \label{Subsec.SLT}
 
 We now formulate our full scaling limit theorems for  multivariate Hawkes processes that are asymptotically critical.  
 The weak convergence results in this paper are mainly established in two function spaces $D_{J_1}(\mathbb{R}_+;\mathbb{C}^d)$ and $D_{S}(\mathbb{R}_+;\mathbb{C}^d)$ that are defined by endowing the space $D(\mathbb{R}_+;\mathbb{C}^d)$ with the $J_1$-topology (see \cite[Chapter~VI]{JacodShiryaev2003}) and the $S$-topology (see Appendix~\ref{Sec.TopologiesD}) respectively. 
 The $S$-topology is weaker than the $J_1$-topology but stronger than the pseudo-path topology (see \cite{DellacherieMeyer1975}).  
 Without further mention, we also endow  $C(\mathbb{R}_+;\mathbb{C}^d)$ with the uniform topology. 
 
 For each $n\in\mathbb{Z}_+$, we assume that the $n$-th process  $(\boldsymbol{N}_{n},\boldsymbol{\varLambda}_{n})$ has parameter $(\boldsymbol{\mu}_{n},\boldsymbol{\phi}_{n})$. 
 It is asymptotically critical if $\rho(\|\boldsymbol{\phi}_n  \|_{L^1})\to 1$ as $n\to\infty$. 
 In contrast to the uni-variate case in which the asymptotic criticality is equivalent to the $L^1$-norm of kernel goes to $1$, the convergence of $ \rho(\|\boldsymbol{\phi}_n  \|_{L^1})$ may not be inherited by $\|\boldsymbol{\phi}_n  \|_{L^1}$ and hence is not enough to study the long-run behavior of $(\boldsymbol{N}_{n},\boldsymbol{\varLambda}_{n})$. 
 Here we specify the asymptotic criticality as follows. 
 For each $r>0$ and $k \in\mathtt{D}$, let $\mathbb{R}_{+,r,k}^{d\times d}$ be the space of all non-negative $d\times d$-matrices with spectral radius $ r$ and the algebraic multiplicity of eigenvalue $r$ equal to $k$; see Proposition~\ref{Appendix.Prop.PFT}.
 \begin{assumption}[Asymptotic criticality]\label{AsymCriticality}
 For some $\ell \in\mathtt{D}$ and matrix $\boldsymbol{K} \in \mathbb{R}_{+,1,\ell}^{d\times d}$, assume that
 	\beqlb\label{eqn.AsymCriticality}
 	\lim_{n\to\infty} \big\|\boldsymbol{\phi}_n \big\|_{L^1} = \boldsymbol{K}.
 	\eeqlb
 \end{assumption}
 
 Let $\boldsymbol{R}_{n}$ be the resolvent associated to $\boldsymbol{\phi}_{n}$ and $\widetilde{\boldsymbol{N}}_{n} $ be the compensated point process of $\boldsymbol{N}_{n}$ defined as in (\ref{Resolvent}) and (\ref{eqn.CompensatedN}) respectively. 
 By (\ref{SVR}), the density process $\boldsymbol{\varLambda}_{n}$ satisfies the following equation
 \beqlb\label{SVR.n}
 \boldsymbol{\varLambda}_{n}(t)
 =\boldsymbol{H}_n(t)  + \boldsymbol{R}_{n} *d\widetilde{\boldsymbol{N}}_{n}(t)
 \quad \mbox{with} \quad
 \boldsymbol{H}_n(t):= \boldsymbol{\mu}_{n}(t) +    \boldsymbol{R}_{n}* \boldsymbol{\mu}_{n} (t)
 ,\quad t\geq 0 . 
 \eeqlb
 The moment formula of compensated point processes induces that 
 \beqnn 
 \mathbf{E} \Big[\big|\widetilde{\boldsymbol{N}}_{n}(nt)\big|^2 \Big] = \mathbf{E} \big[\boldsymbol{N}_{n}(nt)\big]= \mathbf{E}\big[ \mathcal{I}_{\boldsymbol{\varLambda}_{n}}(nt) \big] = n\cdot \int_0^t \mathbf{E}\big[\boldsymbol{\varLambda}_{n}(ns)\big] ds
 \sim \mathbf{E}\big[\boldsymbol{\varLambda}_{n}(nt)\big] \cdot nt,
 \eeqnn
 as $n\to\infty$ for any $t>0$.
 Thus, a natural scaling in time and space leads us to consider the convergence of the following rescaled processes:
 \beqlb \label{eqn.ScaledProcess}
 \boldsymbol{\varLambda}^{(n)} (t) := \frac{\boldsymbol{\varLambda}_{n}(nt)}{\theta_n},\quad 
 \mathcal{I}_{\boldsymbol{\varLambda}^{(n)}}(t):= \frac{\mathcal{I}_{\boldsymbol{\varLambda}_{n}}(nt)}{n\cdot\theta_n}, \quad
 \boldsymbol{N}^{(n)} (t):= \frac{\boldsymbol{N}_{n}(nt)}{n\cdot\theta_n},
 \quad  
 \widetilde{\boldsymbol{N}}^{(n)} (t) := \frac{\widetilde{\boldsymbol{N}}_{n}(nt)}{\sqrt{n\cdot\theta_n}},
 \quad t\geq 0, 
 \eeqlb
 where $\{ \theta_n \}_{n\geq 1}$  is a sequence of positive scaling parameters such that $\theta_n\to \infty$ as $n\to\infty$.
 
 We now give some sufficient conditions on the parameter $(\boldsymbol{\mu}_{n},\boldsymbol{\phi}_{n})$ for the existence of non-degenerated scaling limits.  
 Let $ \boldsymbol{\mu}^{(n)}$ and $ \boldsymbol{R}^{(n)}$ be the rescaled exogenous intensity and rescaled resolvent given by 
 \beqnn
  \boldsymbol{\mu}^{(n)}(t):= \sqrt{\frac{n}{\theta_n}}\cdot \boldsymbol{\mu}_{n}(nt)
  \quad \mbox{and}\quad
  \boldsymbol{R}^{(n)} (t):= \sqrt{\frac{n}{\theta_n}}\cdot\boldsymbol{R}_{n}(nt),
  \quad t\geq 0. 
 \eeqnn
 In view of (\ref{SVR.n}), the rescaled process $\boldsymbol{\varLambda}^{(n)}$ satisfies 
 \beqlb\label{IntegralLambda}
 \boldsymbol{\varLambda}^{(n)} (t)
 \ar=\ar \boldsymbol{H}^{(n)}(t) 
 +\boldsymbol{R}^{(n)}*d\widetilde{\boldsymbol{N}}^{(n)} (t)
 \quad \mbox{with}\quad
  \boldsymbol{H}^{(n)}(t):= 
 \frac{\boldsymbol{\mu}^{(n)}(t)}{\sqrt{n\cdot \theta_n}} +  \boldsymbol{R}^{(n)}*\boldsymbol{\mu}^{(n)}(t),\quad t\geq 0. \quad\ 
 \eeqlb
 Its long-run behavior is fully determined by the asymptoticsof  $(\boldsymbol{H}^{(n)},\boldsymbol{R}^{(n)})$  or $(\boldsymbol{\mu}^{(n)}, \boldsymbol{R}^{(n)})$. 
 
 We start to study  the asymptotics of $\boldsymbol{R}^{(n)}$ via its Laplace transform. 
 Denote by $\mathcal{L}_{\boldsymbol{\phi}_{n}}$ and $\mathcal{L}_{\boldsymbol{R}_{n}}$ the Laplace transforms of $\boldsymbol{\phi}_{n}$ and $\boldsymbol{R}_{n}$ respectively, i.e.,  
 \beqlb
 \mathcal{L}_{\boldsymbol{\phi}_{n}}(\lambda) := \int_0^\infty e^{-\lambda t} \boldsymbol{\phi}_{n}(t)dt
 \quad \mbox{and} \quad 
 \mathcal{L}_{\boldsymbol{R}_{n}}(\lambda) := \int_0^\infty e^{-\lambda t} \boldsymbol{R}_{n}(t)dt,\quad \lambda >0. 
 \eeqlb 
 Taking the Laplace transform of each term in  (\ref{Resolvent}) gives the following equality
 \beqlb\label{eqn.ConnectionLapPhiR}
 \mathcal{L}_{\boldsymbol{R}_{n}}(\lambda) = \mathcal{L}_{\boldsymbol{\phi}_{n}} (\lambda) + \mathcal{L}_{\boldsymbol{R}_{n}}(\lambda)\cdot \mathcal{L}_{\boldsymbol{\phi}_{n}} (\lambda) 
 =\mathcal{L}_{\boldsymbol{\phi}_{n}} (\lambda) +  \mathcal{L}_{\boldsymbol{\phi}_{n}} (\lambda)\cdot  \mathcal{L}_{\boldsymbol{R}_{n}}(\lambda).
 \eeqlb  
 As $\lambda\to \infty$, the matrix $\mathcal{L}_{\boldsymbol{\phi}_{n}} (\lambda/n)$ decreases to zero and hence its spectral radius continuously decreases to $0$; see Proposition~\ref{Appendix.Prop.LimitMatrix}(1).  
 Let $\lambda_{1,n}^{-}$ denotes the first passage point of its spectral radius into $[0,1]$, i.e.,
 \beqnn
 \lambda_{1,n}^{-} :=  \inf\big\{ \lambda \geq 0:  \rho\big(\mathcal{L}_{\boldsymbol{\phi}_{n}} (\lambda/n)\big) \leq 1 \big\}<\infty. 
 \eeqnn 
 Consider an important auxiliary matrix-valued function ${\boldsymbol{\it\Psi}}^{(n)} $ on $\mathbb{R}_+$ defined by
 \beqlb\label{eqn.varPsin}
 {\boldsymbol{\it\Psi}}^{(n)} (\lambda) := \sqrt{n\cdot \theta_n} \cdot \big(\mathbf{Id} -\mathcal{L}_{\boldsymbol{\phi}_{n}} (\lambda/n) \big) , \quad \lambda \geq 0 .
 \eeqlb 
 When $\lambda> \lambda_{1,n}^{-}$, the matrix $ {\boldsymbol{\it\Psi}}^{(n)} (\lambda)$ is a $M$-matrix\footnote{A matrix $\boldsymbol{A} \in \mathbb{R}^{d\times d}$ is called a \textsl{$M$-matrix} if it is of the form 
 $ \boldsymbol{A}= \lambda \cdot \mathbf{Id}-\boldsymbol{B}$ with $\boldsymbol{B} \in \mathbb{R}^{d\times d}_+$ and $\lambda >\rho (\boldsymbol{B})$.} and hence invertible. 
 By using the change of variables and then (\ref{eqn.ConnectionLapPhiR}) along with (\ref{eqn.varPsin}), the Laplace transform of $\boldsymbol{R}^{(n)}$ has the following representation
 \beqlb\label{eqn.LaplaceIntR} 
  \mathcal{L}_{\boldsymbol{R}^{(n)}}(\lambda )
  :=  \frac{ 
 	\mathcal{L}_{\boldsymbol{R}_n}(\lambda /n)}{ \sqrt{n\cdot \theta_n}}       
  =  \big(  \boldsymbol{\it\Psi}^{(n)}(\lambda) \big)^{-1}  \cdot  \mathcal{L}_{\boldsymbol{\phi}_n} (\lambda/n) 
  =  \mathcal{L}_{\boldsymbol{\phi}_n} (\lambda/n) \cdot   \big(  \boldsymbol{\it\Psi}^{(n)}(\lambda) \big)^{-1} ,
  \quad \lambda > \lambda_{1,n}^{-}.
 \eeqlb 
 Hence the convergence of rescaled resolvent $\boldsymbol{R}^{(n)}$ is fully determined by the asymptotics of $\mathcal{L}_{\boldsymbol{\phi}_n} (\lambda/n)$, which can be exactly established under the next mild condition.  
 
 \begin{condition}\label{Main.Condition}
 	There exists a function $\boldsymbol{\varPhi}\in  C(\mathbb{R}_+;\mathbb{R}^{d\times d})$ such that as $n\to\infty$,
 	\beqlb\label{eqn.Main.Condition} 
 	\boldsymbol{\varPhi}^{(n)}(\lambda):=\sqrt{n\cdot \theta_n} \cdot \big( \boldsymbol{K}-  \mathcal{L}_{\boldsymbol{\phi}_n} (\lambda/n) \big) \to \boldsymbol{\varPhi}(\lambda), 
 	\quad \lambda \geq 0.
 	\eeqlb 
 \end{condition} 
 
 Let $\mathcal{EBF}$ be the space of all \textsl{extended Bernstein functions} on $\mathbb{R}_+$ and $\mathcal{EBF}^{d\times d}$ the space of all functions $\boldsymbol{f} \in C(\mathbb{R}_+;\mathbb{R}^{d\times d})$ with $f_{ij} \in \mathcal{EBF}$; see Appendix~\ref{Appendix--EBF}. 
 By Proposition~\ref{Prop.EBF01}, the function $\boldsymbol{f}$ admits the following representation
 \beqlb 
 \boldsymbol{f}(\lambda)= \boldsymbol{b}^{\boldsymbol{f}} + \boldsymbol{\sigma}^{\boldsymbol{f}}  \cdot \lambda 
 + \int_0^\infty \big(1-e^{-\lambda t} \big) \boldsymbol{\nu}^{\boldsymbol{f}} (dt),\quad \lambda \geq 0,
 \eeqlb
 for some \textsl{L\'evy triplet} $(\boldsymbol{b}^{\boldsymbol{f}}
  ,\boldsymbol{\sigma}^{\boldsymbol{f}} ,\boldsymbol{\nu}^{\boldsymbol{f}} )$ with $\boldsymbol{b}^{\boldsymbol{f}}  \in \mathbb{R}^{d\times d}$, $\boldsymbol{\sigma}^{\boldsymbol{f}}  \in \mathbb{R}_+^{d\times d}$ and $\boldsymbol{\nu}^{\boldsymbol{f}}(dt) \in M_{\rm loc}\big((0,\infty);\mathbb{R}_+^{d\times d} \big) $ satisfying $
 \int_0^\infty (1\wedge t) \boldsymbol{\nu}^{\boldsymbol{f}} (dt) <\infty. $
 By the change of variables,  we can write $\boldsymbol{\varPhi}^{(n)}$ as  
 \beqnn
 \boldsymbol{\varPhi}^{(n)}(\lambda)
 = \sqrt{n\cdot \theta_n} \big( \boldsymbol{K}- \|\boldsymbol{\phi}_n\|_{L^1}\big) +\int_0^\infty \big(1-e^{-\lambda t} \big) \cdot  \sqrt{n^3\cdot \theta_n} \boldsymbol{\phi}_n(nt)dt  \in\mathcal{EBF}^{d\times d}. 
 \eeqnn 
 The next proposition is a direct consequence of Condition~\ref{Main.Condition},  Proposition~\ref{Prop.EBF02} and Corollary~\ref{Corollary.EBF01}. 
 
 \begin{proposition} \label{Prop.RepPhi}
 	The following three claims hold.
 	\begin{enumerate} 
 		\item[(1)] Under Condition~\ref{Main.Condition}, we have $ \boldsymbol{\varPhi} \in  \mathcal{EBF}^{d\times d}$ with L\'evy triplet $(\boldsymbol{b}^{\boldsymbol{\varPhi}},\boldsymbol{\sigma}^{\boldsymbol{\varPhi}},\boldsymbol{\nu}^{\boldsymbol{\varPhi}})$. 
 		
 		\item[(2)] Condition~\ref{Main.Condition} holds if and only if the following limits hold as $n\to\infty$:
 		\beqlb\label{EquaviMainCondition.01}
 		\sqrt{n\cdot \theta_n} \big( \boldsymbol{K}- \|\boldsymbol{\phi}_n\|_{L^1}\big) \to \boldsymbol{b}^{\boldsymbol{\varPhi}}, \quad 
 		\sqrt{n\cdot \theta_n} \int_0^\infty \Big(\frac{t}{n}\wedge 1\Big) \boldsymbol{\phi}_n(t)dt \to \boldsymbol{\sigma}^{\boldsymbol{\varPhi}}+ \int_0^\infty \big(1\wedge t\big) \boldsymbol{\nu}^{\boldsymbol{\varPhi}}(dt), 
 		\eeqlb
 		and for any bounded, continuous function $g $ on $\mathbb{R}_+$ satisfying that $g(x)=O(x)$ as $x\to0+$, 
 		\beqlb\label{EquaviMainCondition.02}
 		\sqrt{n\cdot \theta_n}  \int_0^\infty g\Big(\frac{t}{n}\Big)\boldsymbol{\phi}_n(t)dt \to  \int_0^\infty g(t) \boldsymbol{\nu}^{\boldsymbol{\varPhi}}(dt). 
 		\eeqlb 
 		
 		\item[(3)]  For any $\boldsymbol{K} \in \mathbb{R}_{+,1,\ell}^{d\times d}$ and $\boldsymbol{\varPhi} \in \mathcal{EBF}^{d\times d}$,   
 		there exist a sequence of Hawkes processes with parameters $\{ (\boldsymbol{\mu}_n,\boldsymbol{\phi}_n) \}_{n\geq 1}$ and a sequence of positive scaling parameters $\{\theta_n\}_{n\geq 1}$ satisfying Condition~\ref{Main.Condition}. 
 	\end{enumerate}
 \end{proposition}

 To get our scaling limit theorem, it remains to  study the asymptotics of the matrix inverse $\big( \boldsymbol{\varPsi}^{(n)} \big)^{-1}$. 
 Since $\mathbf{Id} -\mathcal{L}_{\boldsymbol{\phi}_{n}} (\lambda/n)$ may be non-invertible as $n\to\infty$, its inverse is quite elusive.   
 Hence some limitations on  $\boldsymbol{\varPhi}$ are necessary to overcome this difficulty.  
 Let $\mathtt{I}:= \{ 1,\cdots, \ell \} $ and $\mathtt{J}:= \{ \ell +1 ,\cdots, d \}$.  
 For a vector $\boldsymbol{x}= (x_i)_{i\in\tt{D}}$   and a $d\times d$-matrix $\boldsymbol{A}=(A_{ij})_{i,j\in\mathtt{D}}$, we write
 \beqnn
 \boldsymbol{x}=\begin{pmatrix}\boldsymbol{x}_{\mathtt{I}} &
 	\boldsymbol{x}_{\mathtt{J}} \end{pmatrix}
 \quad\mbox{and}\quad 
 \boldsymbol{A}=\begin{pmatrix}
 	\boldsymbol{A}_{\mathtt{I}\mathtt{I}} & \boldsymbol{A}_{\mathtt{I}\mathtt{J}}\\
 	\boldsymbol{A}_{\mathtt{J}\mathtt{I}} & \boldsymbol{A}_{\mathtt{J}\mathtt{J}}
 \end{pmatrix}  
 \eeqnn
 with $\boldsymbol{x}_{\mathtt{I}}:=  (x_i)_{i\in\mathtt{I}}$, $\boldsymbol{A}_{\mathtt{I}\mathtt{J}}:=(A_{ij})_{i\in\mathtt{I},j\in \mathtt{J}}$ and other blocks defined in the same way. 
 By the real Schur decomposition; see Proposition~\ref{Appendix.Prop.RSD}, there always exist (not unique) an orthogonal  matrix $\boldsymbol{Q} $ and an upper triangular matrix $\boldsymbol{U}$ with all diagonal entries of $\boldsymbol{U}_{\mathtt{I}\mathtt{I}}$ equal to $1$ such that $\boldsymbol{K} = \boldsymbol{Q} \boldsymbol{U} \boldsymbol{Q}^{\rm T}$. 
 
 \begin{definition}\label{Def.Admissible}
 The limit function $\boldsymbol{\varPhi}$ is called {\rm $\boldsymbol{K}$-admissible} if there exist two auxiliary matrices $\boldsymbol{Q},\boldsymbol{U} \in \mathbb{R}^{d\times d}$ and an auxiliary constant $\lambda^+\geq 0$ such that 
 \begin{enumerate}
 \item[(1)] $\boldsymbol{Q}$ is orthogonal, $\boldsymbol{U}_{\mathtt{I}\mathtt{I}} = \mathbf{Id}\in \mathbb{R}^{\ell\times \ell}$ and $\boldsymbol{U}_{\mathtt{J}\mathtt{I}}=\mathbf{0} \in \mathbb{R}^{(d-\ell)\times \ell}$ such that $\boldsymbol{K} = \boldsymbol{Q} \boldsymbol{U} \boldsymbol{Q}^{\rm T}$;
 		
 \item[(2)]	$\displaystyle \boldsymbol{\varphi}(\lambda):=\big( 	\boldsymbol{Q}^{\rm T}  \boldsymbol{\varPhi}(\lambda )\boldsymbol{Q} \big)_{\mathtt{I}\mathtt{I}}	+\boldsymbol{U}_{\mathtt{I}\mathtt{J}} \big(  \mathbf{Id}- \boldsymbol{U}_{\mathtt{J}\mathtt{J}}\big)^{-1}\big( \boldsymbol{Q}^{\rm T}  \boldsymbol{\varPhi}(\lambda ) \boldsymbol{Q} \big)_{\mathtt{J}\mathtt{I}} \in \mathbb{R}^{\ell \times \ell} $ is invertible for any $\lambda \geq \lambda^+ $. 		
 \end{enumerate}  
 	
 \end{definition}
  
 Since $\boldsymbol{K} \in \mathbb{R}_{+,1,\ell}^{d\times d}$, all eigenvalues of  $\boldsymbol{U}_{\mathtt{J}\mathtt{J}}  $ have real parts less than one and hence  the matrix $\mathbf{Id}-\boldsymbol{U}_{\mathtt{J}\mathtt{J}} $ is invertible.  
 The existence of auxiliary constant $\lambda^+ $ prevents the integrated rescaled resolvent $\mathcal{I}_{\boldsymbol{R}^{(n)}}$ from growing super-exponentially  when $n$ is large. 
 More precisely, if $\displaystyle \boldsymbol{\varphi}(\lambda)$ is non-invertible for any $\lambda \geq 0$, it would be impossible to find a constant $\beta>0$ such that 
 \beqnn
 \limsup_{t\to\infty} \sup_{n\geq 1} \sup_{i,j\in\mathcal{D}} e^{-\beta t} \int_0^t R_{ij}^{(n)}(s)ds <\infty. 
 \eeqnn
 
 \begin{remark}
 In the uni-variate case $(d=1)$, since $\varPhi$ is non-decreasing, it is admissible if and only if the first passage time of $ \varPhi$ into $(0,\infty)$ is finite, i.e., $\lambda_0^+ := \inf\{\lambda \geq 0: \varPhi(\lambda)>0\}<\infty$, which holds if and only if $\sigma^{\varPhi} >0$ or $b^{\varPhi}+\nu^{\varPhi}  (\mathbb{R}_+) \in (0,\infty]$. 
 Moreover, Condition~\ref{Main.Condition} with $\lambda_0^+ <\infty$ is the necessary and sufficient conditions for the convergence of rescaled uni-variate Hawkes processes.
 \end{remark}

 \begin{lemma} \label{Lemma.Pi}
 If $\boldsymbol{\varPhi}$ is  $\boldsymbol{K}$-admissible, then the following hold.
 \begin{enumerate}
 \item[(1)] There exists a unique measure $\boldsymbol{\varPi}(dt) \in  M_{\rm loc}(\mathbb{R}_+;\mathbb{R}_+^{d\times d}) $  such that for any auxiliary matrices $\boldsymbol{Q},\boldsymbol{U} \in \mathbb{R}^{d\times d}$ and auxiliary constant $\lambda^+ \geq 0$, 
 \beqlb \label{eqn.LimitResolvent}
 \mathcal{L}_{\boldsymbol{\varPi}}(\lambda):= \int_0^\infty e^{-\lambda t} \boldsymbol{\varPi} (dt)
 =  
 \boldsymbol{Q}
 \begin{pmatrix}
 \big(\boldsymbol{\varphi}(\lambda)\big)^{-1} & 
 \big(\boldsymbol{\varphi}(\lambda)\big)^{-1} \boldsymbol{U}_{\mathtt{I}\mathtt{J}}\big(\mathbf{Id}-\boldsymbol{U}_{\mathtt{J}\mathtt{J}} \big)^{-1} \vspace{7pt} \\
 \mathbf{0} & \mathbf{0}
 \end{pmatrix}
 \boldsymbol{Q}^{\rm T},
 \quad \lambda \geq \lambda^+ .
 \eeqlb
  	
 \item[(2)] Under Condition~\ref{Main.Condition}, the measure with density $\boldsymbol{R}^{(n)}$ converges vaguely to  $\boldsymbol{\varPi}(dt)$ as $n\to\infty$, i.e.,
 \beqnn
 \int_{[0,T]} f(t)\boldsymbol{R}^{(n)}(t)dt \to \int_{[0,T]} f(t) \boldsymbol{\varPi}(dt), 
 \quad T\geq 0,\, f \in C(\mathbb{R}_+;\mathbb{R}).
 \eeqnn  
 \end{enumerate} 
 \end{lemma}
 
 The proof of this lemma is provided in Section~\ref{Sec.AsymResolvent}.
 Claim (1) shows that every $\boldsymbol{K}$-admissible function $\boldsymbol{\varPhi}$ uniquely determines a $\sigma$-finite measure $\boldsymbol{\varPi}(dt)$ on $\mathbb{R}_+$. 
 Drawing from the potential theory; see Chapter~II in \cite{Bertoin1996}, we may refer to  $\boldsymbol{\varPi}(dt)$  as the \textsl{$(\boldsymbol{K},\boldsymbol{\varPhi})$-potential measure}.  
 Its distribution function is denoted as $\boldsymbol{\underline\varPi}$, i.e.,
 \beqnn
 \boldsymbol{\underline\varPi}(t):= \boldsymbol{\varPi}\big([0,t]\big),\quad t\geq 0.
 \eeqnn  
 If $\boldsymbol{\varphi}(\infty)=\infty$, we have $  \boldsymbol{\underline\varPi} (0)= \mathbf{0}$; otherwise, we may see that $  \boldsymbol{\underline\varPi} (0) \neq \mathbf{0}$.   
 Two direct consequences of Lemma~\ref{Lemma.Pi}(2) are
 \beqlb\label{Eqn.BoundConvIRn}
 \sup_{n\geq 1}\mathcal{I}_{\boldsymbol{R}^{(n)}}(s)  <\infty 
 \quad \mbox{and}\quad 
 \lim_{n\to\infty}\mathcal{I}_{\boldsymbol{R}^{(n)}}(t)  =  \boldsymbol{\underline\varPi}(t) ,
 \eeqlb
 for all $s\geq 0$ and all continuous points $t$ of $\boldsymbol{\underline\varPi}$. 
 This limit holds uniformly on compacts if $\boldsymbol{\underline\varPi}$ is continuous on $\mathbb{R}$. 
 Without further mention, in the sequel we always assume the next assumption. 
 
 \begin{assumption}\label{MainAssumption02}
 The limit function $\boldsymbol{\varPhi}$ is $\boldsymbol{K}$-admissible and $\boldsymbol{\underline\varPi}(0)$ is a diagonal matrix. 
 \end{assumption}
 
 The diagonal assumption on  $\boldsymbol{\underline\varPi}(0)$ is not always necessary in the following results; see Remark~\ref{Remark.2.23}. 
 Based on the preceding asymptotic analyses and notation, we are ready to formulate our first main theorem that concerns the weak convergence of the sequence of rescaled processes $\big\{(\mathcal{I}_{\boldsymbol{\varLambda}^{(n)}}, \boldsymbol{N}^{(n)}, \widetilde{\boldsymbol{N}}^{(n)})\big\}_{n\geq 1}$. 
 Its proof is provided in Section~\ref{ProofMainThm}.
 
 \begin{theorem} \label{MainThm.S-WeakConvergence}
 Under Condition~\ref{Main.Condition}, if $\mathcal{I}_{\boldsymbol{H}^{(n)}}\to \boldsymbol{\varUpsilon} \in D(\mathbb{R};\mathbb{R}_+^d)$ along all continuity points of $\boldsymbol{\varUpsilon}$, then the following hold.
 \begin{enumerate}
 	\item[(1)] There exist a non-decreasing process $\boldsymbol{\varXi} \in D(\mathbb{R}_+;\mathbb{R}_+^d) $ and  a square-integrable martingale $\boldsymbol{M} \in D(\mathbb{R}_+;\mathbb{R}^d)$ such that  as $n\to\infty$, 
    \beqlb\label{MainWeakConv}
 	\big(\mathcal{I}_{\boldsymbol{\varLambda}^{(n)}}, \boldsymbol{N}^{(n)}, \widetilde{\boldsymbol{N}}^{(n)}\big)  
 	\overset{\rm d}\to 
 	\big(\boldsymbol{\varXi},\boldsymbol{\varXi},\boldsymbol{M}\big),
 	\eeqlb
 	in $D_S(\mathbb{R}_+;\mathbb{R}_+^d\times\mathbb{R}_+^d\times\mathbb{R}^d)$. Moreover, $ (\boldsymbol{\varXi},\boldsymbol{M})$ is stochastically continuous if $\boldsymbol{\varUpsilon}$ is continuous on $\mathbb{R}_+$.
 	
 	\item[(2)]  For a countable set $Q_0 \subset\mathbb{R}_+$, the following two claims hold.
 	\begin{enumerate}
 		\item[(2.i)] The weak convergence (\ref{MainWeakConv}) holds in the sense of finite-dimensional distributions along $\mathbb{R}_+\setminus Q_0$.
 		
 		\item[(2.ii)] For any $\boldsymbol{f},\boldsymbol{g} \in C^1(\mathbb{R}_+;\mathbb{C}_-^d)$  and  $\boldsymbol{h}\in C^1(\mathbb{R}_+;\mathtt{i}\mathbb{R}^d)$, we have 
 		\beqnn
 		\big( \boldsymbol{f} * \boldsymbol{\varLambda}^{(n)} ,\boldsymbol{g}* d\boldsymbol{N}^{(n)} , \boldsymbol{h}* d\widetilde{\boldsymbol{N}}^{(n)}   \big)
 		\overset{\rm f.f.d.}\longrightarrow 
 		\big( \boldsymbol{f} * d\boldsymbol{\varXi},\boldsymbol{g} * d\boldsymbol{\varXi}, \boldsymbol{h}* d\boldsymbol{M}  \big),
 		\eeqnn
 		along $ \mathbb{R}_+\setminus Q_0$	as $n\to\infty$. 
 		Moreover, it holds along $ \mathbb{R}_+$ if $\boldsymbol{f}(0)=\boldsymbol{g}(0)=\boldsymbol{h}(0) = \mathbf{0}$. 
 		
 	\end{enumerate} 
 \end{enumerate}

 \end{theorem} 
 
 \begin{remark}\label{Remark.01}
 The assumption $\mathcal{I}_{\boldsymbol{H}^{(n)}}\to \boldsymbol{\varUpsilon} \in D(\mathbb{R};\mathbb{R}_+^d)$  is natural, e.g.,
 \begin{enumerate}
 \item[(1)]  For $j\in\mathtt{D}$, $\boldsymbol{a}\in \mathbb{R}_+^d$ and $\boldsymbol{A} \in \mathbb{R}_+^{d\times d}$, let
 $\boldsymbol{\mu}_n =  \sqrt{n\theta_n}\cdot \boldsymbol{A}\cdot (\phi_{n,ij})_{i\in\mathtt{D}} + \sqrt{\theta_n/n} \cdot \boldsymbol{a}$ be the combination of total impacts of all type-$j$ events prior to time $0$  and the constant external excitation respectively. 
 For all continuous points $t$ of $\underline{\boldsymbol{\varPi}}$, we have
 \beqnn
 \mathcal{I}_{\boldsymbol{H}^{(n)}}(t)= \boldsymbol{A}\cdot 
 \big(\mathcal{I}_{R^{(n)}_{ij}}(t)\big)_{i\in\mathtt{D}}  + \frac{1}{\sqrt{n\theta_n}} \cdot  \boldsymbol{a}\cdot t  +  \mathcal{I}_{\boldsymbol{R}^{(n)}}* \boldsymbol{a}  (t)\to  \boldsymbol{A}\cdot \big(\underline{\varPi}_{ij}(t)\big)_{i\in\mathtt{D}} + \underline{\boldsymbol{\varPi}}*\boldsymbol{a}  (t).
 \eeqnn
 		
 \item[(2)]  If $\mathcal{I}_{\boldsymbol{\mu}^{(n)}} \to \boldsymbol{\varGamma} \in D(\mathbb{R};\mathbb{R}_+^d)$ along all continuous points of $ \boldsymbol{\varGamma}$, then  along all continuous points of $(\boldsymbol{\varPi}, \boldsymbol{\varGamma})$,
 \beqnn
 \mathcal{I}_{\boldsymbol{H}^{(n)}}  = \frac{\mathcal{I}_{\boldsymbol{\mu}^{(n)}}}{\sqrt{n\theta_n}} +  \boldsymbol{R}^{(n)} * \mathcal{I}_{\boldsymbol{\mu}^{(n)}} \to \boldsymbol{\varPi}*\boldsymbol{\varGamma} \in D(\mathbb{R};\mathbb{R}_+^d).
 \eeqnn  
 \end{enumerate}
 \end{remark}
 
%

 Since the $S$-topology is much weaker than the $J_1$-topology,  the $S$-convergence usually has a few of disadvantages and fails to provide plenty of explicit characteristics of limit processes. 
 For instance, the continuity of pre-limit processes $\{\mathcal{I}_{\boldsymbol{\varLambda}^{(n)}}\}_{n\geq 1}$ may not be inherited by the limit process $\boldsymbol{\varXi}$, e.g., it jumps when $\boldsymbol{\varUpsilon}$ jumps. 
 Although the martingale $\widetilde{\boldsymbol{N}}^{(n)}$ has predictable quadratic co-variation ${\rm diag}(\mathcal{I}_{\boldsymbol{\varLambda}^{(n)}})$ and quadratic co-variation ${\rm diag}(\boldsymbol{N}^{(n)})$, one cannot claim that the martingale $\boldsymbol{M}$ has (predictable) quadratic co-variation ${\rm diag}(\boldsymbol{\varXi})$ unless  it is continuous;  see the last example in \cite[p.371]{MeyerZheng1984}.  
 Nevertheless, it is gratifying to note that  $(\boldsymbol{\varXi},\boldsymbol{M})$  enjoys some properties of exponential martingales. 
 For instance,  a direct consequence of Theorem~\ref{MainThm.S-WeakConvergence} and the continuous mapping theorem (see \cite[Theorem~2.7]{Billingsley1999}) is  
 \beqnn
 \exp\Big\{ \mathtt{i} \boldsymbol{z}\cdot\widetilde{\boldsymbol{N}}^{(n)} -   n\theta_n \Big(\exp\Big\{\frac{\mathtt{i}\boldsymbol{z}}{\sqrt{n\theta_n}}\Big\}-1- \frac{\mathtt{i}\boldsymbol{z}}{\sqrt{n\theta_n}}\Big) \cdot \mathcal{I}_{\boldsymbol{\varLambda}^{(n)}} \Big\} 
 \overset{\rm d}\to \exp\Big\{\mathtt{i} \boldsymbol{z}\cdot\boldsymbol{M} +\frac{1}{2}\cdot  \boldsymbol{z}^2\cdot \boldsymbol{\varXi}  \Big\},\quad \boldsymbol{z}  \in \mathbb{R}^d,
 \eeqnn
 in $D_S(\mathbb{R}_+;\mathbb{C})$ as $n\to\infty$.  
 The standard stopping time argument tells that the pre-limit processes are local martingales, which along with Proposition~\ref{Appendix.Prop.S07} induces that so is the limit process. 
 However, this is still not enough to assert that the martingale $\boldsymbol{M}$ has (predictable) quadratic co-variation ${\rm diag}(\boldsymbol{\varXi})$ by using Exercise~3.14 in \cite[p.153]{RevuzYor2005}  unless either $\boldsymbol{M}$ or $\boldsymbol{\varPi}$ is continuous; see again the last example in \cite[p.371]{MeyerZheng1984}. 
 The next remark provides a sufficient condition but very difficult to be verified, under which  ${\rm diag}(\boldsymbol{\varXi})$ is the quadratic co-variation of $\boldsymbol{M}$. 
 
 \begin{remark}
 For each $i,j\in\mathtt{D}$, if there are no oscillations of $\widetilde{N}_i^{(n)}$ preceding oscillations of $\widetilde{N}_j^{(n)}$, by Theorem~1 in \cite{Jakubowski1996} and Theorem~\ref{MainThm.S-WeakConvergence}, the stochastic equation
 \beqnn
 \mathbf{1}_{\{i=j\}}\cdot N_i^{(n)}(t) = \widetilde{N}_i^{(n)}(t)\widetilde{N}_j^{(n)}(t)-\int_0^t \widetilde{N}_i^{(n)}(s-)d\widetilde{N}_j^{(n)}(s)-\int_0^t \widetilde{N}_j^{(n)}(s-)d\widetilde{N}_i^{(n)}(s),\quad t\geq 0,
 \eeqnn
 converges weakly in the sense of $S$-topology to
 \beqnn
 \mathbf{1}_{\{i=j\}}\cdot \varXi_i(t)= M_i(t)M_j(t)-\int_0^t M_i(s-)dM_j(s)-\int_0^t M_j(s-)dM_i(s), \quad t\geq 0,
 \eeqnn  
 as $n\to\infty$. 
 This tells that $\boldsymbol{M}$ has quadratic co-variation ${\rm diag}(\boldsymbol{\varXi})$; see Definition~4.45 in \cite[p.51]{JacodShiryaev2003}.  
 \end{remark}
 
 In the sequel of this section, we accurately characterize the limit process $(\boldsymbol{\varXi},\boldsymbol{M})$ in different ways and perspectives.
 The next theorem provides an explicit exponential-affine representation of their Fourier-Laplace transform in term of the unique solution to a Riccati-Volterra equation with measure kernel. The proof can be found in Section~\ref{ProofMainThm}.

 \begin{theorem}\label{MainThm.FourierLaplaceFunctional}
 For any row vector functions $\boldsymbol{f}\in C^1(\mathbb{R}_+;\mathbb{C}_-^d)$ and $\boldsymbol{h}\in C^1(\mathbb{R}_+;\mathtt{i}\mathbb{R}^d)$, we have
 \begin{enumerate}
  \item[(1)] There exists a unique global solution $ \boldsymbol{\mathcal{V}}\in D(\mathbb{R}_+;\mathbb{C}^d_-) $ to the  Riccati-Volterra equation 
  \beqlb \label{Eqn.AffineVolterra}
  \boldsymbol{\mathcal{V}}  =   \boldsymbol{\mathcal{W}}* \boldsymbol{\varPi}  
  \quad \mbox{with} \quad
  \boldsymbol{\mathcal{W}}:=  \boldsymbol{f}  + \frac{1}{2} \cdot (\boldsymbol{\mathcal{V}}+\boldsymbol{h})^2.
  \eeqlb
 Moreover, the solution is continuous if either $\boldsymbol{f}(0)=\boldsymbol{h}(0) = \mathbf{0}$ or $\boldsymbol{\underline\varPi}$ is continuous on $\mathbb{R}_+$. 
  		
 \item[(2)] The Fourier-Laplace functional of  $(\boldsymbol{\varXi} ,\boldsymbol{M})$ admits the representation
 \beqlb \label{eqn.FL}
 \mathbf{E}\Big[\exp \big\{ \boldsymbol{f} * d\boldsymbol{\varXi} (T) + \boldsymbol{h}* d\boldsymbol{M} (T)   \big\} \Big]
 = \exp \big\{ \boldsymbol{\mathcal{W}} * d\boldsymbol{\varUpsilon} (T) \big\},
 \quad T\geq 0.
 \eeqlb 
 \end{enumerate}    
 \end{theorem}
 
 Since  methods and tools developed in the theory of Volterra equations are out of work, in contrast to the abundant literature on the classic (non-)linear Volterra equations, only a few results about linear Volterra equations with measure kernel can be found; see \cite{GripenbergLondenStaffans1990} and references therein. 
 To the best of our knowledge, it is the first time to consider the existence and uniqueness of solutions of Riccati-Volterra integral equations with measure kernels. 
 
 In the setting of Remark~\ref{Remark.01}(1) with $\boldsymbol{A}=x\cdot \mathbf{Id}$ and $\boldsymbol{a}\in\mathbb{R}_+^d$, the representation (\ref{eqn.FL}) turns to be 
 \beqnn 
 \mathbf{E}\Big[\exp \big\{ \boldsymbol{f} * d\boldsymbol{\varXi} (T) + \boldsymbol{h}* d\boldsymbol{M} (T)   \big\} \Big]
 = \exp \big\{ x \cdot  \mathcal{V}_j(T) + \boldsymbol{\mathcal{V}} *\boldsymbol{a}(T)\big\},
 \quad T\geq 0. 
 \eeqnn
 In comparison with (2.2) in \cite{DuffieFilipovicSchachermayer2003}, this equality tells that, in a certain sense, the differentiation process $d\boldsymbol{\varXi}$ enjoys the affine property analogous to that of Markov affine processes; see \cite{DuffieFilipovicSchachermayer2003}. 
 Recently, stochastic processes with Fourier-Laplace functional being of the form (\ref{eqn.FL}) are referred to as \textsl{affine Volterra processes} in \cite{Jaber2021,JaberLarssonPulido2019}. 
 We should highlight that our limit process $(\boldsymbol{\varXi},\boldsymbol{M})$ exceeds the existing framework of affine Volterra processes because of the measure kernel $\boldsymbol{\varPi}(dt)$ in  (\ref{Eqn.AffineVolterra}).

 Different to Markov models, it is not enough to clarify the evolution dynamic of $(\boldsymbol{\varXi},\boldsymbol{M})$ through its Fourier-Laplace functional. 
 As our third main result, the next theorem provides another deep view on the limit process via a representation in term of stochastic Volterra equation. The proof is given in Section~\ref{ProofMainThm}.

 \begin{theorem}\label{MainThm.SVERepresentation}
 The limit process $(\boldsymbol{\varXi},\boldsymbol{M})$ satisfies the Volterra integration equation 
 \beqlb\label{Eqn.HawkesVoletrra}
 \boldsymbol{\varXi}(t)= \boldsymbol{\varUpsilon}(t)+ \boldsymbol{\varPi}*\boldsymbol{M}(t)
 ,\quad t\geq 0.
 \eeqlb
 \end{theorem}
 
 With the continuity of the limit process, in the next theorem we establish a convergence result stronger than that in Theorem~\ref{MainThm.S-WeakConvergence} and also provide more explicit characterizations of the limit processes. 
 In the sequel, we denote by  $\boldsymbol{B}:=(B_i)_{i\in\mathtt{D}}^{\rm T}$ a standard $d$-dimensional Brownian motion and define 
 \beqnn
 \boldsymbol{B}\circ \boldsymbol{f}(t):= \big( B_i(f_i(t)) \big)_{i\in\mathtt{D}}^{\rm T},\quad t\geq 0,
 \eeqnn
 for a $d$-dimensional non-negative function $\boldsymbol{f}$ on $\mathbb{R}_+$. 
    
 \begin{theorem}\label{MainThm.J1WeakConvergence}
 If either $\boldsymbol{\varXi} $ or $\boldsymbol{M}$ is continuous, then the following hold.
 \begin{enumerate}
  \item[(1)] The weak convergence in (\ref{MainWeakConv}) also holds in $D_{J_1}(\mathbb{R}_+;\mathbb{R}_+^d\times\mathbb{R}_+^d\times\mathbb{R}^d)$.
  		
  \item[(2)] The limit process $(\boldsymbol{\varXi},\boldsymbol{M}) \in C(\mathbb{R}_+;\mathbb{R}_+^d\times\mathbb{R}^d)$ is the unique continuous weak solution to 
  \beqlb  \label{Eqn.HawkesVoletrra.02}
  \left\{ \begin{split}
   \boldsymbol{\varXi}(t) &= \boldsymbol{\varUpsilon}(t) + \boldsymbol{\varPi}*\boldsymbol{M}(t),\cr
  \boldsymbol{M}(t) &=  \boldsymbol{B} \circ \boldsymbol{\varXi}(t) .
  \end{split} 
  \right.
  \eeqlb
 \end{enumerate}
 \end{theorem}
  
 By Theorem~4.52 in \cite[p.55]{JacodShiryaev2003}, the continuity can be mutually inherited between $ \boldsymbol{\varXi}$ and $ \boldsymbol{M}$.   
 In particular, the continuity of $\boldsymbol{\varUpsilon} $ and $\boldsymbol{\varPi}(dt)$ also can be passed on to $\boldsymbol{\varXi}$; see Proposition~\ref{Prop.A1}(1).  
  
 \begin{corollary}\label{Corollary.ContinuityXiM}
  If $\boldsymbol{\varUpsilon} \in C(\mathbb{R}_+;\mathbb{R}_+^d)$ and $\boldsymbol{\varPi}(dt)$ has density $\boldsymbol{\pi}$, then $(\boldsymbol{\varXi},\boldsymbol{M}) \in C(\mathbb{R}_+;\mathbb{R}_+^d\times \mathbb{R}^d)$ and all claims in Theorem~\ref{MainThm.FourierLaplaceFunctional}  hold for any  $\boldsymbol{f}\in L^\infty_{\rm loc}(\mathbb{R}_+;\mathbb{C}_-^d)$ and $\boldsymbol{h}\in L^\infty_{\rm loc}(\mathbb{R}_+;\mathtt{i}\mathbb{R}^d)$.  
 \end{corollary}
 
 In the next corollary, we show that the path regularity of limit process $\boldsymbol{\varXi}$ can be improved when the density function $\boldsymbol{\pi}$ is square-integrable on compacts.  
 Let $L(dt)$ be the Lebesgue measure on $\mathbb{R}$. 
  
 \begin{corollary}\label{Corollary.2001}
 Assume that  $\boldsymbol{\pi}  \in L^2_{\rm loc}(\mathbb{R}_+;\mathbb{R}^{d\times d}_+)$ and $\boldsymbol{\varUpsilon}$ is differentiable on $\mathbb{R}$ with derivative $\boldsymbol{\varUpsilon}'$, then $\boldsymbol{\varXi} $ is differentiable and has a predictable derivative $\boldsymbol{\xi}$ that is the unique non-negative weak solution to 
 \beqlb  \label{eqn.SVE04}
 \left\{ \begin{split}
 	\boldsymbol{\xi} (t) &=  \boldsymbol{\varUpsilon}'(t)
 	+   \boldsymbol{\pi} *d\boldsymbol{M}(t)  , \quad  \mathbf{P} \otimes L(dt)\mbox{-a.s.} ,  \cr 
 	\boldsymbol{M}(t) &=	\int_0^t {\rm diag } \big(\sqrt{\boldsymbol{\xi} (s)}\big)d\boldsymbol{B}(s),\quad t\geq 0.
 \end{split} 
 \right.
 \eeqlb
 \end{corollary}
 
  
 A sufficient condition for  $\boldsymbol{\pi}  \in L^2_{\rm loc}(\mathbb{R}_+;\mathbb{R}^{d\times d}_+)$ is that the sequence $\{  \boldsymbol{R}^{(n)} \}_{n\geq 1}$ is uniformly square-integrable on compacts.
 The weak uniqueness of solutions to  (\ref{eqn.SVE04}) is defined as the {\rm pseudo-paths} of any two non-negative weak solutions $ \boldsymbol{\xi}$ and $ \boldsymbol{\xi}'$ having the same distribution law, i.e., the two random image measures of Lebesgue measure under the two mappings $t\mapsto (t,\boldsymbol{\xi})$ and $t\mapsto (t,\boldsymbol{\xi}')$ are equal in distribution. 
 We refer to  \cite[Chapter IV]{DellacherieMeyer1975} and \cite{MeyerZheng1984} for more details about pseudo-paths. 
 

 \begin{corollary}\label{Corollary.2002}
 Under assumptions in Corollary~\ref{Corollary.2001},  we have the following.
 
 \begin{enumerate}
 	\item[(1)]   If $\boldsymbol{\pi}$ is continuous on $\mathbb{R}_+$ and differentiable on $(0,\infty)$ with derivative $\boldsymbol{\pi}'$, then \eqref{eqn.SVE04} is  equivalent to
 	 \beqlb\label{eqn.SVE03} 
 		\boldsymbol{\xi}(t)
 		=  \boldsymbol{\varUpsilon}' (t) +   \boldsymbol{\pi}'*	\boldsymbol{M}(t)  +  \boldsymbol{\pi}(0) \cdot 	\boldsymbol{M}(t)  , \quad  \mathbf{P} \otimes L(dt)\mbox{-a.s.}
 	\eeqlb 
 	
 	\item[(2)] Additionally, if  $\boldsymbol{\pi}'$ is continuous on $\mathbb{R}_+$ and differentiable on $(0,\infty)$ with derivative $\boldsymbol{\pi}''$,  then  \eqref{eqn.SVE04} is equivalent to 
 	\beqlb\label{eqn.SVE03001} 
 		\boldsymbol{\xi}(t)
 		=  \boldsymbol{\varUpsilon}' (t) 
 		+ \int_0^t\big(\boldsymbol{\pi}'(0)\cdot \boldsymbol{M}(s) + \boldsymbol{\pi}''*	\boldsymbol{M}(s) \big)\, ds  +  \boldsymbol{\pi}(0) \cdot 	\boldsymbol{M}(t)  , \quad  \mathbf{P} \otimes L(dt)\mbox{-a.s.}  
 	\eeqlb 
 	and pathwise uniqueness holds. 
 \end{enumerate}

 \end{corollary}
 
 The preceding two equivalent equations can be obtained immediately by applying Fubini's theorem and the stochastic Fubini's theorem to \eqref{eqn.SVE04}. 
 The pathwise uniqueness of solutions to \eqref{eqn.SVE03001} can be proved similarly as in the proof of Theorem~5.3 in \cite{PromelScheffels2023}. The detailed proof of this corollary is omitted. 
 Obviously, the derivative process $\boldsymbol{\xi}$ has the same  regularity of $\boldsymbol{\varUpsilon}'$ in case (1) and is a  semimartingale in case (2). 

 The limit process $\boldsymbol{\varXi}$ is mainly determined by the $(\boldsymbol{K},\boldsymbol{\varPhi})$-potential measure $\boldsymbol{\varPi}(dt)$.
 In particular, in the setting of Remark~\ref{Remark.01}(1) with $\boldsymbol{a}=\mathbf{0}$ and $\boldsymbol{A}=\mathbf{Id}$, we have
 \beqnn  
 \boldsymbol{\varXi}(t)= \boldsymbol{\underline{\varPi}}_{\cdot j}(t)+ \boldsymbol{\varPi}*\boldsymbol{M}(t) 
 \quad\mbox{and hence}\quad
 \mathbf{E}\big[\boldsymbol{\varXi}(t)\big] = \boldsymbol{\underline{\varPi}}_{\cdot j}(t),\quad j\in\mathtt{D}. 
 \eeqnn
 Without external excitation, Hawkes systems are self-closed  and tend to be either freezing or boiling, e.g. $\boldsymbol{\underline{\varPi}}_{\cdot j} (\infty) \in \{ 0,\infty \}$; see \cite[Section~2]{HorstXu2022} for detailed explanations in the Markovian case. 
 Let $\lambda_{\boldsymbol{\varPi}}$ denotes the  exponential growth rate of $\boldsymbol{\varPi}$ defined by
 \beqnn
 \lambda_{\boldsymbol{\varPi}}:= \inf \big\{ \lambda \geq 0: \mathcal{L}_{\boldsymbol{\varPi}}(\lambda)<\infty  \big\} \in [0,\lambda^+ ]. 
 \eeqnn
 When $\lambda_{\boldsymbol{\varPi}}>0$, the function $\underline{\boldsymbol{\varPi}}$  has an exponential growth with rate $\lambda_{\boldsymbol{\varPi}}$ and hence the growth rate of $\boldsymbol{\varXi}$ increases exponentially  to infinity with a positive probability. 
 In contrast, when $\lambda_{\boldsymbol{\varPi}}=0$ and $\underline{\boldsymbol{\varPi}} (\infty) <\infty$,  
 the growth rate of $\boldsymbol{\varXi}$ decreases to $0$ quickly as time goes. 
 Finally, when $\lambda_{\boldsymbol{\varPi}}=0$ and $\underline{\boldsymbol{\varPi}}(\infty) =\infty$, the growth rate of $\boldsymbol{\varXi}$  changes very slowly. 
 These give arise to a criticality criterion for the limit process $\boldsymbol{\varXi}$, i.e., it is said to be \textsl{subcritical} if  $ \underline{\boldsymbol{\varPi}}(\infty)<\infty$, \textsl{critical} if $\lambda_{\boldsymbol{\varPi}}=0$ and $\underline{\boldsymbol{\varPi}}(\infty)=\infty$, or \textsl{supercritical} if $\lambda_{\boldsymbol{\varPi}}>0$.
 
 \begin{remark}
 In the uni-variate case ($d=1$),  the criticality criterion  turns to be much clear and concise, i.e.,  $\varXi$ is \textsl{subcritical} if $b^\varPhi>0$, \textsl{critical} if $b^\varPhi=0$ or \textsl{supercritical} if $b^\varPhi<0$.    
 \end{remark}

 The next theorem provides another stochastic Volterra equation equivalent to (\ref{Eqn.HawkesVoletrra}). 
 This alternate representation separates the impact of parameter $\boldsymbol{b}^{\boldsymbol{\varPhi}}$ on the dynamics of $\boldsymbol{\varXi}$ from that of  $(\boldsymbol{\sigma}^{\boldsymbol{\varPhi}},\boldsymbol{\nu}^{\boldsymbol{\varPhi}})$.  
 Define
 \beqnn
  \boldsymbol{\varPhi}_0(\lambda):=  \boldsymbol{\varPhi}(\lambda)-  \boldsymbol{\varPhi}(0) =\boldsymbol{\varPhi}(\lambda)- \boldsymbol{b}^{\boldsymbol{\varPhi}}
  \quad\mbox{and}\quad 
  \boldsymbol{\varphi}_0(\lambda):=  \boldsymbol{\varphi}(\lambda)-   \boldsymbol{\varphi}(0) \in \mathbb{R}^{\ell\times\ell},\quad \lambda \geq 0.
 \eeqnn 
 It is easy to identify that $\boldsymbol{\varPhi}_0 \in \mathcal{EBF}^{d\times d}$ with L\'evy triplet $(\mathbf{0}, \boldsymbol{\sigma}^{\boldsymbol{\varPhi}},\boldsymbol{\nu}^{\boldsymbol{\varPhi}})$ and 
 \beqlb\label{eqn.224}
 \boldsymbol{\varphi}_0(\lambda)
 =\big( 	\boldsymbol{Q}^{\rm T}  \boldsymbol{\varPhi}_0(\lambda )\boldsymbol{Q}  \big)_{\mathtt{I}\mathtt{I}}	+\boldsymbol{U}_{\mathtt{I}\mathtt{J}} \big(  \mathbf{Id}- \boldsymbol{U}_{\mathtt{J}\mathtt{J}}\big)^{-1}\big( 	\boldsymbol{Q}^{\rm T}  \boldsymbol{\varPhi}_0(\lambda ) \boldsymbol{Q} \big)_{\mathtt{J}\mathtt{I}}. 
 \eeqlb

 \begin{theorem}\label{Thm.ConvergenceDS01}	
 For some constant  $\lambda_0^+\geq 0$, assume that  $ \boldsymbol{\varphi}_0(\lambda)$ is invertible for  any $\lambda > \lambda_0^+$ and  
 \beqlb\label{eqn.225}
 \big( 	\boldsymbol{Q}^{\rm T}  \boldsymbol{b}^{\boldsymbol{\varPhi}} \boldsymbol{Q} \big)_{\mathtt{I}\mathtt{J}}	+\boldsymbol{U}_{\mathtt{I}\mathtt{J}} \big(  \mathbf{Id}- \boldsymbol{U}_{\mathtt{J}\mathtt{J}}\big)^{-1}\big( 	\boldsymbol{Q}^{\rm T} \boldsymbol{b}^{\boldsymbol{\varPhi}} \boldsymbol{Q} \big)_{\mathtt{J}\mathtt{J}} =0,
 \eeqlb
 we have the following.
 \begin{enumerate}
 \item[(1)] The $(\boldsymbol{K},\boldsymbol{\varPhi}_0)$-potential measure $\boldsymbol{\varPi}_0(dt) \in  M_{\rm loc}(\mathbb{R}_+;\mathbb{R}_+^{d\times d}) $ exists and is uniquely determined by the Laplace transform (\ref{eqn.LimitResolvent}) with $\boldsymbol{\varphi}$ replaced by $\boldsymbol{\varphi}_0$.   
 		
 \item[(2)] The two potential measures $\boldsymbol{\varPi}(dt)$ and $\boldsymbol{\varPi}_0(dt)$ satisfy the following resolvent equation
 \beqlb\label{eqn.ResolventEquation}
 \boldsymbol{\varPi}_0 
 =\boldsymbol{\varPi} + \boldsymbol{\varPi}_0 * \big(\boldsymbol{b}^{\boldsymbol{\varPhi}}\cdot \boldsymbol{\varPi} \big) 
 = \boldsymbol{\varPi} + \boldsymbol{\varPi} * \big(\boldsymbol{b}^{\boldsymbol{\varPhi}}\cdot \boldsymbol{\varPi}_0 \big) . 
 \eeqlb
 		
 \item[(3)] If $\boldsymbol{\varUpsilon}= \boldsymbol{\varPi}* \boldsymbol{\varGamma}$  with $\boldsymbol{\varGamma} \in D(\mathbb{R}_+;\mathbb{R}_+^d)$, then the stochastic Volterra equation (\ref{Eqn.HawkesVoletrra}) is equivalent to 
 \beqlb \label{eqn.SVE0101}
 \boldsymbol{\varXi} \ar=\ar   \boldsymbol{\varPi}_0  * \boldsymbol{\varGamma} - \boldsymbol{\varPi}_0 * \big(\boldsymbol{b}^{\boldsymbol{\varPhi}}\cdot \boldsymbol{\varXi} \big) + \boldsymbol{\varPi}_0  * \boldsymbol{M} . 
 \eeqlb  
 \end{enumerate}	 
 \end{theorem}

 By using Fubini's theorem to \eqref{eqn.ResolventEquation}, we see that the absolute continuity is passed between the two measures $\boldsymbol{\varPi}(dt)$ and $\boldsymbol{\varPi}_0(dt)$. Moreover, their densities, denoted by $\boldsymbol{\pi}$ and $\boldsymbol{\pi}_0$, enjoy the same regularity and integrability. 
 The next result is a  direct consequence of Theorem~\ref{Thm.ConvergenceDS01} and Corollary~\ref{Corollary.2001}-\ref{Corollary.2002}. The proof is omitted.
 
 \begin{corollary}\label{Corollary.2.20}
 Under conditions in Theorem~\ref{Thm.ConvergenceDS01} with $\boldsymbol{\varGamma}$ having derivative $\boldsymbol{\varGamma}'$ on $\mathbb{R}$, the following hold.
 \begin{enumerate}
  \item[(1)] If $\boldsymbol{\pi}_0  \in L^2_{\rm loc}(\mathbb{R}_+;\mathbb{R}^{d\times d}_+)$, then the stochastic Volterra equation (\ref{eqn.SVE04}) is equivalent to 
  \beqlb\label{eqn.SVE0401}
  \boldsymbol{\xi} (t) 
  = \int_0^t \boldsymbol{\pi}_0 (t-s) \big( \boldsymbol{\varGamma}'(s) - \boldsymbol{b}^{\boldsymbol{\varPhi}}\cdot \boldsymbol{\xi} (s) \big)ds + \int_0^t \boldsymbol{\pi}_0 (t-s)\cdot  {\rm diag } \big(\sqrt{\boldsymbol{\xi} (s)}\big)d\boldsymbol{B}(s)  , \quad  \mathbf{P} \otimes L(dt)\mbox{-a.s.}  
  \eeqlb
 		
  \item[(2)] If $\boldsymbol{\pi}_0 $ has derivative $\boldsymbol{\pi}_0'$ on $(0,\infty)$, then the stochastic Volterra equation (\ref{eqn.SVE03}) is equivalent to 
  \beqlb\label{eqn.SVE0301}
   \left\{\begin{split}
  \boldsymbol{\xi} (t) 
  &=  \int_0^t \boldsymbol{\pi}'_0(t-s)\boldsymbol{M}_0(s)ds +\boldsymbol{\pi}(0) \cdot  \boldsymbol{M}_0(t)  , \quad  \mathbf{P} \otimes L(dt)\mbox{-a.s.} ,\cr
  \boldsymbol{M}_0(t) &=	\boldsymbol{\varGamma}'(t)-\int_0^t \boldsymbol{b}^{\boldsymbol{\varPhi}}\cdot \boldsymbol{\xi}(s) ds +\int_0^t {\rm diag } \big(\sqrt{\boldsymbol{\xi} (s)}\big)d\boldsymbol{B}(s) . 
  \end{split}\right. 
  \eeqlb
 \end{enumerate} 
 \end{corollary}
 
 \begin{remark}\label{Remark.2.23}
 When $\boldsymbol{\underline\varPi}(0)$ is non-diagonal, we have the following. 
 \begin{enumerate}
 	\item[(1)] Our proof of Lemma~\ref{Lemma.ExistUniqueV} about the uniqueness of solutions to \eqref{Eqn.AffineVolterra} will no longer work, but the proof of existence  remains valid; see Lemma~\ref{Lemma.ExitenceRV}.  
 	Fortunately, the sequence $\{\big(\mathcal{I}_{\boldsymbol{\varLambda}^{(n)}}, \boldsymbol{N}^{(n)}, \widetilde{\boldsymbol{N}}^{(n)}\big)\}_{n\geq 1}$ is still relatively compact in $D_S(\mathbb{R}_+;\mathbb{R}_+^d\times\mathbb{R}_+^d\times\mathbb{R}^d)$ and all results still hold for any convergent sub-sequence and accumulation point. 
 	If all accumulation points are continuous, by Lemma~\ref{Lemma.FLFXiM} they have the same Fourier-Laplace functional \eqref{eqn.FL} and hence are identical in law.
 	
 	\item[(2)] All results remain hold if the function $\boldsymbol{\varUpsilon}$ is strictly increasing and there is a continuous accumulation point. 
 	Indeed, by Lemma~\ref{Lemma.FLFXiM}, the Fourier-Laplace functional \eqref{eqn.FL} holds for the continuous accumulation point associated with any solution $\boldsymbol{\mathcal{V}}$ of \eqref{Eqn.AffineVolterra}. This along with the strict monotonicity of $\boldsymbol{\varUpsilon}$ induces the uniqueness of $\boldsymbol{\mathcal{W}}$ and hence $\boldsymbol{\mathcal{V}}$.

 \end{enumerate} 
 \end{remark}
 
 \begin{remark}
 In the uni-variate case,  $\varPi_0(dt)$ is the potential measure of a subordinator with Laplace exponent $\varPhi_0$; see \cite{Bertoin1996}. 
 By Lemma~11.15 in \cite[p.166]{SchillingSongVondracek2012}, it has density $\pi_0$ if one of the following holds:
 \medskip\smallskip \\ \medskip\smallskip 
 \centerline{(i) $\sigma^\varPhi>0$; \quad (ii)  $\nu^\varPhi(\mathbb{R}_+)=\infty$ and $\nu^\varPhi(dt)<< L(dt) $.} 
 Let $\bar\nu^\varPhi(t):= \nu^\varPhi([t,\infty)) \in L^1_{\rm loc}(\mathbb{R}_+;\mathbb{R}_+)$ be the tail function of measure $\nu^\varPhi(dt)$.  
 \begin{enumerate}
 \item[(1)]   In case (i), by Proposition 1 and Corollary 2 in \cite{DoringSavov2011} with $q=0$, the density $\pi_0$ is differentiable and  the following representations hold:
 \beqnn
 \pi_0(t) = \frac{1}{\sigma^\varPhi} +  \sum_{k=1}^\infty \frac{(-1)^k}{|\sigma^\varPhi|^{k+1}}\cdot\int_0^t\big({\bar\nu}^\varPhi\big)^{*k}(s)ds
 \quad \mbox{and}\quad 
 \pi_0'(t) =     \sum_{k=1}^\infty \frac{(-1)^k}{|\sigma^\varPhi|^{k+1}}\cdot \big({\bar\nu}^\varPhi\big)^{*k}(t)
 ,\quad t\geq 0. 
 \eeqnn
 
 \item[(2)] In case (ii), by (3.1) in \cite{DoringSavov2011} with $q=0$ and $\delta=0$ we have  $ \pi_0 * \bar\nu^\varPhi \equiv 1$ on $(0,\infty)$. Moreover, the density $\pi_0 \in L^2_{\rm loc}((0,\infty);\mathbb{R}_+)$  if 
 \beqnn
 \int_0^1 \Big|\int_0^t\bar\nu^\varPhi(s)ds\Big|^{-2} dt<\infty .
 \eeqnn
 \end{enumerate} 
 \end{remark}
  
 \subsection{Examples}
 
 In this section, we provide a specific example of process obtained through the scaling limit of Hawkes processes.  
 Assume that $(\boldsymbol{N}_{n},\boldsymbol{\varLambda}_{n})$ has a constant exogenous rate $\boldsymbol{\mu}_{n}(t) \equiv \boldsymbol{\mu}_{n} \in \mathbb{R}_+^d $.   
 Consider a vector $\boldsymbol{a}\in\mathbb{R}_+^d$, a  matrix $\boldsymbol{b} \in \mathbb{R}^{d\times d}$ and an invertible diagonal matrix $\boldsymbol{c} \in\mathbb{R}_+^{d\times d} $. 
 For some constants  $\alpha \in(0,1]$ and  $\beta\geq 0$, 
 assume that the following limits hold as $n\to\infty$,
 \beqlb\label{eqn.201}
  n^{1-\alpha}\cdot \boldsymbol{\mu}_{n}\to \boldsymbol{a} ,\quad  n^{\alpha} \cdot \bigg( \mathbf{Id} -\int_0^\infty  e^{ \frac{\beta}{n}\cdot t} \boldsymbol{\phi}_{n}(s)ds \bigg) \to \boldsymbol{b}  
 \eeqlb
 and
 \beqnn  
 \begin{cases}
 \displaystyle n^{\alpha} \int_{nt}^\infty  e^{ \frac{\beta}{n}\cdot t} \boldsymbol{\phi}_{n}(s)ds \to     \frac{\boldsymbol{c}}{\alpha}\cdot t^{-\alpha}, & \mbox{if }\alpha \in(0,1); \vspace{7pt} \\
 \displaystyle	 \int_0^\infty s   e^{ \frac{\beta}{n}\cdot s} \boldsymbol{\phi}_{n}(s)ds \to  \boldsymbol{c}  , & \mbox{if }\alpha=1.  
 \end{cases}
 \eeqnn 
 
 In this setting, Assumption~\ref{AsymCriticality} and Condition~\ref{Main.Condition} hold with $
 \boldsymbol{K}=\mathbf{Id}$ and 
 $\boldsymbol{\varPhi}(\lambda)= \boldsymbol{b}   +  \boldsymbol{c}  \cdot (\lambda +\beta)^\alpha$ for $ \lambda \geq  -\beta $.  
 Moreover, we also have $\boldsymbol{\varphi} =\boldsymbol{\varPhi} $ is $\mathbf{Id}$-admissible with $ \boldsymbol{Q}=\boldsymbol{U}=\mathbf{Id}$ and $\lambda^+ >\sup_{i\in\mathtt{D}} \big|b_i/c_i \big|.$
 By Lemma~\ref{Lemma.Pi}(1) and Theorem~\ref{Thm.ConvergenceDS01}(1), the potential measure  $\boldsymbol{\varPi}_0(dt)$ has the  density 
 \beqnn
 \boldsymbol{\pi}_0(t)=   \frac{t^{\alpha-1}}{e^{\beta t}\Gamma(\alpha)}\cdot \boldsymbol{c}^{-1},\quad t>0. 
 \eeqnn 
 Applying Theorem~\ref{MainThm.J1WeakConvergence}, \ref{Thm.ConvergenceDS01}(3) and Corollary~\ref{Corollary.2.20} to the rescaled processes $\big(\mathcal{I}_{\boldsymbol{\varLambda}^{(n)}}, \boldsymbol{N}^{(n)}, \widetilde{\boldsymbol{N}}^{(n)} \big)$ defined by (\ref{eqn.ScaledProcess}) with $\theta_n=n^{2\alpha-1}$, we can get the following result. 

 \begin{corollary} \label{Coro.Example}
  We have  
 $ \big(\mathcal{I}_{\boldsymbol{\varLambda}^{(n)}}, \boldsymbol{N}^{(n)}, \widetilde{\boldsymbol{N}}^{(n)} \big)\overset{\rm w}\to \big(\boldsymbol{\varXi} ,\boldsymbol{\varXi} ,\boldsymbol{M} \big)$ in $D_{J_1}(\mathbb{R}_+;\mathbb{R}^d_+\times\mathbb{R}^d_+\times \mathbb{R}^d )$ as $n\to\infty$ with $(\boldsymbol{\varXi} ,\boldsymbol{M})$ being the unique weak solution to 
 \beqnn
 \boldsymbol{\varXi}(t) =   \int_0^t \frac{(t-s)^{\alpha-1}}{e^{\beta (t-s)}\Gamma(\alpha)}  \cdot \boldsymbol{c}^{-1}\cdot 
 \Big( \boldsymbol{a}\cdot s -\boldsymbol{b}  \cdot \boldsymbol{\varXi}(s) + \boldsymbol{M}(s)\Big) ds
 \quad \mbox{and}\quad 
 \boldsymbol{M}(t)= \boldsymbol{B}\circ\boldsymbol{\varXi}(t). 
 \eeqnn
 In particular, when $\alpha\in(1/2,1]$, the process $\boldsymbol{\varXi}$ is absolutely continuous with a continuous derivative   $\boldsymbol{\xi}$ admitting the following representation.
 \begin{enumerate}
 \item[(1)] If $\alpha\in(1/2,1)$, we have 
 \beqlb\label{MultiRoughCIR}
 \boldsymbol{\xi}(t) 
 \ar=\ar  \int_0^t \frac{(t-s)^{\alpha-1}}{e^{\beta (t-s)}\Gamma(\alpha)} \boldsymbol{c}^{-1} \Big( \boldsymbol{a}  -\boldsymbol{b}  \cdot \boldsymbol{\xi}(s) \Big)ds   + \int_0^t \frac{(t-s)^{\alpha-1}}{e^{\beta (t-s)}\Gamma(\alpha)} \boldsymbol{c}^{-1}  \sqrt{{\rm diag}(\boldsymbol{\xi}(s))}d\boldsymbol{B}(s) . 
 \eeqlb

 \item[(2)] If $\alpha=1$, we have 
 \beqlb\label{MultiCIR}
 \boldsymbol{\xi}(t) 
 \ar=\ar  \int_0^t \Big( \boldsymbol{c}^{-1}\boldsymbol{a}\cdot s -\big(\boldsymbol{c}^{-1}\boldsymbol{b} +	 \beta \cdot \mathbf{Id}\big) \cdot \boldsymbol{\xi}(s) \Big)ds   + \int_0^t  \boldsymbol{c}^{-1} \sqrt{{\rm diag}(\boldsymbol{\xi}(s))}d\boldsymbol{B}(s) .
 \eeqlb	
 \end{enumerate}
 \end{corollary}
 
 \begin{remark}
 If $\boldsymbol{\phi}_{n}$ has the form of $\boldsymbol{\phi}_{n}(t) = \boldsymbol{b}_n \cdot \kappa (1+t)^{-\kappa-1}\cdot e^{-\frac{\beta}{n}\cdot t} $  for $t\geq 0$ and satisfies (\ref{eqn.201}). 
 By generalizing the proof of Theorem~2.5 in \cite{HorstXuZhang2023a} if $\kappa \in (1/2,1)$ or using Theorem~3.6 in \cite{Xu2021} if $\kappa>1$, we can prove the weak convergence of the rescaled intensity processes $\{ \boldsymbol{\varLambda}^{(n)} \}_{n\geq 1}$ in $D_{J_1}(\mathbb{R}_+;\mathbb{R}^d_+)$ to the unique solution of (\ref{MultiRoughCIR}) with $\alpha=\kappa$ if $\kappa \in (1/2,1)$  or (\ref{MultiCIR}) if $\kappa >1$ respectively.  
 \end{remark}

  \section{Asymptotics of rescaled resolvent}
 \label{Sec.AsymResolvent}
 \setcounter{equation}{0}

 In this section, we study the long-run behavior of rescaled resolvent sequence $\{\boldsymbol{R}^{(n)}\}_{n\geq 1}$ to give the proof of Lemma~\ref{Lemma.Pi}. 
 As a preparation, in the next two propositions we first consider the invertibility of the following two matrix-valued functions:
 \beqnn
  \big(\mathbf{Id} -\boldsymbol{Q}^{\rm T}\mathcal{L}_{\boldsymbol{\phi}_n}(\lambda /n)\boldsymbol{Q} \big)_{\mathtt{J}\mathtt{J}} \in\mathbb{R}^{(d-\ell)\times(d-\ell)}
 \quad \mbox{and}\quad
 \boldsymbol{\varphi}^{(n)}(\lambda):=\sqrt{n\cdot \theta_n} \cdot \boldsymbol{ \varphi}_n(\lambda/n) \in\mathbb{R}^{\ell \times \ell} ,
 \quad \lambda \geq 0,
 \eeqnn
  with $\boldsymbol{ \varphi}_n(\lambda/n):= \boldsymbol{D}^{(n)}_{\mathtt{I}\mathtt{I}}(\lambda)-\boldsymbol{D}^{(n)}_{\mathtt{I}\mathtt{J}}(\lambda)\big(\boldsymbol{D}^{(n)}_{\mathtt{J}\mathtt{J}}(\lambda)\big)^{-1}\boldsymbol{D}^{(n)}_{\mathtt{J}\mathtt{I}}(\lambda)$ and $\boldsymbol{D}^{(n)}(\lambda):=\mathbf{Id}- \boldsymbol{Q}^{\rm T}\mathcal{L}_{\boldsymbol{\phi}_n}(\lambda/n )\boldsymbol{Q}$.

 \begin{proposition}\label{Prop.AuxAsymR.01}
	For each $\lambda\geq  0$, we have 	$\boldsymbol{D}^{(n)}(\lambda) \to  \mathbf{Id}-  \boldsymbol{U}$ as $n\to\infty$.
	Moreover, there exists an integer  $n_1 \geq 1$ such that $ \boldsymbol{D}^{(n)}_{\mathtt{J}\mathtt{J}}(\lambda)  $ is invertible for any $n\geq n_1$.  
 \end{proposition}
 \proof By Condition~\ref{Main.Condition} and then Definition~\ref{Def.Admissible}(1), we have $\mathcal{L}_{\boldsymbol{\phi}_n}(\lambda /n) \to \boldsymbol{K}$ and hence $\boldsymbol{Q}^{\rm T}\mathcal{L}_{\boldsymbol{\phi}_n}(\lambda /n)\boldsymbol{Q} \to \boldsymbol{U}$ as $n\to\infty$,
 which induces the first claim. 
 Since $( \mathbf{Id}-  \boldsymbol{U})_{\mathtt{J}\mathtt{J}}$ is invertible; see the paragraph below Definition~\ref{Def.Admissible},
 so is the matrix $\boldsymbol{D}^{(n)}_{\mathtt{J}\mathtt{J}}(\lambda)  $ for large $n\geq 1$. 
 \qed

\begin{proposition}\label{Prop.AuxAsymR.02}
 Recall $\lambda^+$ given in Definition~\ref{Def.Admissible}.	
 For each $\lambda \geq  \lambda^+$,  we have $\boldsymbol{\varphi}^{(n)}(\lambda) \to \boldsymbol{\varphi} (\lambda)$ as $n\to\infty$. 
	Moreover, there exists an integer  $n_2 \geq 1$ such that  $\boldsymbol{\varphi}^{(n)}(\lambda)$ and $\boldsymbol{ \varphi}_n(\lambda/n)$ are invertible for any $n\geq n_2$. 
\end{proposition}
 \proof The facts that $\boldsymbol{U}=\boldsymbol{Q}^{\rm T} \boldsymbol{K}\boldsymbol{Q}$,  $\boldsymbol{U}_{\mathtt{I}\mathtt{I}}=\mathbf{Id}$ and  $\boldsymbol{U}_{\mathtt{J}\mathtt{I}}=\mathbf{0}$ induce   that 
 \beqnn
 \boldsymbol{D}^{(n)}_{\mathtt{D}\mathtt{I}}(\lambda)
 =\big(  \boldsymbol{Q}^{\rm T}(\boldsymbol{K}-\mathcal{L}_{\boldsymbol{\phi}_n}(\lambda /n))\boldsymbol{Q}\big)_{\mathtt{D}\mathtt{I}}.
 \eeqnn
 By this and Condition~\ref{Main.Condition}, we have as $n\to\infty$,
 \beqnn
 \ar\ar\sqrt{n\cdot \theta_n} \cdot   
 \boldsymbol{D}^{(n)}_{\mathtt{I}\mathtt{I}}(\lambda)   \to \big(  \boldsymbol{Q}^{\rm T} \boldsymbol{\varPhi}(\lambda) \boldsymbol{Q}\big)_{\mathtt{I}\mathtt{I}},\quad 
 \boldsymbol{D}^{(n)}_{\mathtt{I}\mathtt{J}}(\lambda) \to   \boldsymbol{U}_{\mathtt{I}\mathtt{J}},\cr
 \ar\ar\cr
 \ar\ar  \sqrt{n\cdot \theta_n} \cdot\boldsymbol{D}^{(n)}_{\mathtt{J}\mathtt{I}}(\lambda) \to \big(  \boldsymbol{Q}^{\rm T} \boldsymbol{\varPhi}(\lambda) \boldsymbol{Q}\big)_{\mathtt{J}\mathtt{I}} ,    \quad  
\boldsymbol{D}^{(n)}_{\mathtt{J}\mathtt{J}}(\lambda) \to \mathbf{Id} - \boldsymbol{U}_{\mathtt{J}\mathtt{J}}. 
 \eeqnn
 Plugging them back into $\boldsymbol{\varphi}^{(n)}(\lambda)$   induces our first claim.
 Our second claim follows directly from the first  one as well as the invertibility of $\boldsymbol{\varphi}(\lambda)$ and Proposition~\ref{Appendix.Prop.LimitMatrix}(2).
\qed

 \begin{lemma}\label{Lemma.AuxAsymR.03}
  For each $\lambda\geq \lambda^+$, the matrix $ \boldsymbol{\varPsi}^{(n)} (\lambda)$ is invertible for any  $n\geq  (n_1\vee n_2)$ and 
  \beqlb\label{eqn.LimInveversePsi}
  \lim_{n\to\infty} \big(  \boldsymbol{\varPsi}^{(n)} (\lambda) \big)^{-1} 
  = \boldsymbol{Q}
    \begin{pmatrix}
    \big(\boldsymbol{\varphi}(\lambda)\big)^{-1} & 
    \big(\boldsymbol{\varphi}(\lambda)\big)^{-1} \boldsymbol{U}_{\mathtt{I}\mathtt{J}}\big(\mathbf{Id}-\boldsymbol{U}_{\mathtt{J}\mathtt{J}} \big)^{-1}\vspace{3pt} \\
    \mathbf{0} & \mathbf{0}
	\end{pmatrix}
  \boldsymbol{Q}^{\rm T}.
  \eeqlb
 \end{lemma}

 \proof By (\ref{eqn.varPsin}) and the fact that $\boldsymbol{Q}\boldsymbol{Q}^{\rm T}=\mathbf{Id}$, we see that  
 \beqlb\label{eqn.301}
 \boldsymbol{\varPsi}^{(n)} (\lambda) = \boldsymbol{Q}\cdot \sqrt{n\cdot\theta_n} \boldsymbol{D}^{(n)} (\lambda) \cdot \boldsymbol{Q}^{\rm T}. 
 \eeqlb
 This means that $ \boldsymbol{\varPsi}^{(n)} (\lambda)$ is  invertible if and only if $ \boldsymbol{D}^{(n)} (\lambda) $ is invertible.   
 By using Proposition~\ref{Appendix.Prop.BlockInversion}(2) along with Proposition~\ref{Prop.AuxAsymR.01} and \ref{Prop.AuxAsymR.02},  we have that $\boldsymbol{D}^{(n)} (\lambda)$ is invertible for any $n\geq n_1\vee n_2$ if and only if the matrix $  \boldsymbol{A}^{(n)}(\lambda)$ well-defined by
 \beqlb\label{eqn.302}
 \begin{pmatrix}
 \big(\boldsymbol{ \varphi}_n(\lambda/n)\big)^{-1} & -\big(\boldsymbol{ \varphi}_n(\lambda/n)\big)^{-1}\boldsymbol{D}^{(n)}_{\mathtt{I}\mathtt{J}}(\lambda)\big(\boldsymbol{D}^{(n)}_{\mathtt{J}\mathtt{J}}(\lambda)\big)^{-1} \vspace{7pt}\\ 
 -\big(\boldsymbol{D}^{(n)}_{\mathtt{J}\mathtt{J}}(\lambda)\big)^{-1} \boldsymbol{D}^{(n)}_{\mathtt{J}\mathtt{I}}(\lambda)\big(\boldsymbol{ \varphi}_n(\lambda/n)\big)^{-1} &
 \big(\boldsymbol{D}^{(n)}_{\mathtt{J}\mathtt{J}}(\lambda)\big)^{-1}\big( \mathbf{Id} +  \boldsymbol{D}^{(n)}_{\mathtt{J}\mathtt{I}}(\lambda)\big(\boldsymbol{ \varphi}_n(\lambda/n)\big)^{-1}\boldsymbol{D}^{(n)}_{\mathtt{I}\mathtt{J}}(\lambda)\big(\boldsymbol{D}^{(n)}_{\mathtt{J}\mathtt{J}}(\lambda)\big)^{-1}\big)
 \end{pmatrix}
 \eeqlb
 is invertible. 
 Indeed, the matrix  $  \boldsymbol{A}^{(n)}(\lambda)$ is invertible for any $n\geq n_1\vee n_2$ because  by Proposition~\ref{Appendix.Prop.BlockInversion}(1),
 \beqnn
 {\rm det} \big(\boldsymbol{A}^{(n)}(\lambda) \big) 
 = {\rm det} \Big(\big(\boldsymbol{ \varphi}_n(\lambda/n)\big)^{-1} \Big)  {\rm det} \Big(  \big(\boldsymbol{D}_{\mathtt{J}\mathtt{J}}^{(n)}(\lambda) \big)^{-1} \Big) \neq 0.
 \eeqnn
 
 We now start to prove the limit (\ref{eqn.LimInveversePsi}). 
 By (\ref{eqn.301}) and (\ref{eqn.302}), we can write $\big(\boldsymbol{\varPsi}^{(n)} (\lambda)\big)^{-1}$ as 
 \beqlb\label{eqn.303}
 \big(\boldsymbol{\varPsi}^{(n)} (\lambda)\big)^{-1}
 = \boldsymbol{Q}\cdot \frac{\big(\boldsymbol{D}^{(n)} (\lambda)\big)^{-1} }{\sqrt{n\cdot\theta_n}}\cdot \boldsymbol{Q}^{\rm T} 
 = \boldsymbol{Q}\cdot \frac{ \boldsymbol{A}^{(n)}(\lambda)}{\sqrt{n\cdot\theta_n}} \cdot \boldsymbol{Q}^{\rm T}
 \eeqlb
 and the matrix  $ \boldsymbol{A}^{(n)}(\lambda)/\sqrt{n\cdot\theta_n} $ has the representation
 \beqnn \begin{pmatrix}
 	\big( \boldsymbol{ \varphi}^{(n)}(\lambda) \big)^{-1} & -	\big( \boldsymbol{ \varphi}^{(n)}(\lambda) \big)^{-1}\boldsymbol{D}^{(n)}_{\mathtt{I}\mathtt{J}}(\lambda)\big(\boldsymbol{D}^{(n)}_{\mathtt{J}\mathtt{J}}(\lambda)\big)^{-1} \vspace{7pt} \\
 	-\big(\boldsymbol{D}^{(n)}_{\mathtt{J}\mathtt{J}}(\lambda)\big)^{-1} \boldsymbol{D}^{(n)}_{\mathtt{J}\mathtt{I}}(\lambda)	\big( \boldsymbol{ \varphi}^{(n)}(\lambda) \big)^{-1} &
 	\big(\boldsymbol{D}^{(n)}_{\mathtt{J}\mathtt{J}}(\lambda)\big)^{-1}\big( \frac{\mathbf{Id}}{\sqrt{n\cdot\theta_n} } +  \boldsymbol{D}^{(n)}_{\mathtt{J}\mathtt{I}}(\lambda)	\big( \boldsymbol{ \varphi}^{(n)}(\lambda) \big)^{-1}\boldsymbol{D}^{(n)}_{\mathtt{I}\mathtt{J}}(\lambda)\big(\boldsymbol{D}^{(n)}_{\mathtt{J}\mathtt{J}}(\lambda)\big)^{-1}\big)
 \end{pmatrix} 
 .
 \eeqnn
 Applying the limits in Proposition~\ref{Prop.AuxAsymR.01}, \ref{Prop.AuxAsymR.02} and Condition~\ref{Main.Condition} to the preceding blacks,  
 \beqnn
 \lim_{n\to\infty} \frac{ \boldsymbol{A}^{(n)}(\lambda)}{\sqrt{n\cdot\theta_n}}
 = \begin{pmatrix}
  \big(\boldsymbol{\varphi}(\lambda)\big)^{-1} & 
  \big(\boldsymbol{\varphi}(\lambda)\big)^{-1} \boldsymbol{U}_{\mathtt{I}\mathtt{J}}\big(\mathbf{Id}-\boldsymbol{U}_{\mathtt{J}\mathtt{J}} \big)^{-1} \vspace{3pt} \\
  \mathbf{0} & \mathbf{0}
 \end{pmatrix}.
 \eeqnn
 The limit (\ref{eqn.LimInveversePsi}) follows immediately by   taking this back into (\ref{eqn.303}). 
 \qed

 \textit{\textbf{Proof of Lemma~\ref{Lemma.Pi}.}} 
 By applying Condition~\ref{Main.Condition} and the limit (\ref{eqn.LimInveversePsi}) to  (\ref{eqn.LaplaceIntR}), we have as $n\to\infty$,
 \beqlb \label{eqn.401}
 \mathcal{L}_{\boldsymbol{R}^{(n)}}(\lambda ) 
 \to \boldsymbol{K} \boldsymbol{Q}
 \begin{pmatrix}
	\big(\boldsymbol{\varphi}(\lambda)\big)^{-1} & 
	\big(\boldsymbol{\varphi}(\lambda)\big)^{-1} \boldsymbol{U}_{\mathtt{I}\mathtt{J}}\big(\mathbf{Id}-\boldsymbol{U}_{\mathtt{J}\mathtt{J}} \big)^{-1} \vspace{3pt} \\
	\mathbf{0} & \mathbf{0}
 \end{pmatrix}
 \boldsymbol{Q}^{\rm T}, \quad  \lambda \geq \lambda^+. 
 \eeqlb 
 This yields that there exists a unique $\sigma$-finite measure on $\mathbb{R}_+$, denoted by $\boldsymbol{\varPi}(dt)$, that is independent of the auxiliary matrices $\boldsymbol{Q},\boldsymbol{U} \in \mathbb{R}^{d\times d}$ and auxiliary constant $\lambda_+\geq 0$ such that for large $\lambda \geq 0$,
 \beqlb\label{eqn.403}
 \lim_{n\to\infty} \mathcal{L}_{\boldsymbol{R}^{(n)}}(\lambda )  
 = \int_{\mathbb{R}_+} e^{-\lambda t} \boldsymbol{\varPi}(dt),
 \eeqlb
 which is equal to the limit function in (\ref{eqn.401}). 
 Consequently, the measure with density  $\boldsymbol{R}^{(n)}$ converges vaguely to $\boldsymbol{\varPi}(dt)$ as $n\to\infty$. 
 By the fact that $\boldsymbol{K}=\boldsymbol{Q} \boldsymbol{U}\boldsymbol{Q}^{\rm T} $, 
 \beqlb\label{eqn.402}
 \int_{\mathbb{R}_+} e^{-\lambda t} \boldsymbol{\varPi}(dt)=  \boldsymbol{Q} 
 \begin{pmatrix}
	\mathbf{Id} &	\boldsymbol{U}_{\mathtt{I}\mathtt{J}}\\
	\mathbf{0} & \boldsymbol{U}_{\mathtt{J}\mathtt{J}}
 \end{pmatrix} \begin{pmatrix}
	\big(\boldsymbol{\varphi}(\lambda)\big)^{-1} & 
	\big(\boldsymbol{\varphi}(\lambda)\big)^{-1} \boldsymbol{U}_{\mathtt{I}\mathtt{J}}\big(\mathbf{Id}-\boldsymbol{U}_{\mathtt{J}\mathtt{J}} \big)^{-1} \vspace{3pt} \\
	\mathbf{0} & \mathbf{0}
 \end{pmatrix}\boldsymbol{Q}^{\rm T}  
 \eeqlb  
 and then (\ref{eqn.LimitResolvent}) holds. 
 \qed

  \begin{corollary}\label{Prop.503}
 	For any row vector function $\boldsymbol{f}\in L^\infty_{\rm loc}(\mathbb{R}_+;\mathbb{C}^d)$ and $T\geq 0$, we have as $n\to\infty$,
 	\beqnn
 	\big\|\boldsymbol{f} * \boldsymbol{R}^{(n)}    -\boldsymbol{f}* \boldsymbol{\varPi}\big\|_{L^1_T} \to 0.
 	\eeqnn 
 \end{corollary}
 \proof  For any $T\geq 0$, one can always find a sequence of functions  $\{\boldsymbol{f}_\epsilon\}_{\epsilon>0}\subset C(\mathbb{R}_+;\mathbb{C}^d)$ such that 
 \beqnn 
 \sup_{\epsilon>0}\|\boldsymbol{f}_\epsilon \|_{L^\infty_T}  \leq 2\cdot \| \boldsymbol{f}\|_{L^\infty_T}
 \quad \mbox{and}\quad 
 \lim_{\epsilon\to 0+} \|\boldsymbol{f}_\epsilon-\boldsymbol{f}\|_{L^1_T} = 0.
 \eeqnn
 By the triangle inequality, 
 \beqnn
 \big\|\boldsymbol{f} * \boldsymbol{R}^{(n)}    -\boldsymbol{f}* \boldsymbol{\varPi}\big\|_{L^1_T} 
 \leq 	
 \big\| (\boldsymbol{f}  -\boldsymbol{f}_\epsilon ) *  \boldsymbol{R}^{(n)}\big\|_{L^1_T}
 + 	\big\| ( \boldsymbol{f}_\epsilon    -\boldsymbol{f} ) * \boldsymbol{\varPi}\big\|_{L^1_T} 
 + 	\big\|\boldsymbol{f}_\epsilon * \boldsymbol{R}^{(n)}    -\boldsymbol{f}_\epsilon* \boldsymbol{\varPi}\big\|_{L^1_T}.
 \eeqnn 
 Applying Proposition~\ref{Prop.A1}(2) to the first two terms on the right side of this inequality and then using (\ref{Eqn.BoundConvIRn}),
 \beqnn
 \big\| (\boldsymbol{f}  -\boldsymbol{f}_\epsilon ) *  \boldsymbol{R}^{(n)}\big\|_{L^1_T} \leq \big\|\boldsymbol{f}_\epsilon-\boldsymbol{f} \big\|_{L^1_T} \cdot \big\|\boldsymbol{R}^{(n)} \big\|_{L^1_T}
 \quad \mbox{and}\quad 
 \big\| (\boldsymbol{f}_\epsilon  -\boldsymbol{f} ) * \boldsymbol{\varPi}\big\|_{L^1_T} \leq  \big\| \boldsymbol{f}_\epsilon-\boldsymbol{f} \big\|_{L^1_T} \cdot   \underline{\boldsymbol{\varPi}}(T),
 \eeqnn
 which go to $0$ as $\epsilon \to 0 +$. 
 On the other hand, by the continuity of $\boldsymbol{f}_\epsilon$ and the vague convergence of $\boldsymbol{R}^{(n)}(t)dt$ to $\boldsymbol{\varPi}(dt)$,  we have for each $t\geq 0$ and $\epsilon>0$,
 \beqnn
 \lim_{n\to\infty}\Big| \int_0^t\boldsymbol{f}_\epsilon(t-s) \boldsymbol{R}^{(n)} (s)ds -  \int_0^t\boldsymbol{f}_\epsilon(t-s) \boldsymbol{\varPi}(ds)  \Big| = 0.
 \eeqnn 
 By the dominated convergence theorem along with the local
 boundedness of $\boldsymbol{f}_\epsilon$, we have 
 \beqnn
 \lim_{n\to\infty} 	\big\|\boldsymbol{f}_\epsilon * \boldsymbol{R}^{(n)}    -\boldsymbol{f}_\epsilon* \boldsymbol{\varPi}\big\|_{L^1_T} \leq  \int_0^T \lim_{n\to\infty}\Big| \int_0^t\boldsymbol{f}_\epsilon(t-s) \boldsymbol{R}^{(n)} (s)ds -  \int_0^t\boldsymbol{f}_\epsilon(t-s) \boldsymbol{\varPi}(ds) \Big|dt =0 .
 \eeqnn
 Putting all preceding estimates together, we can get the desired result immediately.
 \qed  
 
 \begin{remark}\label{Remark.3.5}
 Assume that $\{ \boldsymbol{\nu}^{(n)} \}_{n\geq 1} \subset M_{\rm loc}(\mathbb{R}_+;\mathbb{R}_+^{d\times k})$ converges vaguely to $\boldsymbol{\nu} \in M_{\rm loc}(\mathbb{R}_+;\mathbb{R}_+^{d\times k})$. 
 For any row vector functions $\boldsymbol{f}\in L^\infty_{\rm loc}(\mathbb{R}_+;\mathbb{C}^d)$ and $T\geq 0$, we have $ \big\|\boldsymbol{f} * \boldsymbol{\nu}^{(n)}    -\boldsymbol{f}* \boldsymbol{\nu}\big\|_{L^1_T} \to 0$ as $n\to\infty$.  
 \end{remark}

 \section{Fourier-Laplace functionals}
 \label{Sec.ConvergenceFLF}
 \setcounter{equation}{0}
  
  In this section,  we establish the explicit exponential-affine representations of Fourier-Laplace functionals of the Hawkes process $(\boldsymbol{\widetilde{N}},\boldsymbol{\varLambda})$, the rescaled process $(\boldsymbol{\widetilde{N}}^{(n)}, \boldsymbol{\varLambda}^{(n)})$ and the solution $(\boldsymbol{\varXi},\boldsymbol{M})$ of (\ref{Eqn.HawkesVoletrra.02}). 
  Meanwhile, the existence and uniqueness of solutions to the Riccati-Volterra equation (\ref{Eqn.AffineVolterra}) are also considered. 
  They will play an important role in the proofs of our main results. 
  As a preparation, we refer the reader to Appendix~\ref{Appendix-Volterra} for elements and terminologies of Volterra integrals and equations.
  
 For row vector functions $\boldsymbol{f} \in L^\infty_{\rm loc}(\mathbb{R}_+;\mathbb{C}_-^d) $ and $\boldsymbol{h}\in L^\infty_{\rm loc}(\mathbb{R}_+;\mathtt{i}\mathbb{R}^d)$, by Proposition~\ref{Prop.A2} there exits a unique non-continuable solution $(T_\infty, \boldsymbol{\tilde{V}})\in (0,\infty] \times  C([0,T_\infty);\mathbb{C}^d)$  to the following two equivalent nonlinear Volterra equations 
 \beqlb\label{VolRiccati}
 \boldsymbol{V} =  \boldsymbol{W}   *\boldsymbol{R}
 \quad \mbox{and}\quad
 \boldsymbol{V} = \big(  \boldsymbol{W} +\boldsymbol{V} \big) * \boldsymbol{\phi} 
 \quad \mbox{with}\quad 
 \boldsymbol{W} := \boldsymbol{f}  + \mathcal{W}(\boldsymbol{V} + \boldsymbol{h}) ,
 \eeqlb  
 where $\mathcal{W}(\boldsymbol{x}):= \big(e^{x_i}-1- x_i \big)_{i\in\mathtt{D}}$ 
 for $\boldsymbol{x}\in\mathbb{C}^d$. 
 The existence of a unique global solution will be established later. 
 Based on the martingale representation~(\ref{SVR}), in the next proposition we establish an exponential-affine representation of Fourier-Laplace functional of $(\boldsymbol{N},\boldsymbol{\varLambda})$.


   \begin{proposition}\label{Lemma.FourLapGHP}
   	For row vector functions $\boldsymbol{f} \in L^\infty_{\rm loc}(\mathbb{R}_+;\mathbb{C}_-^d) $ and $\boldsymbol{h}\in L^\infty_{\rm loc}(\mathbb{R}_+;\mathtt{i}\mathbb{R}^d)$, let $(T_\infty,\boldsymbol{V})\in (0,\infty]  \times C([0,T_\infty);\mathbb{C}^d) $ be the non-continuable solution of (\ref{VolRiccati}), we have 
   	\beqlb\label{FourLapFunHP}
   	\mathbf{E}\Big[\exp\big\{   \boldsymbol{f}*\boldsymbol{\varLambda}(T)   + \boldsymbol{h}* d\widetilde{\boldsymbol{N}}(T) \big\}\Big]
   	= \exp \big\{  \boldsymbol{W}  * \boldsymbol{H}(T) \big\},  \quad T\in[0,T_\infty) . 
   	\eeqlb

   \end{proposition}
   \proof 
 For each $T\in(0,T_\infty)$, we introduce a semimartingale $ \{Z_T(t):t\in[0,T]\}$ by
   \beqlb\label{eqn.3002}
   Z_T(t)
   \ar:=\ar  \boldsymbol{W} * \boldsymbol{H}  (T) -  \int_0^t  \big(\boldsymbol{W}-\boldsymbol{f} \big) (T-s) \boldsymbol{\varLambda}(s)ds + \int_0^t  \big(  \boldsymbol{V} +\boldsymbol{h} \big)(T-s) d\widetilde{\boldsymbol{N}}(s).
   \eeqlb
   Applying It\^{o}'s formula to $\exp\{Z_T(t) \}$ and then using the definition of $\boldsymbol{W}$, we obtain that
   \beqlb\label{eqn.3008}
   e^{Z_T(t)}= e^{\mathbf{E}[Y_T]} + \int_0^t e^{Z_T(s-)} \big(e^{\boldsymbol{V}+  \boldsymbol{h} } -\mathbf{1}\big) (T-s)  d\widetilde{\boldsymbol{N}}(s),\quad t\in[0,T] , 
   \eeqlb
   which is a $(\mathscr{F}_t)$-local martingale.  If it were a true martingale, then the desired result follows from 
   \beqnn
   Z_T(T)=  
   \boldsymbol{f}*\boldsymbol{\varLambda}(T)   + \boldsymbol{h} * d\widetilde{\boldsymbol{N}}(T)
   \eeqnn 
   by taking expectations on both sides of (\ref{eqn.3008}) with $t=T$. 
   
   It remains to prove that $\{ \exp\{Z_T(t) \}: t\in[0,T]  \}$ is a true $(\mathscr{F}_t)$-martingale. 
   We first define the process
   \beqnn
   U_T(t):=  \int_0^t \big(e^{\boldsymbol{V} + \boldsymbol{h} } -\mathbf{1}\big) (T-s)  d\widetilde{\boldsymbol{N}}(s), \quad t\geq 0.
   \eeqnn
   In view of the local integrability of $\boldsymbol{H}$ and the local boundedness of $(\boldsymbol{V},\boldsymbol{g})$, 
   it follows from the Burkholder-Davis-Gundy inequality that 
   \beqnn
   \sup_{t\in[0,T]}\mathbf{E}\Big[ \big| U_T(t) \big|^2\Big]
   \ar\leq\ar \sum_{i=1}^d \sup_{t\in[0,T]} \int_0^t\mathbf{E} \big[\varLambda_i(s)\big] \cdot \Big| e^{(V_i +h_i)(T-s)}-1\Big|^2 ds \cr 
   \ar\leq\ar  \sum_{i=1}^d \sup_{t\in[0,T]} \Big| e^{(V_i +h_i)(t)}-1\Big|^2 \cdot \mathcal{I}_{H_i}(T) <\infty
   \eeqnn
   and hence that $ \{U_T(t):t\in[0,T]\}$ is a locally uniformly square integrable $(\mathscr{F}_t)$-martingale. Let
   \beqnn
   \mathcal{E}_{U_T}:= \big\{ \mathcal{E}_{U_T}(t): t\geq 0 \big\} 
   \eeqnn
   be the Dol\'ean-Dade exponential of $U_T$. By It\^o's formula,
   \beqnn
   \mathcal{E}_{U_T}(t)
   =\exp\Big\{ -\int_0^t  \big(\boldsymbol{W}-\boldsymbol{f} \big)(T-s)\boldsymbol{\Lambda}(s) ds + \int_0^t  (\boldsymbol{V} +\boldsymbol{h})(T-s) d\widetilde{\boldsymbol{N}}(s)  \Big\}, \quad t\geq 0.
   \eeqnn
 From this and (\ref{eqn.3002}), we see that $ e^{Z_T(t)}=e^{\mathbf{E}[Y_T]} \cdot \mathcal{E}_{U_T}(t) $ for $t\in[0,T]$.
 Since $ \mathcal{E}_{U_T}$ is a non-negative local martingale and hence a super-martingale, so $\mathbf{E}[\mathcal{E}_{U_T}(t)] \leq 1$.
 Thus it suffices to prove that 
 \beqnn
 \mathbf{E}\big[ \mathcal{E}_{U_T}(t) \big]=1, \quad t\geq 0.
 \eeqnn
   
   To this end, we introduce, for each $t_0\geq 0$ and $n\geq 1$ the quantities
   \beqnn
   \tau_n:=\inf\Big\{ s\geq 0: \max_{i\in\mathtt{D}}\mathcal{I}_{\Lambda_i}(s)\geq n \Big\} \wedge t_0
   \quad\mbox{and}\quad
   \mathcal{E}_{U_T}^{(n)}(t):=  \mathcal{E}_{U_T}(\tau_n \wedge t), 
   \quad t\geq 0. 
   \eeqnn  
   Since  $(\boldsymbol{V},\boldsymbol{h})$ is locally bounded, there exists a constant $C>0$ such that uniformly in $t \in [0,T]$ and $i\in\mathtt{D}$ the following holds:
   \beqnn
   \int_0^t \mathbf{1}_{\{s\leq \tau_n\}}  \varLambda_i (s) \bigg| 1-\big( 1-(V_i +h_i)(T-s) \big)      e^{(V_i +h_i)(T-s) }  \bigg|  ds 
   \leq C\cdot \int_0^{\tau_n}\varLambda_i(s) ds \leq C\cdot n. 
   \eeqnn
   Hence, the process $\mathcal{E}_{U_T}^{(n)}$ is a martingale for each $n\geq 1$, due to Theorem~IV.3 in \cite{LepingleMemin1978}. 
   Thus, 
   \beqnn 
   1 = \mathbf{E}\big[\mathcal{E}_{U_T}^{(n)}(t_0)\big]
   \ar=\ar \mathbf{E}\big[\mathcal{E}_{U_T}^{(n)}(t_0); \tau_n = t_0\big] + \mathbf{E}\big[\mathcal{E}_{U_T}^{(n)}(t_0);\tau_n<t_0\big]  \\
   \ar=\ar \mathbf{E}\big[\mathcal{E}_{U_T}(t_0);\tau_n= t_0\big] + \mathbf{E}\big[\mathcal{E}_{U_T}^{(n)}(\tau_n);\tau_n<t_0\big] . 
   \eeqnn
   By the monotone convergence theorem and the fact that $\tau_n \to t_0$ a.s.~as $n\to\infty$, 
   \beqnn
   \mathbf{E}\big[\mathcal{E}_{U_T}(t_0); \tau_n=t_0 \big]\to  \mathbf{E}\big[ \mathcal{E}_{U_T}(t_0) \big]
   \eeqnn
   and so it suffices to prove that $\mathbf{E}[\mathcal{E}_{U_T}^{(n)}(\tau_n);\tau_n<t_0] \to 0$ as $n\to\infty$. 
   
   To this end, we define a probability law  $\mathbf{Q}_T^{(n)}$ on $(\Omega,\mathscr{F},\mathscr{F}_t)$ by 
   \beqnn
   \frac{d  \mathbf{Q}_T^{(n)} }{ d \mathbf{P} } = \mathcal{E}_{U_T}^{(n)} (\tau_n). 
   \eeqnn
   By the definition of $\tau_n$ and Chebyshev's inequality it follows that
   \beqnn
   \mathbf{E}\big[\mathcal{E}_{U_T}^{(n)}(\tau_n); \tau_n<t_0 \big]  =\mathbf{Q}_T^{(n)} \big(\tau_n<t_0 \big) = \mathbf{Q}_T^{(n)} \Big( \max_{i\in\mathtt{D}}\mathcal{I}_{\varLambda_i} (t_0) \geq n\Big) \leq \frac{1}{n}\sum_{i=1}^d \mathbf{E}^{\mathbf{Q}_T^{(n)}}\big[ \mathcal{I}_{\varLambda_i}(t_0)\big]  
   \eeqnn
   and hence the desired result holds if we can establish a uniform upper bound on $\mathbf{E}^{\mathbf{Q}_T^{(n)}}\big[ \mathcal{I}_{\varLambda_i}(t_0)\big]$. In what follows we prove that there exists $\beta > 0$ such that
   \beqnn
   \sup_{n\geq 1}  \mathbf{E}^{\mathbf{Q}^{(n)}_T}\big[  \mathcal{I}_{\boldsymbol{\varLambda} }(t_0) \big] 
   \leq  2 e^{\beta t_0}  \mathcal{I}_{\boldsymbol{H}} (t_0).
   \eeqnn
   
   As  in~\cite[Section~2]{HorstXu2021}, on an extension of the original probability space we can define $d$ independent Poisson random measures (PRMs) $\{N_{0,i}(ds,dz): i\in\mathtt{D}\}$ on $(0,\infty)^2$ with intensity $ds dz$ such that 
   \beqlb
   N_i(t)= \int_0^t \int_0^{\varLambda_i(s-)} N_{0,i}(ds,dz)
   \quad{and} \quad
   \varLambda_i(t) =  H_i (t)+ \sum_{j=1}^d \int_0^t  \int_0^{\varLambda_i(s-)}R_{ij}(t-s)  \widetilde{N}_{0,j}(ds,dz ), \label{eqn.3001}
   \eeqlb
   for any $t\geq0$ and $i\in\mathtt{D}$, where $\widetilde{N}_{0,j}(ds,dz ):= N_{0,j}(ds,dz )-dsdz$ is the compensated random measure.  By Girsanov's theorem for random measures as stated in, e.g.~Theorem~3.17 in \cite[p.170]{JacodShiryaev2003}, the PRM $N_{0,i}(ds,dz)$ is a random point measure under $\mathbf{Q}_T^{(n)}$ with intensity 
   \beqnn
   \mathbf{1}_{\{ s\leq \tau_n \}}\cdot \exp \{ (V_i+h_i)(T-s) \}  ds dz.
   \eeqnn
   Moreover, the second equation in (\ref{eqn.3001}) - which holds a.s.~under the measure $\mathbf{P}$ - induces the following equality in distribution under the measure $\mathbf{Q}_T^{(n)}$: 
   \beqnn
   \varLambda_i (t)\ar \overset{\rm d}=\ar   H_i (t) + \sum_{j=1}^d \int_0^t  \int_0^{\varLambda_j(s-)}R_{ij}(t-s)  \widetilde{N}_{0,j}(ds,dz ) \cr
   \ar\ar    + \sum_{j=1}^d \int_0^t \mathbf{1}_{\{ s\leq \tau_n \}} R_{ij}(t-s) \varLambda_j (s) \big(e^{(V_j +h_j)(T-s)}-1 \big)ds , \quad i\in\mathtt{D}, t \geq 0.  
   \eeqnn 
   Taking expectations on both sides of the above equation and then integrating them over $[0,t]$ yields that
   \beqnn
   \int_0^t \mathbf{E}^{\mathbf{Q}^{(n)}_T}\big[ \boldsymbol{\varLambda} (s) \big]ds  
   \ar\leq\ar \mathcal{I}_{\boldsymbol{H}}(t) + C_0  \int_0^t dr \int_0^r  \boldsymbol{R}(t-s)   \mathbf{E}^{\mathbf{Q}^{(n)}_T}\big[  \boldsymbol{\varLambda} (s) \big] ds, 
   \eeqnn
   for some constant $C_0>0$ that is independent of $t$ and $n$. 
   Applying Fubini's theorem to the term on the left-hand side of the equality and to the double integral on the right-hand side of the inequality and then using the fact that $ \mathcal{I}_{\Lambda }$ is non-decreasing we see that
   \beqnn
   \mathbf{E}^{\mathbf{Q}^{(n)}_T}\big[  \mathcal{I}_{\boldsymbol{\varLambda} }(t) \big] 
   \ar\leq\ar   \mathcal{I}_{\boldsymbol{H}}(t) + C_0 \int_0^t\boldsymbol{R}(t-s)   \mathbf{E}^{\mathbf{Q}^{(n)}_T}\big[  \mathcal{I}_{\boldsymbol{\varLambda} }( s) \big]  ds.  
   \eeqnn  
   From (\ref{Resolvent}),  for any constant $\beta>0$ such that the norm  
   \beqnn
   \sum_{i,j=1}^d\int_0^\infty e^{-\beta t}  \phi_{ij}(t)dt \leq \frac{1/2}{1+C_0},
   \eeqnn
   we have 
   \beqnn
   \sum_{i,j=1}^d \int_0^\infty e^{-\beta t}R_{ij}(t)dt   \leq   \sum_{i,j=1}^d\int_0^\infty e^{-\beta t}  \phi_{ij}(t)dt \Big(1-\sum_{i,j=1}^d\int_0^\infty e^{-\beta t}  \phi_{ij}(t)dt\Big)^{-1} \leq \frac{1}{1+2C_0}.
   \eeqnn 
   and hence 
   \beqnn
   \sup_{t\in[0,t_0]} e^{-\beta t} \mathbf{E}^{\mathbf{Q}^{(n)}_T}\big[  \mathcal{I}_{\boldsymbol{\varLambda} }(t) \big]  
   \ar\leq\ar \sup_{t\in[0,t_0]} e^{-\beta t} \mathcal{I}_{\boldsymbol{H}} (t) +\sup_{t\in[0,t_0]}  C_0 \int_0^t  e^{-\beta (t-s)} \boldsymbol{R}(t-s) e^{-\beta s}  \mathbf{E}^{\mathbf{Q}^{(n)}_T}\big[  \mathcal{I}_{\boldsymbol{\varLambda} }( s) \big]  ds\cr 
   \ar\leq\ar \sup_{t\in[0,t_0]} e^{-\beta t} \mathcal{I}_{\boldsymbol{H}} (t) + C_0\cdot    \int_0^\infty e^{-\beta t}\boldsymbol{R}(t)dt   \cdot \sup_{t\in[0,t_0]} e^{-\beta t}  \mathbf{E}^{\mathbf{Q}^{(n)}_T}\big[  \mathcal{I}_{\boldsymbol{\varLambda} }(t) \big]\cr
   \ar\leq\ar  \sup_{t\in[0,t_0]} e^{-\beta t} \mathcal{I}_{\boldsymbol{H}} (t) + C_0\cdot    \sum_{i,j=1}^d \int_0^\infty e^{-\beta t}R_{ij}(t)dt \cdot \sup_{t\in[0,t_0]} e^{-\beta t}  \mathbf{E}^{\mathbf{Q}^{(n)}_T}\big[  \mathcal{I}_{\boldsymbol{\varLambda} }(t) \big]\cr
   \ar\leq\ar  \mathcal{I}_{\boldsymbol{H}} (t_0) + \frac{1}{2}\cdot  \sup_{t\in[0,t_0]} e^{-\beta t}  \mathbf{E}^{\mathbf{Q}^{(n)}_T}\big[  \mathcal{I}_{\boldsymbol{\varLambda} }(t) \big] ,
   \eeqnn
   which yields that $\sup_{n\geq 1} \sup_{t\in[0,t_0]} e^{-\beta t}  \mathbf{E}^{\mathbf{Q}^{(n)}_T}\big[  \mathcal{I}_{\boldsymbol{\varLambda} }(t) \big] \leq 2\cdot  \mathcal{I}_{\boldsymbol{H}} (t_0)$.
   Consequently, 
   \beqnn
   \sup_{n\geq 1}  \mathbf{E}^{\mathbf{Q}^{(n)}_T}\big[  \mathcal{I}_{\boldsymbol{\varLambda} }(t_0) \big] 
   \leq e^{\beta t_0} \sup_{n\geq 1} \sup_{t\in[0,t_0]} e^{-\beta t}  \mathbf{E}^{\mathbf{Q}^{(n)}_T}\big[  \mathcal{I}_{\boldsymbol{\varLambda} }(t) \big]
   \leq  2 e^{\beta t_0}  \mathcal{I}_{\boldsymbol{H}} (t_0).
   \eeqnn
   \qed
   
     \begin{lemma}  \label{Lemma.FourierLaplaceHP}
    	For row vector functions $\boldsymbol{f} \in L^\infty_{\rm loc}(\mathbb{R}_+;\mathbb{C}_-^d) $ and $\boldsymbol{h}\in L^\infty_{\rm loc}(\mathbb{R}_+;\mathtt{i}\mathbb{R}^d)$, we have the following.
    	\begin{enumerate}
    		\item[(1)] The equivalent nonlinear Volterra equations in (\ref{VolRiccati}) have a unique solution $\boldsymbol{V} \in C(\mathbb{R}_+;\mathbb{C}^d_-)$.
    		
    		\item[(2)] The Fourier-Laplace functional of  $(\boldsymbol{\varLambda} ,\boldsymbol{\widetilde{N}})$ admits the representation (\ref{FourLapFunHP}) for any $T\geq 0.$ 
    		
    		\item[(3)] $\mathtt{Re}\big(\boldsymbol{f}*\boldsymbol{\phi}\big)-2  \cdot \mathbf{1}*\boldsymbol{\phi}   \leq \mathtt{Re}(\boldsymbol{V}) \leq \mathtt{Re}\big(\boldsymbol{f}*\boldsymbol{\phi}\big) $ 
    		and 
    		$|\mathtt{Im}(\boldsymbol{V})| \leq  |\mathtt{Im}\big((\boldsymbol{f}-\boldsymbol{h})*\boldsymbol{\phi}\big)| + \mathbf{1} *\boldsymbol{\phi} $.
    	\end{enumerate}
    \end{lemma}
   \proof
   The existence and uniqueness of a non-continuable solution $(\boldsymbol{V},T_\infty)$ follows from standard arguments. 
   In view of Lemma~\ref{Lemma.FourLapGHP} it remains to prove that  $T_\infty=\infty$ and  that $ {\tt Re}\boldsymbol{V} \leq 0$. 
   For each fixed $j\in \mathtt{D}$, let 
   \beqnn
   \boldsymbol{\phi}_{\cdot j}:=(\boldsymbol{\phi}_{i j})_{i\in\mathcal{D}}^{\rm T}
   \quad \mbox{and}\quad 
   \boldsymbol{R}_{\cdot j}:=(\boldsymbol{R}_{i j})_{i\in\mathcal{D}}^{\rm T}
   \eeqnn
   be the $j$-columns of $\boldsymbol{\phi}$ and $\boldsymbol{R}$ respectively. 
   Choosing $\boldsymbol{\mu}=\boldsymbol{\phi}_{\cdot j}$ we have
   \beqnn
   \boldsymbol{H}= \boldsymbol{\phi}_{\cdot j}+\boldsymbol{R}* \boldsymbol{\phi}_{\cdot j}=\boldsymbol{R}_{\cdot j}
   \quad\mbox{and}\quad
   \boldsymbol{W} * \boldsymbol{H} = \boldsymbol{W} * \boldsymbol{R}_{\cdot j}= V_j.
   \eeqnn
   Taking these back into equation (\ref{FourLapFunHP}) shows that 
   \beqnn
   \mathbf{E}\Big[\exp \big\{  \boldsymbol{f} *  \boldsymbol{\Lambda}(T) + \boldsymbol{h} * d\widetilde{\boldsymbol{N}}(T) \big\}\Big]
   = e^{  V_j(T) },
   \quad  T\in[0,T_\infty). 
   \eeqnn 
   It hence follows that 
   \beqnn
   \exp\big\{{\tt Re}V_j(T)\big\} = \big|\exp\big\{  V_j(T) \big\}\big|
   \ar\leq\ar \mathbf{E}\Big[\Big|\exp\Big\{  \boldsymbol{f} *  \boldsymbol{\Lambda}(T) + \boldsymbol{h} * d\widetilde{\boldsymbol{N}}(T) \Big\}\Big|\Big]\cr
   \ar=\ar\mathbf{E}\Big[ \exp\Big\{ {\tt Re}\Big( \boldsymbol{f} *  \boldsymbol{\Lambda}(T) + \boldsymbol{h} * d\widetilde{\boldsymbol{N}}(T) \Big) \Big\} \Big]\leq 1 , 
   \eeqnn
   and hence that ${\tt Re}V_j(T) \leq 0$. Furthermore, from (\ref{VolRiccati}) we can see that 
   \beqnn
   {\tt Re}(V_j)= {\tt Re}\big((\boldsymbol{f}-\boldsymbol{h} )*\boldsymbol{\phi}_{\cdot j}\big)+ \Big(\exp\{{\tt Re}(V_i +h_i)\}\cos ({\tt Im}(V_i +h_i)) -1\Big)_{i\in\mathtt{D}}*\boldsymbol{\phi}_{\cdot j}.
   \eeqnn
   Since ${\tt Re}(V_i +h_i)={\tt Re}(V_i )\leq 0$, we have that $ {\tt Re}\big((\boldsymbol{f}-\boldsymbol{h} )*\boldsymbol{\phi}_{\cdot j}\big)={\tt Re}\big(\boldsymbol{f} *\boldsymbol{\phi}_{\cdot j}\big)$ and
   \beqnn
   {\tt Re}\big(\boldsymbol{f} *\boldsymbol{\phi}_{\cdot j}\big)
   \geq {\tt Re}(V_j)\geq
   {\tt Re}\big(\boldsymbol{f} *\boldsymbol{\phi}_{\cdot j}\big) -2\cdot\big(\mathbf{1}*\boldsymbol{\phi}_{\cdot j} \big).
   \eeqnn
   Similarly, we also have 
   \beqnn
   {\tt Im}(V_j) = {\tt Im} \big((\boldsymbol{f}-\boldsymbol{h})*\boldsymbol{\phi}_{\cdot j}\big)+ \Big(\exp\{{\tt Re}(V_i +h_i)\}\sin \big({\tt Im}(V_i +h_i)\big) \Big)_{i\in\mathtt{D}}*\boldsymbol{\phi}_{\cdot j},
   \eeqnn
   and $ |{\tt Im}(V_j)| \leq \big|{\tt Im} \big((\boldsymbol{f}-\boldsymbol{h})*\boldsymbol{\phi}_{\cdot j}\big) \big|+\mathbf{1} * \boldsymbol{\phi}_{\cdot j}$.
   These also yield that  $T_\infty=\infty$ and $\boldsymbol{V}\in C(\mathbb{R}_+;\mathbb{C}_-^d)$. 
   \qed

   By using Lemma~\ref{Lemma.FourierLaplaceHP} and the change of variables,  the next lemma provides an explicit exponential-affine representation of Fourier-Laplace functional of $\big(  \widetilde{\boldsymbol{N}}^{(n)}, \boldsymbol{\varLambda}^{(n)} \big)$ in term of the unique solution to the following two equivalent nonlinear Volterra equations
  \beqlb\label{eqn.504}
  \boldsymbol{V}^{(n)} =   \boldsymbol{W}^{(n)}   * \boldsymbol{R}^{(n)}
  \quad \mbox{and}\quad 
  \boldsymbol{V}^{(n)} = \big(  \boldsymbol{W}^{(n)} + \sqrt{n \theta_n}\boldsymbol{V}^{(n)} \big) * \boldsymbol{\phi}^{(n)}  
  \eeqlb
  for row vector functions $\boldsymbol{f}\in L^\infty_{\rm loc}(\mathbb{R}_+;\mathbb{C}_-^d)$ and $\boldsymbol{h}\in L^\infty_{\rm loc}(\mathbb{R}_+;\mathtt{i}\mathbb{R}^d)$, where
  \beqlb\label{eqn.5041}
  \boldsymbol{W}^{(n)}:= \boldsymbol{f} + n\theta_n\cdot \mathcal{W}\Big(\frac{\boldsymbol{V}^{(n)}  + \boldsymbol{h}}{\sqrt{n\cdot\theta_n} }  \Big)
  \quad \mbox{and}\quad 
  \boldsymbol{\phi}^{(n)} (t):= \sqrt{\frac{n}{\theta_n}}\cdot \boldsymbol{\phi}_n(nt),\quad t\geq 0.
  \eeqlb

 \begin{lemma}\label{Lemma.FLFLambdaNn}
 For row vector functions $\boldsymbol{f} \in L^\infty_{\rm loc}(\mathbb{R}_+;\mathbb{C}_-^d)$ and $\boldsymbol{h}\in L^\infty_{\rm loc}(\mathbb{R}_+;\mathtt{i}\mathbb{R}^d)$, we have the following.
  	
 \begin{enumerate}
 \item[(1)]  The  equivalent  Volterra equations in (\ref{eqn.504}) have a unique global solution  $ \boldsymbol{V}^{(n)}\in C(\mathbb{R}_+;\mathbb{C}^d_-)$.
  		
 \item[(2)] The Fourier-Laplace functional of  $\big( \boldsymbol{\varLambda}^{(n)} , \widetilde{\boldsymbol{N}}^{(n)}\big)$ admits the representation 
 \beqlb\label{eqn.503} 
 \mathbf{E} \Big[ \exp\Big\{\boldsymbol{f} *\boldsymbol{\varLambda}^{(n)}(T) + \boldsymbol{h}* d \widetilde{\boldsymbol{N}}^{(n)}(T)  \Big\} \Big]
 \ar=\ar  \exp \big\{   \boldsymbol{W}^{(n)}  * \boldsymbol{H}^{(n)}(T) \big\} 
 , \quad T\geq 0. 
 \eeqlb
  		
 \item[(3)] For each $T\geq 0$, there exist an integer $n_0\geq 1$ and a constant $C>0$  independent of $(\boldsymbol{f}, \boldsymbol{h})$ such that  
 \beqlb \label{eqn.506}  
 \sup_{n\geq n_0}\big\|\boldsymbol{V}^{(n)}\big\|_{L^\infty_T} 
 \ar\leq\ar C \cdot \big( \|\boldsymbol{f}\|_{L^\infty_T}  + \|\boldsymbol{h}\|_{L^\infty_T}^2 \big)  
 \quad \mbox{and} \quad 
 \sup_{n\geq n_0}\big\|\boldsymbol{W}^{(n)}\big\|_{L^\infty_T} 
  \leq  C \cdot \big( 1+ \|\boldsymbol{f}\|_{L^\infty_T}^2    + \|\boldsymbol{h}\|_{L^\infty_T}^4 \big) . \qquad 
 \eeqlb
 \end{enumerate}
 \end{lemma}
  \proof  
  For each $n\geq 1$, we define the following two functions
  \beqnn
  \boldsymbol{f}_n(t):= \frac{\boldsymbol{f}(t/n)}{n \cdot\theta_n}
  \quad\mbox{and}\quad 
  \boldsymbol{h}_n(t):= \frac{\boldsymbol{h}(t/n)}{\sqrt{n\cdot\theta_n}} ,
  \quad t\geq 0.
  \eeqnn
  Let $\boldsymbol{V}_n \in C(\mathbb{R}_+;\mathbb{C}^d_-)$ be the unique global solution to the  equivalent nonlinear Volterra integral equations in (\ref{VolRiccati}) with $(\boldsymbol{f}, \boldsymbol{h})= (\boldsymbol{f}_n, \boldsymbol{h}_n)$, i.e., 
  \beqlb\label{eqn.505}
  \boldsymbol{V}_n =  \boldsymbol{W}_n  * \boldsymbol{R}_n
  \quad \mbox{and}\quad 
  \boldsymbol{V}_n = ( \boldsymbol{W}_n+\boldsymbol{V}_n) * \boldsymbol{\phi}_n
  \quad \mbox{with}\quad 
  \boldsymbol{W}_n = \boldsymbol{f}_n  +\mathcal{W}(\boldsymbol{V}_n + \boldsymbol{h}_n) .
  \eeqlb
  Rescale the two functions $\boldsymbol{V}_n$ and $\boldsymbol{W}_n$ as follows
  \beqlb\label{eqn.515}
  \boldsymbol{V}^{(n)}(t) = \sqrt{n \theta_n} \cdot \boldsymbol{V}_n(nt)
  \quad \mbox{and}\quad 
  \boldsymbol{W}^{(n)}(t) = n \theta_n \cdot \boldsymbol{W}_n(nt) ,\quad t\geq 0. 
  \eeqlb
 Using the change of variable to (\ref{eqn.505}), one can identify that 
  \beqnn
  \boldsymbol{W}^{(n)} =\boldsymbol{f} + n\theta_n\cdot \mathcal{W}\Big(\frac{\boldsymbol{V}^{(n)}  + \boldsymbol{h}}{\sqrt{n\cdot\theta_n} } \Big)
  \eeqnn
 and $ \boldsymbol{V}^{(n)}$ solves the equivalent nonlinear Volterra integral equations in (\ref{eqn.504}). The uniqueness of solutions of (\ref{eqn.504}) follows directly from that of (\ref{eqn.505}). Claim (1) has been proved.
  
 Recall the function $ \boldsymbol{H}_n$ defined in (\ref{SVR.n}). For claim (2), by the change of variables and  Lemma~\ref{Lemma.FourierLaplaceHP}(2),
  \beqlb\label{eqn.502} 
  \mathbf{E} \Big[ \exp\Big\{\boldsymbol{f} *\boldsymbol{\varLambda}^{(n)}(T) + \boldsymbol{h}* d \widetilde{\boldsymbol{N}}^{(n)}(T)  \Big\} \Big]
  \ar=\ar\mathbf{E} \Big[ \exp\Big\{\boldsymbol{f}_n *\boldsymbol{\varLambda}_n(nT) + \boldsymbol{h}_n* d \widetilde{\boldsymbol{N}}_n(nT)  \Big\} \Big]\cr
  \ar\ar\cr
  \ar=\ar \exp \Big\{   \boldsymbol{W}_n  * \boldsymbol{H}_n(nT) \Big\}\cr
  \ar=\ar \exp \Big\{   \boldsymbol{W}^{(n)}  * \boldsymbol{H}^{(n)}(T) \Big\} . 
  \eeqlb
  Here the last equality follows by using the change of variables along with (\ref{IntegralLambda}) and (\ref{eqn.515}).
  
  We now start to prove claim (3). 
  For $\lambda >\lambda_+$ and $t\geq 0$, let 
   \beqlb\label{eqn.511}
  \boldsymbol{V}^{(n)}_\lambda(t):= e^{-\lambda t}  \boldsymbol{V}^{(n)}(t), \quad 
  \boldsymbol{f}_\lambda(t):=e^{-\lambda t}  \boldsymbol{f}(t),\quad
  \boldsymbol{W}^{(n)}_\lambda(t):= e^{-\lambda t}  \boldsymbol{W}^{(n)}(t),\quad 
  \boldsymbol{\phi}^{(n)}_\lambda(t):= e^{-\lambda t}  \boldsymbol{\phi}^{(n)}(t). 
  \eeqlb  
  Multiplying both sides of the second nonlinear Volterra equation in (\ref{eqn.504}) by the exponential function $e^{-\lambda t}$, we can obtain that 
  \beqnn
   \boldsymbol{V}^{(n)}_\lambda(t)
  \ar=\ar \Big(   \boldsymbol{W}^{(n)}_\lambda +  \sqrt{n \theta_n}\boldsymbol{V}^{(n)}_\lambda \Big) * \boldsymbol{\phi}^{(n)}_\lambda(t),\quad t\geq 0.
  \eeqnn
  By (\ref{eqn.505}), the right-hand side can be separated into the following three parts: $\boldsymbol{I}^{(n)}_1(t) = \boldsymbol{f}_\lambda  * \boldsymbol{\phi}^{(n)}_\lambda(t) $ and
  \beqnn
    \boldsymbol{I}^{(n)}_2(t) \ar=\ar  \int_0^t n\theta_n\bigg(\exp\Big\{ \frac{ \boldsymbol{h}(s)}{\sqrt{n\theta_n}}   \Big\}-1 -\frac{ \boldsymbol{h}(s)}{\sqrt{n\theta_n}}  \bigg)e^{-\lambda s} \cdot \boldsymbol{\phi}^{(n)}_\lambda(t-s)ds, \cr
    \boldsymbol{I}^{(n)}_3(t) \ar=\ar \int_0^t n\theta_n\bigg(\exp\Big\{ \frac{\boldsymbol{V}^{(n)}(s) }{\sqrt{n\theta_n}} +\frac{ \boldsymbol{h}(s)}{\sqrt{n\theta_n}}   \Big\}  -\exp\Big\{ \frac{ \boldsymbol{h}(s)}{\sqrt{n\theta_n}}   \Big\} \bigg)e^{-\lambda s} \cdot \boldsymbol{\phi}^{(n)}_\lambda(t-s)ds.
  \eeqnn
  Firstly, we have $ \| \boldsymbol{f}_\lambda  * \boldsymbol{\phi}^{(n)}_\lambda \|_{L^\infty_T} \leq   \|\boldsymbol{f}_\lambda\|_{L^\infty_T}  \cdot \| \boldsymbol{\phi}^{(n)}_\lambda \|_{L^1} $. 
  Secondly, by the inequality that $ |e^x-1-x| \leq |x|^2$ for any $x\in \mathbb{C}_-$,  
  \beqnn
   \big\|\boldsymbol{I}^{(n)}_2 \big\|_{L^\infty_T} 
   \ar\leq\ar n\theta_n \cdot \bigg\| \exp\Big\{ \frac{ \boldsymbol{h}}{\sqrt{n\theta_n}} \Big\}-1 -\frac{ \boldsymbol{h}}{\sqrt{n\theta_n}} \bigg\|_{L^\infty_T} \cdot \| \boldsymbol{\phi}^{(n)}_\lambda\|_{L^1}
   \leq   \|\boldsymbol{h}\|_{L^\infty_T}^2  \cdot \| \boldsymbol{\phi}^{(n)}_\lambda\|_{L^1}.
  \eeqnn
  Thirdly, by using the equality  $ |e^{x+y}-e^x| \leq |y| $ for any $x,y\in \mathbb{C}_-$ and then the change of variables,  
  \beqnn
  \big\|\boldsymbol{I}^{(n)}_3 \big\|_{L^\infty_T} \leq \big\|\boldsymbol{V}^{(n)}_\lambda \big\|_{L^\infty_T}\cdot \sqrt{n\theta_n} \big\| \boldsymbol{\phi}^{(n)}_\lambda \big\|_{L^1} = \big\|\boldsymbol{V}^{(n)}_\lambda \big\|_{L^\infty_T}\cdot  \int_0^\infty e^{-\frac{\lambda }{n}\cdot t}\boldsymbol{\phi}_{n}(t)dt 
  = \big\|\boldsymbol{V}^{(n)}_\lambda \big\|_{L^\infty_T}\cdot \mathcal{L}_{\boldsymbol{\phi}_{n}}(\lambda /n). 
  \eeqnn  
  Putting all estimates above together,  we have 
  \beqnn
  \big\| \boldsymbol{V}^{(n)}_\lambda \big\|_{L^\infty_T} 
  \ar\leq\ar \big( \|\boldsymbol{f}\|_{L^\infty_T}  +\|\boldsymbol{h}\|_{L^\infty_T}^2 \big)\cdot \big\| \boldsymbol{\phi}^{(n)}_\lambda\big\|_{L^1} +\big\|\boldsymbol{V}^{(n)}_\lambda\big\|_{L^\infty_T}\cdot \mathcal{L}_{\boldsymbol{\phi}_{n}}(\lambda /n)
  \eeqnn
  Recall the function $\boldsymbol{\varPsi}^{(n)}$ defined in (\ref{eqn.varPsin}). Moving the last term in the preceding inequality to the left hand side and then multiplying both sides by $\sqrt{n\theta_n}$, 
  \beqlb\label{eqn.507}
  \big\| \boldsymbol{V}^{(n)}_\lambda \big\|_{L^\infty_T} \cdot \boldsymbol{\varPsi}^{(n)}(\lambda) 
  \ar\leq\ar 
  \big(\|\boldsymbol{f}\|_{L^\infty_T} +\|\boldsymbol{h}\|_{L^\infty_T}^2 \big)\cdot \sqrt{n\theta_n} \big\| \boldsymbol{\phi}^{(n)}_\lambda\big\|_{L^1}
  = \big( \|\boldsymbol{f}\|_{L^\infty_T}  +\|\boldsymbol{h}\|_{L^\infty_T}^2 \big)\cdot \mathcal{L}_{\boldsymbol{\phi}_{n}}(\lambda /n).
  \eeqlb
  Since $(\boldsymbol{\varPsi}^{(n)}(\lambda))^{-1} \to \mathcal{L}_{\boldsymbol{\varPi}}(\lambda)\in\mathbb{R}_+^{d\times d}$ as $n\to\infty$; see (\ref{eqn.LimInveversePsi}), there exist an integer $n_0\geq 1$ and a constant $C>0$ independent of $\boldsymbol{f}$ and $\boldsymbol{h}$ such that for any $n\geq n_0$,
  \beqnn
  \boldsymbol{\varPsi}^{(n)}(\lambda)
  \mbox{ is invertible} \mbox{ and }  
  \big(\boldsymbol{\varPsi}^{(n)}(\lambda)\big)^{-1}\leq C .
  \eeqnn 
  Multiplying both sides of (\ref{eqn.507}) by $\big(\boldsymbol{\varPsi}^{(n)}(\lambda)\big)^{-1}$ and then using the preceding upper bound as well as Condition~\ref{Main.Condition}, there exists a constant $C>0$ independent of $(\boldsymbol{f},\boldsymbol{h})$  such that for any $n\geq n_0$,
  \beqnn
  \big\| \boldsymbol{V}^{(n)}_\lambda \big\|_{L^\infty_T} 
  \leq \big( \|\boldsymbol{f}\|_{L^\infty_T}  +\|\boldsymbol{h}\|_{L^\infty_T}^2 \big)\cdot \mathcal{L}_{\boldsymbol{\phi}_{n}}(\lambda /n) \cdot \big(\boldsymbol{\varPsi}^{(n)}(\lambda)\big)^{-1}  \leq C\cdot  \big( \|\boldsymbol{f}\|_{L^\infty_T}  +\|\boldsymbol{h}\|_{L^\infty_T}^2 \big) . 
  \eeqnn
  This immediately yields  that 
  \beqnn
  \big\| \boldsymbol{V}^{(n)}  \big\|_{L^\infty_T} \leq e^{\lambda t} \cdot\| \boldsymbol{V}^{(n)}_\lambda \|_{L^\infty_T} \leq C\cdot e^{\lambda t} \big( \|\boldsymbol{f}\|_{L^\infty_T}  +\|\boldsymbol{h}\|_{L^\infty_T}^2 \big).
  \eeqnn 
  By using the inequality that $ |e^x-1-x| \leq |x|^2 $ for any $x\in \mathbb{C}_-$, the triangle inequality  and then the preceding upper bound estimate to $ \boldsymbol{W}^{(n)}$, 
  \beqnn
  \big\|\boldsymbol{W}^{(n)} \big\|_{L^\infty_T}  
  \ar\leq\ar \|\boldsymbol{f}\|_{L^\infty_T} + \big( \big\| \boldsymbol{V}^{(n)}\big\|_{L^\infty_T} + \|\boldsymbol{h} \|_{L^\infty_T}  \big)^2 
  \leq  C \cdot \big(\|\boldsymbol{f}\|_{L^\infty_T}  +\|\boldsymbol{f}\|_{L^\infty_T}^2 + \|\boldsymbol{h} \|_{L^\infty_T}^2 +\|\boldsymbol{h}\|_{L^\infty_T}^4 \big) .
  \eeqnn
  We have finished the whole proof. 
  \qed
  
  In the sequel of this section, we consider the well-posedness of (\ref{Eqn.AffineVolterra}) by studying the asymptotics of the sequence $\{ \boldsymbol{V}^{(n)} \}_{n\geq 1}$ and then establish the exponential-affine representation of Fourier-Laplace functional of solutions to (\ref{Eqn.HawkesVoletrra}). 
  
  \begin{lemma} \label{Lemma.UniquenessRV}
  	For row vector functions $\boldsymbol{f} \in L^\infty_{\rm loc}(\mathbb{R}_+;\mathbb{C}_-^d)$ and $\boldsymbol{h}\in L^\infty_{\rm loc}(\mathbb{R}_+;\mathtt{i}\mathbb{R}^d)$, the uniqueness of global solutions in $L^\infty_{\rm loc}(\mathbb{R}_+;\mathbb{C}_-^d)$ holds for (\ref{Eqn.AffineVolterra}).
  	
  \end{lemma}
  \proof  
  Assume that $  \boldsymbol{\mathcal{V}}, \boldsymbol{\mathcal{V}}^* \in L^\infty_{\rm loc}(\mathbb{R}_+;\mathbb{C}^d_-)$ are two global solutions of  (\ref{Eqn.AffineVolterra}). 
  We first prove the case of $\boldsymbol{h}=0$.
  By Assumption~\ref{MainAssumption02}, we have for $i\in\mathtt{D}$, 
  \beqlb \label{eqn.3004}
  \frac{\underline{\varPi}_{ii}(0)}{2}\cdot (\mathcal{V}_{i}^*(t))^2-\mathcal{V}_{i}^*(t)-\Big(  \frac{\underline{\varPi}_{ii}(0)}{2}\cdot (\mathcal{V}_{i}(t))^2-\mathcal{V}_{i}(t) \Big)  =  \frac{1}{2}\sum_{j=1}^d \int_{(0,t]}\big( (\mathcal{V}_j)^2- (\mathcal{V}_j^* )^2 \big)(t-s) \varPi_{ji}(ds).
  \eeqlb
  Since $\mathcal{V}_i,\mathcal{V}_i^* \in L^\infty_{\rm loc}(\mathbb{R}_+;\mathbb{C}_-)$,  there exists a constant $C_0>0$ independent of $i,j\in\mathtt{D}$ such that uniformly in $t\in[0,T]$,
  \beqlb\label{eqn.6003}
  \int_{(0,t]} \Big|\big(\mathcal{V}_j(t-s)\big)^2 - \big(\mathcal{V}_j^*(t-s)\big)^2 \Big| \varPi_{ji}(ds)
  \ar\leq\ar  C_0      \int_{(0,t]} \big|\mathcal{V}_j(t-s)-  \mathcal{V}_j^*(t-s)\big| \varPi_{ji}(ds).
  \eeqlb
  For $i\in \mathtt{D}$ with $\underline{\varPi}_{ii}(0)=0$, taking (\ref{eqn.6003}) back into (\ref{eqn.3004}) we have
  \beqlb\label{eqn.6005}
  \big| \mathcal{V}_i(t)-  \mathcal{V}_i^*(t) \big|
  \leq \frac{C_0}{2}\cdot  \sum_{j=1}^d \int_{(0,t]} \big|\mathcal{V}_j(t-s)-  \mathcal{V}_j^*(t-s)\big| \varPi_{ji}(ds).
  \eeqlb 
  On the other hand,  for $i\in \mathtt{D}$ with $\underline{\varPi}_{ii}(0)>0$, we can write the left-hand side of (\ref{eqn.3004}) as
  \beqnn
  \frac{\underline{\varPi}_{ii}(0)}{2}\cdot \Big(\mathcal{V}_{i}^*(t)+\mathcal{V}_{i}(t)-2\cdot\big|\underline{\varPi}_{ii}(0)\big|^{-1}\Big) \cdot \Big(\mathcal{V}_{i}^*(t)- \mathcal{V}_{i}(t) \Big).
  \eeqnn
  Since $\mathtt{Re}\mathcal{V}_i,\mathtt{Re}\mathcal{V}_i^* \leq 0$, we have 
  \beqnn
  \bigg|\frac{\underline{\varPi}_{ii}(0)}{2}\cdot (\mathcal{V}_{i}^*(t))^2-\mathcal{V}_{i}^*(t)-\Big(  \frac{\underline{\varPi}_{ii}(0)}{2}\cdot (\mathcal{V}_{i}(t))^2-\mathcal{V}_{i}(t) \Big)\bigg|
  \geq \big|\mathcal{V}_{i}^*(t)- \mathcal{V}_{i}(t) \big|  
  \eeqnn
  and by (\ref{eqn.6003}),
  \beqnn
  \big|\mathcal{V}_{i}^*(t)- \mathcal{V}_{i}(t) \big|
  \leq \frac{1}{2} \sum_{j=1}^d \int_{(0,t]} \big|\mathcal{V}_j(t-s)-  \mathcal{V}_j^*(t-s)\big| \varPi_{ji}(ds).
  \eeqnn
  Putting this together with (\ref{eqn.6005}), one can find a constant $C_1>0$ depend on $T$ such that 
  \beqlb\label{eqn.6006}
  \big|   \boldsymbol{\mathcal{V}}(t)-\boldsymbol{\mathcal{V}}^*(t)| \leq 
  \frac{C_1}{2} \cdot \int_{(0,t]} \big|\boldsymbol{\mathcal{V}}(t-s)-\boldsymbol{\mathcal{V}}^*(t-s)\big|   \boldsymbol{\varPi} (ds) , 
  \quad t\in[0,T]. 
  \eeqlb
  For some constant $\lambda>0$ to be specified later, let 
  \beqlb \label{eqn.6007}
  \boldsymbol{\mathcal{V}}_\lambda(t):= e^{-\lambda t}\cdot \boldsymbol{\mathcal{V}}(t),\quad
  \boldsymbol{\mathcal{V}}^*(t):=  e^{-\lambda t}\cdot \boldsymbol{\mathcal{V}}^*(t),\quad  \boldsymbol{\varPi}_\lambda(dt):=  e^{-\lambda t}\cdot \boldsymbol{\varPi}(dt),\quad t\geq 0. 
  \eeqlb
  Multiplying both sides of (\ref{eqn.6006}) by $e^{-\lambda t}$, we have 
  \beqnn
  \big|   \boldsymbol{\mathcal{V}}_\lambda(t)-\boldsymbol{\mathcal{V}}^*_\lambda(t)| \leq 
  \frac{C_1}{2} \cdot \int_{(0,t]} \big|\boldsymbol{\mathcal{V}}_\lambda(t-s)-\boldsymbol{\mathcal{V}}^*_\lambda(t-s)\big|   \boldsymbol{\varPi}_\lambda (ds) , \quad   t\in[0,T].
  \eeqnn 
  By the fact that $\|\boldsymbol{\mathcal{V}}\|_{L^\infty_T},\|\boldsymbol{\mathcal{V}}^*\|_{L^\infty_T} <\infty$, there exists a constant $C_2>0$ independent of $\lambda $ such that
  \beqnn
  \big\|   \boldsymbol{\mathcal{V}}_\lambda-\boldsymbol{\mathcal{V}}^*_\lambda  \big\|_{L^\infty_T}  
  \leq \frac{C_2}{2} \cdot \big\| \boldsymbol{\mathcal{V}}_\lambda-\boldsymbol{\mathcal{V}}^*_\lambda\big\|_{L^\infty_T}\cdot  \boldsymbol{\varPi}_\lambda\big((0,T]\big).
  \eeqnn
  Notice  that  $\boldsymbol{\varPi}_\lambda\big((0,T]\big) \to \mathbf{0}$ as $\lambda \to \infty$, we choose $\lambda >0$ large enough such that $\boldsymbol{\varPi}_\lambda\big((0,T]\big)<1/C_2$ and then  $\|\boldsymbol{\mathcal{V}}_\lambda-\boldsymbol{\mathcal{V}}^*_\lambda  \|_{L^\infty_T} =\|\boldsymbol{\mathcal{V}}-\boldsymbol{\mathcal{V}}^* \|_{L^\infty_T} =0$.
  This induces the uniqueness for  (\ref{Eqn.AffineVolterra}).   
  For the general case of $\boldsymbol{h}\neq\mathbf{0}$, let $\boldsymbol{\mathcal{V}}^{\boldsymbol{h}}:=\boldsymbol{\mathcal{V}} +\boldsymbol{h}$ and $\boldsymbol{\mathcal{V}}^{*,\boldsymbol{h}}:= \boldsymbol{\mathcal{V}}^*+ \boldsymbol{h}$, which satisfy
  \beqnn
  \boldsymbol{\mathcal{V}}^{\boldsymbol{h}}= \boldsymbol{f} +\boldsymbol{h}+ \frac{1}{2}\cdot \big(  \boldsymbol{\mathcal{V}}^{\boldsymbol{h}} \big)^2 *\boldsymbol{\varPi}
  \quad \mbox{and}\quad 
  \boldsymbol{\mathcal{V}}^{*,\boldsymbol{h}}=  \boldsymbol{f} +\boldsymbol{h}+ \frac{1}{2}\cdot\big(  \boldsymbol{\mathcal{V}}^{*,\boldsymbol{h}} \big)^2 *\boldsymbol{\varPi}. 
  \eeqnn 
  Thus $\boldsymbol{\mathcal{V}}= \boldsymbol{\mathcal{V}}^{*}$ if and only if $\boldsymbol{\mathcal{V}}^{\boldsymbol{h}}=\boldsymbol{\mathcal{V}}^{*,\boldsymbol{h}}$, which follows directly from the previous result.
  \qed 
  
  \begin{lemma} \label{Lemma.ExitenceRV}
  	For a subsequence $\{ \boldsymbol{V}^{(n_k)} \}_{k\geq 1}$, if  $\boldsymbol{V}^{(n_k)} \overset{\rm a.e.}\to \boldsymbol{\mathcal{V}}$ as $k\to\infty$, then the following hold.
  	\begin{enumerate}
  		\item[(1)]	$\boldsymbol{\mathcal{V}}  \in L^\infty_{\rm loc}(\mathbb{R}_+;\mathbb{C}_-^d)$ and is the unique global solution of (\ref{Eqn.AffineVolterra}). 
  		
  		\item[(2)] For each $T\geq 0$, we have $\|\boldsymbol{W}^{(n_k)}-\boldsymbol{\mathcal{W}}\|_{L^\infty_T}\to 0$ as $k\to\infty$. 
  	\end{enumerate}

  \end{lemma}
  \proof  By Lemma~\ref{Lemma.FLFLambdaNn}(3), we have $\boldsymbol{\mathcal{V}}  \in L^\infty_{\rm loc}(\mathbb{R}_+;\mathbb{C}_-^d)$ and $ \big\| \boldsymbol{V}^{(n_k)} -\boldsymbol{\mathcal{V}}  \big\|_{L^\infty_T}\to 0$ as $n\to\infty $ for any $T\geq 0$.   
   By using the inequality $|e^x-1-x-x^2/2|\leq |x|^2$ for any $x\in \mathbb{C}_-$,
   \beqlb
   \big\|\boldsymbol{W}^{(n_k)}-\boldsymbol{\mathcal{W}} \big\|_{L^\infty_T}
   \ar=\ar \bigg\| n\theta_n\cdot \mathcal{W}\Big(\frac{\boldsymbol{V}^{(n)} + \boldsymbol{h}}{\sqrt{n\theta_n} } \Big) - \frac{1}{2} (\boldsymbol{V}^{(n)} + \boldsymbol{h})^2 \bigg\|_{L^\infty_T} +\frac{1}{2}\big\|  (\boldsymbol{V}^{(n)} + \boldsymbol{h})^2- (\boldsymbol{\mathcal{V}} + \boldsymbol{h})^2 \big\|_{L^\infty_T} \cr
   \ar\leq\ar \frac{\big\|\boldsymbol{V}^{(n)} + \boldsymbol{h} \big\|_{L^\infty_T}^3}{\sqrt{n\cdot \theta_n}}  + \frac{1}{2}\big\|  \boldsymbol{V}^{(n)}  + \boldsymbol{\mathcal{V}}+ 2 \boldsymbol{h} \big\|_{L^\infty_T} \cdot \big\|  \boldsymbol{V}^{(n)} - \boldsymbol{\mathcal{V}}  \big\|_{L^\infty_T},
   \eeqlb
   which goes to $0$ as $n\to\infty$ and hence claim (2) holds. 
   
  We now prove that $\boldsymbol{\mathcal{V}}  $ solves (\ref{Eqn.AffineVolterra}), which holds if along a subsequence  $\{ n_l \}_{l\geq 1}$ of $\{ n_k \}_{k\geq 1}$,
  \beqlb\label{eqn.4013}
  \boldsymbol{f} * \boldsymbol{R}^{(n_l)} \overset{\rm a.e.}\to  \boldsymbol{f} * \boldsymbol{\varPi} 
  \quad \mbox{and}\quad 
   n\theta_n\cdot \mathcal{W}\Big(\frac{\boldsymbol{V}^{(n_l)} + \boldsymbol{h}}{\sqrt{n\theta_n} } \Big) * \boldsymbol{R}^{(n_l)} \overset{\rm a.e.}\to  \frac{1}{2}\cdot \big(  \boldsymbol{\mathcal{V}} +\boldsymbol{h}\big)^2 * \boldsymbol{\varPi},
  \eeqlb
  as $l\to\infty$. 
  A simple transform shows that
   \beqnn
  n\theta_n\cdot \mathcal{W}\Big(\frac{\boldsymbol{V}^{(n)} + \boldsymbol{h}}{\sqrt{n\theta_n} } \Big) * \boldsymbol{R}^{(n)} - \frac{1}{2}\cdot \big(  \boldsymbol{\mathcal{V}} +\boldsymbol{h}\big)^2 * \boldsymbol{\varPi}
  \ar=\ar \bigg( n\theta_n\cdot \mathcal{W}\Big(\frac{\boldsymbol{V}^{(n)} + \boldsymbol{h}}{\sqrt{n\theta_n} } \Big) - \frac{1}{2} (\boldsymbol{V}^{(n)} + \boldsymbol{h})^2 \bigg)* \boldsymbol{R}^{(n)} \cr
  \ar\ar +  \frac{1}{2} \cdot \Big(  (\boldsymbol{V}^{(n)} + \boldsymbol{h})^2-   \big(  \boldsymbol{\mathcal{V}} +\boldsymbol{h}\big)^2 \Big)* \boldsymbol{R}^{(n)}\cr
  \ar\ar 
  + \frac{1}{2}\cdot \big( \big(  \boldsymbol{\mathcal{V}} +\boldsymbol{h}\big)^2 *  \boldsymbol{R}^{(n)} - \big(\boldsymbol{\mathcal{V}} +\boldsymbol{h}\big)^2 * \boldsymbol{\varPi} \Big).
  \eeqnn
   By using the inequality $|e^x-1-x-x^2/2|\leq |x|^2$ for any $x\in \mathbb{C}_-$ along with (\ref{Eqn.BoundConvIRn}) and (\ref{eqn.506}),
  \beqnn
  \bigg| n\theta_n\cdot \mathcal{W}\Big(\frac{\boldsymbol{V}^{(n)} + \boldsymbol{h}}{\sqrt{n\theta_n} } \Big) - \frac{1}{2} (\boldsymbol{V}^{(n)} + \boldsymbol{h})^2 \bigg|* \boldsymbol{R}^{(n)} (T)
  \leq \frac{1}{\sqrt{n\cdot \theta_n}} \big\|\boldsymbol{V}^{(n)} + \boldsymbol{h} \big\|_{L^\infty_T}^3 \cdot   \mathcal{I}_{\boldsymbol{R}^{(n)} }(T) \to 0 
  \eeqnn
  and also 
  \beqnn
  \Big| (\boldsymbol{V}^{(n)} + \boldsymbol{h})^2-   \big(  \boldsymbol{\mathcal{V}} +\boldsymbol{h}\big)^2 \Big|* \boldsymbol{R}^{(n)}(T) 
  \ar\leq\ar C\cdot \big\| \boldsymbol{V}^{(n)} -\boldsymbol{\mathcal{V}}  \big\|_{L^\infty_T}\cdot \mathcal{I}_{\boldsymbol{R}^{(n)} }(T) \to 0,  
  \eeqnn
  as $n\to\infty$ for any $T\geq 0$. 
  Moreover, by Corollary~\ref{Prop.503} we also have 
  \beqnn
  \big\|  \boldsymbol{f} * \boldsymbol{R}^{(n)} -  \boldsymbol{f} * \boldsymbol{\varPi} \big\|_{L^1_T} \to 0
  \quad \mbox{and}\quad
   \big\| \big(  \boldsymbol{\mathcal{V}} +\boldsymbol{h}\big)^2 *  \boldsymbol{R}^{(n)} - \big(\boldsymbol{\mathcal{V}} +\boldsymbol{h}\big)^2 * \boldsymbol{\varPi}  \big\|_{L^1_T} \to 0,
  \eeqnn
  as $n\to\infty$ for any $T\geq 0$. Consequently,  we have that along a subsequence $\{n_l\}_{l\geq1}$ of $\{n_k\}_{k\geq1}$, 
  \beqnn
   \boldsymbol{f} * \boldsymbol{R}^{(n_l)} \overset{\rm a.e.}\to  \boldsymbol{f} * \boldsymbol{\varPi} 
   \quad \mbox{and} \quad
  \big(  \boldsymbol{\mathcal{V}} +\boldsymbol{h}\big)^2 *  \boldsymbol{R}^{(n_l)} \overset{\rm a.e.}\to \big(\boldsymbol{\mathcal{V}} +\boldsymbol{h}\big)^2 * \boldsymbol{\varPi},
  \eeqnn
    as $l\to\infty$.  The two limits in (\ref{eqn.4013}) hold.  
    The uniqueness follows from Lemma~\ref{Lemma.UniquenessRV}. 
   \qed

 \begin{lemma}\label{Lemma.ConV}
 For row vector functions $\boldsymbol{f}\in L^\infty_{\rm loc}(\mathbb{R}_+;\mathbb{C}_-^d)$ and $\boldsymbol{h}\in L^\infty_{\rm loc}(\mathbb{R}_+;\mathtt{i}\mathbb{R}^d)$,
  assume that $(T_\infty,\boldsymbol{\mathcal{V}}) \in (0,\infty] \times L^\infty_{\rm loc}([0,T_\infty);\mathbb{C}^d)$  is a non-continuable solution of  (\ref{Eqn.AffineVolterra}). If $\underline{\boldsymbol{\varPi}}(0)=\mathbf{0}$, we have the following.    
  \begin{enumerate}
  	\item[(1)] The function $ \boldsymbol{\mathcal{V}} \in L^\infty_{\rm loc}([0,T_\infty);\mathbb{C}^d_-) $ is the unique global solution of (\ref{Eqn.AffineVolterra}).

  	\item[(2)]  Let $\boldsymbol{V}^{(n)} \in C(\mathbb{R}_+;\mathbb{C}^d_-)$ be the unique solution of (\ref{eqn.504}). 
  	Recall $ \boldsymbol{W}^{(n)} $ and $ \boldsymbol{\mathcal{W}}$ defined in (\ref{eqn.5041}) and (\ref{Eqn.AffineVolterra}).
  	For  each $T\geq 0$, we have
  	\beqlb\label{eqn.509}
   \lim_{n\to\infty}\Big( \big\| \boldsymbol{V}^{(n)} -\boldsymbol{\mathcal{V}} \big\|_{L^\infty_T} +\big\| \boldsymbol{W}^{(n)} -\boldsymbol{\mathcal{W}} \big\|_{L^\infty_T} \Big) = 0 
  	\eeqlb
 and there exists a constant $C>0$ independent of $\boldsymbol{f}$ and $\boldsymbol{h}$ such that 
  	\beqlb\label{eqn.5010}
  	\| \boldsymbol{\mathcal{V}} \|_{L^\infty_T} \leq  C \cdot\big( \|\boldsymbol{f}\|_{L^\infty_T} + \|\boldsymbol{h}\|_{L^\infty_T}^2 \big) 
  	\quad \mbox{and}\quad
  	\big\|\boldsymbol{\mathcal{W}} \big\|_{L^\infty_T}
  	\leq C\cdot  \big(1+ \|\boldsymbol{f}\|_{L^\infty_T} + \|\boldsymbol{h}\|_{L^\infty_T}^2 \big)^2 . 
  	\eeqlb
  \end{enumerate} 
 \end{lemma}
 \proof   
 It is obvious that claim (1) follows directly from claim (2) and Lemma~\ref{Lemma.FLFLambdaNn}. 
 We now start to prove claim (2). 
 For a constant $\lambda \geq \lambda_+$ to be specified later, recall the functions $\boldsymbol{V}^{(n)}_\lambda$,  $\boldsymbol{\phi}^{(n)}_\lambda$ defined in (\ref{eqn.511}) and $\boldsymbol{\mathcal{V}}_\lambda$, $\boldsymbol{\varPi}_\lambda$ defined in (\ref{eqn.6007}). 
  Multiplying both sides of   (\ref{Eqn.AffineVolterra}) and (\ref{eqn.504}) by $e^{-\lambda t}$,
 \beqlb\label{eqn.3005}
 \boldsymbol{V}^{(n)}_\lambda  -\boldsymbol{\mathcal{V}}_\lambda   
 =\sum_{i=1}^4 \boldsymbol{J}_i^{(n)} 
 \quad \mbox{and hence}\quad 
 \big\|\boldsymbol{V}^{(n)}_\lambda  -\boldsymbol{\mathcal{V}}_\lambda\big\|_{L^1_T}   
 \leq \sum_{i=1}^4 \big\|\boldsymbol{J}_i^{(n)} \big\|_{L^1_T},
 \eeqlb 
 where $\boldsymbol{J}_1^{(n)}(t):= \boldsymbol{f}_\lambda * \boldsymbol{R}^{(n)}_\lambda (t)  -\boldsymbol{f}_\lambda * \boldsymbol{\varPi}_\lambda (t) $ and 
 \beqnn
 \boldsymbol{J}_2^{(n)}(t)\ar:=\ar \int_{[0,t]} \Big(n\theta_n\cdot \mathcal{W}\Big(\frac{\boldsymbol{V}^{(n)}(t-s)+\boldsymbol{h}(t-s)}{\sqrt{n\theta_n}}\Big)- \frac{1}{2}   \big(  \boldsymbol{V}^{(n)}(t-s)+\boldsymbol{h}(t-s) \big)^2 \Big) e^{-\lambda(t-s)}\cdot \boldsymbol{R}^{(n)}_\lambda (s)ds,\cr
 \boldsymbol{J}_3^{(n)}(t)\ar:=\ar\frac{1}{2}   \int_{[0,t]}  \big(  \boldsymbol{\mathcal{V}}(t-s)+\boldsymbol{h}(t-s)\big)^2e^{-\beta(t-s)} \cdot  \big( \boldsymbol{R}^{(n)}_\lambda(s)ds - \boldsymbol{\varPi}_\lambda(ds) \big),\cr
  \boldsymbol{J}_4^{(n)}(t)\ar:=\ar \frac{1}{2}   \int_{[0,t]} \Big( \big( \boldsymbol{V}^{(n)}(t-s)+\boldsymbol{h}(t-s)\big)^2- \big(\boldsymbol{\mathcal{V}}(t-s)+\boldsymbol{h}(t-s) \big)^2\Big)e^{-\lambda(t-s)}\cdot  \boldsymbol{R}^{(n)}_\lambda (s)ds . 
 \eeqnn 
 Moreover, by using Proposition~\ref{Prop.A1}(2) along with
 $
 \sup_{n\geq 1} \big\|\boldsymbol{V}^{(n)}\big\|_{L^\infty_T}+ \big\|\boldsymbol{\mathcal{V}}\big\|_{L^\infty_T} +\big\|\boldsymbol{h}\big\|_{L^\infty_T} <\infty 
 $, there exists a constant $C_0>0$ independent of $\lambda$ such that 
 \beqnn
 \big\|\boldsymbol{J}_4^{(n)} \big\|_{L^1_T}
 \ar\leq\ar \frac{1}{2} \cdot  \big\|\boldsymbol{V}^{(n)}  +\boldsymbol{\mathcal{V}} +\boldsymbol{h}\big\|_{L^1_T} \cdot  \big\|  | \boldsymbol{V}^{(n)}_\lambda   -  \boldsymbol{\mathcal{V}}_{\lambda}   | *    \boldsymbol{R}^{(n)}_\lambda \big\|_{L^1_T} 
 \leq \frac{C_0}{2}\cdot \big\|  \boldsymbol{V}^{(n)}_\lambda   -  \boldsymbol{\mathcal{V}}_{\lambda} \big\|_{L^1_T}\cdot \big\|      \boldsymbol{R}^{(n)}_\lambda \big\|_{L^1_T}.
 \eeqnn
 Notice that $\big\|      \boldsymbol{R}^{(n)}_\lambda \big\|_{L^1_T} = \mathcal{I}_{\boldsymbol{R}^{(n)}_\lambda}(T) \to \underline{\boldsymbol{\varPi}_{\lambda}}(T)$ as $n\to\infty$ and also $\underline{\boldsymbol{\varPi}_{\lambda}}(T)\to 0$ as $\lambda \to \infty$. 
 There exist two constant $n_0\geq 1$ and $\lambda >0$ such that  
 \beqnn
 \sup_{n\geq n_0} \big\|      \boldsymbol{R}^{(n)}_\lambda \big\|_{L^1_T} \leq 1/C_0
 \quad \mbox{and}\quad 
 \big\|\boldsymbol{J}_4^{(n)} \big\|_{L^1_T} \leq \frac{1}{2}\cdot  \big\|  \boldsymbol{V}^{(n)}_\lambda   -  \boldsymbol{\mathcal{V}}_{\lambda} \big\|_{L^1_T}.
 \eeqnn
 Putting all preceding estimates back into (\ref{eqn.3005}),  we have 
 \beqlb \label{eqn.3006}
 \big\|\boldsymbol{V}^{(n)}_\lambda  -\boldsymbol{\mathcal{V}}_\lambda\big\|_{L^1_T}
 \leq 2 \cdot  \sum_{i=1}^3 \big\|\boldsymbol{J}_i^{(n)} \big\|_{L^1_T}.
 \eeqlb
 In view of Corollary~\ref{Prop.503} and Remark~\ref{Remark.3.5},  we have as $n\to\infty$,
 \beqnn
 \big\| \boldsymbol{J}_1^{(n)} \big\|_{L^1_T} 
 \leq \big\| \boldsymbol{f} * \boldsymbol{R}^{(n)} -\boldsymbol{f}* \boldsymbol{\varPi}  \big\|_{L^1_T} \to 0
 \quad\mbox{and}\quad
 \big\| \boldsymbol{J}_3^{(n)} \big\|_{L^1_T} \leq  \big\| (\boldsymbol{\mathcal{V}}+\boldsymbol{h})^2* \boldsymbol{R}^{(n)} - (\boldsymbol{\mathcal{V}}+\boldsymbol{h})^2* \boldsymbol{\varPi}  \big\|_{L^1_T}  \to 0.
 \eeqnn 
  By using the inequality $|e^x-1-x-x^2/2| \leq |x|^3$ for any $x\in\mathbb{C}_-$ and the using  (\ref{Eqn.BoundConvIRn}) as well as
 $
 \sup_{n\geq 1} \big\|\boldsymbol{V}^{(n)}+\boldsymbol{h}\big\|_{L^\infty_T}  <\infty 
 $, we also have as $n\to\infty$,
 \beqnn
 \big\| \boldsymbol{J}_2^{(n)} \big\|_{L^1_T} 
 \ar\leq\ar \frac{1}{\sqrt{n\theta_n}}\big\| (\boldsymbol{V}^{(n)}+\boldsymbol{h}) * \boldsymbol{R}^{(n)} \big \|_{L^1_T} 
 \leq   \frac{1}{\sqrt{n\cdot\theta_n}}\cdot \| \boldsymbol{V}^{(n)}+\boldsymbol{h}\|_{L^\infty_T}^3
 \cdot 
 \mathcal{I}_{\boldsymbol{R}^{(n)}}(T) 
 \to 0. 
 \eeqnn  
 Plugging them into (\ref{eqn.3006}), we have $ \big\|\boldsymbol{V}^{(n)}_\lambda  -\boldsymbol{\mathcal{V}}_\lambda\big\|_{L^1_T}\to0 $ and hence $\big\|\boldsymbol{V}^{(n)}   -\boldsymbol{\mathcal{V}} \big\|_{L^\infty_T} \to 0$ as $n\to\infty$. 
 Additionally, using the triangle inequality and then the inequality $|e^x-1-x-x^2/2| \leq |x|^3$ for any $x\in\mathbb{C}_-$ again we have 
 \beqnn
 \big\| \boldsymbol{W}^{(n)} -\boldsymbol{\mathcal{W}} \big\|_{L^\infty_T} 
 \ar\leq\ar \Big\| n\theta_n\cdot \mathcal{W}\Big(\frac{\boldsymbol{V}^{(n)}+\boldsymbol{h}}{\sqrt{n\theta_n}}\Big)- \frac{1}{2}   \big(  \boldsymbol{V}^{(n)}+\boldsymbol{h} \big)^2 \Big\|_{L^\infty_T} 
 + \frac{1}{2} \cdot \big\|  \big(  \boldsymbol{V}^{(n)}+\boldsymbol{h} \big)^2- \big(  \boldsymbol{\mathcal{V}}+\boldsymbol{h} \big)^2 \big\|_{L^\infty_T}\cr
 \ar\leq\ar  \frac{1}{\sqrt{n\cdot\theta_n}}\cdot \| \boldsymbol{V}^{(n)}+\boldsymbol{h}\|_{L^\infty_T}^3
 + \frac{1}{2} \cdot \big\|  \boldsymbol{V}^{(n)} + \boldsymbol{\mathcal{V}}+2\boldsymbol{h}   \big\|_{L^\infty_T} \cdot \big\|    \boldsymbol{V}^{(n)}- \boldsymbol{\mathcal{V}}  \big\|_{L^\infty_T},
 \eeqnn
 which goes to $0$ as $n\to\infty$. Here we have proved (\ref{eqn.509}).
 For (\ref{eqn.5010}), by the triangle inequality along with (\ref{eqn.506}) and (\ref{eqn.509}), we have for any $T \in [0,T_\infty)$, 
 \beqnn
 \big\| \boldsymbol{\mathcal{V}}  \big\|_{L^\infty_T} \leq \limsup_{n\to\infty}  \big\|\boldsymbol{V}^{(n)}  \big\|_{L^\infty_T} + \limsup_{n\to\infty}  \big\|\boldsymbol{V}^{(n)} - \boldsymbol{\mathcal{V}}  \big\|_{L^\infty_T} \leq C \cdot\big( \|\boldsymbol{f}\|_{L^\infty_T} + \|\boldsymbol{h}\|_{L^\infty_T}^2 \big)  , 
 \eeqnn
 for some constant $C>0$ independent of $\boldsymbol{f}$ and $\boldsymbol{h}$.  
 The desired upper bound for $\big\| \boldsymbol{\mathcal{W}}  \big\|_{L^\infty_T}$ can be proved in the same way. 
 \qed

 \begin{corollary}\label{Corollary.407}
 	Assume $\boldsymbol{ \varPi}(dt)$ has density $\boldsymbol{ \pi}$, for $\boldsymbol{ f} \in L^\infty_{\rm loc}(\mathbb{R}_+;\mathbb{C}^d_-)$ and $\boldsymbol{h} \in L^\infty_{\rm loc}(\mathbb{R}_+;\mathtt{i}\mathbb{R}^d)$, there exists a unique global solution  $ \boldsymbol{\mathcal{V}}\in C(\mathbb{R}_+;\mathbb{C}^d_-) $ to (\ref{Eqn.AffineVolterra}).
 \end{corollary}
 \proof  By Proposition~\ref{Prop.A2}, the Volterra equation (\ref{Eqn.AffineVolterra}) has a unique non-continuous solution $(T_\infty,\boldsymbol{\mathcal{V}}) \in (0,\infty]\times C([0,T_\infty);\mathbb{C}^d)$. 
 By Lemma~\ref{Lemma.ConV}(1), we have $T_\infty =\infty$ and $\boldsymbol{\mathcal{V}}$ is the unique global solution.
 \qed 
 
 \begin{lemma}\label{Lemma.FLFXiM}
 For $\boldsymbol{ f} \in L^\infty_{\rm loc}(\mathbb{R}_+;\mathbb{C}^d_-)$ and $\boldsymbol{h} \in L^\infty_{\rm loc}(\mathbb{R}_+;\mathtt{i}\mathbb{R}^d)$, assume that $(T_\infty,\boldsymbol{\mathcal{V}})\in (0,\infty]\times L^\infty_{\rm loc}([0,T_\infty);\mathbb{C}^d)$ is a non-continuable solution~of (\ref{Eqn.AffineVolterra}) and $(\boldsymbol{ \varXi}, \boldsymbol{M}) \in C(\mathbb{R}_+;\mathbb{R}_+^d\times \mathbb{R}^d)$ is a solution of (\ref{Eqn.HawkesVoletrra.02}), then the representation (\ref{eqn.FL}) holds for any $T\in[0,T_\infty)$.  
 \end{lemma}
 \proof  We define the following  auxiliary process:  for $t\in [0,T]$, 
 \beqlb\label{eqn.4023}
 \begin{split}
 Z_T(t)&:=  \boldsymbol{\mathcal{W}}* d\boldsymbol{\varUpsilon}(T) + \int_0^t \big(  \boldsymbol{f}- \boldsymbol{\mathcal{W}} \big)(T-s) d \boldsymbol{\varXi} (s) + U_T(t),\cr
 U_T(t)&:= \int_0^t  \big(\boldsymbol{\mathcal{V}}+ \boldsymbol{h}\big)(T-s) d \boldsymbol{M}(s).
 \end{split}
 \eeqlb
 It is obvious that $U_T$ is a martingale on $[0,T]$. 
 Similarly as in the proof of Proposition~\ref{Lemma.FourLapGHP},  we can prove that the Dol\'eans-Dade exponential of $U_T$ defined by
 \beqnn
 \mathcal{E}_{U_T}(t)\ar:=\ar \exp\bigg\{ U_T(t)-  \frac{1}{2}\int_0^t \big(\boldsymbol{\mathcal{V}}+ \boldsymbol{h}\big) ^2(T-s) d \boldsymbol{\varXi} (s) \bigg\} \cr
 \ar=\ar  \exp\bigg\{ U_T(t)+\int_0^t \big(  \boldsymbol{f}- \boldsymbol{\mathcal{W}} \big)(T-s)  d \boldsymbol{\varXi} (s) \bigg\},\quad t\in[0,T] , 
 \eeqnn
 is a true martingale. Notice that $\exp\{  Z_T(t) \}= \exp\{\boldsymbol{\mathcal{W}}* d\boldsymbol{\varUpsilon}(T) \} \cdot \mathcal{E}_{U_T}(t)$, we have 
 \beqlb
 \mathbf{E}\Big[\exp\big\{  Z_T(T) \big\}\Big] = \exp\big\{\boldsymbol{\mathcal{W}}* d\boldsymbol{\varUpsilon}(T) \big\},\quad t\in[0,T].
 \eeqlb
 Plugging (\ref{Eqn.AffineVolterra}) into the first equation in  (\ref{eqn.4023}) and then letting $t=T$, we have $$ Z_T(T) =  \boldsymbol{f} * d\boldsymbol{\varXi} (T) + \boldsymbol{h} * d\boldsymbol{M} (T) $$ and hence the representation (\ref{eqn.FL}) holds.
 \qed

  \section{Proofs of main results}
 \label{ProofMainThm}
 \setcounter{equation}{0}
  
 In this section, we provide the detailed proofs for our main results, including Theorem~\ref{MainThm.S-WeakConvergence}, \ref{MainThm.FourierLaplaceFunctional},  \ref{MainThm.SVERepresentation},  \ref{MainThm.J1WeakConvergence}, \ref{Thm.ConvergenceDS01}	and Corollary~\ref{Corollary.ContinuityXiM}, \ref{Corollary.2001}. 
 As a preparation, we start by establishing some uniform moment estimates for    the sequence $\big\{\big( \boldsymbol{\varLambda}^{(n)},\boldsymbol{N}^{(n)},\widetilde{\boldsymbol{N}}^{(n)}\big)\big\}_{n\geq 1}$, which will be used in the proof of their tightness and the characterization of their limit processes. 
 Two direct consequences of our assumptions in Theorem~\ref{MainThm.S-WeakConvergence} are
  \beqlb\label{eqn.500}
 \sup_{n\geq 1} \mathcal{I}_{ \boldsymbol{H}^{(n)}}(T)<\infty
 \quad \mbox{and}\quad 
 \lim_{n\to\infty} \int_0^T f(s) d\mathcal{I}_{ \boldsymbol{H}^{(n)}}(s)=  \lim_{n\to\infty} \int_0^T f(s) \boldsymbol{H}^{(n)}(s)ds =  \int_0^T f(s) d\boldsymbol{\varUpsilon}(s),
 \eeqlb
 for any $T\geq 0$ and any $f\in C(\mathbb{R}_+;\mathbb{C}^d)$. 
 Integrating both sides of (\ref{IntegralLambda}) and then using the stochastic Fubini theorem; see Theorem D.2 in \cite{Xu2024b}, the integrated process $\mathcal{I}_{ \boldsymbol{\varLambda}^{(n)}}$ satisfies the following equation 
 \beqlb\label{eqn.IntLambda}
 \mathcal{I}_{ \boldsymbol{\varLambda}^{(n)}} (t)
 \ar=\ar  \mathcal{I}_{\boldsymbol{H}^{(n)}}(t)
 + \mathcal{I}_{ \boldsymbol{R}^{(n)}}* d\widetilde{\boldsymbol{N}}^{(n)} (t) = \mathcal{I}_{\boldsymbol{H}^{(n)}}(t)
  +  d\mathcal{I}_{ \boldsymbol{R}^{(n)}}*\widetilde{\boldsymbol{N}}^{(n)} (t)  ,\quad t\geq 0.
 \eeqlb

 \begin{proposition}\label{Prop.601}
 	For each $T\geq 0$, we have 
 	\beqlb\label{eqn.601}
 	\sup_{n\geq 1}\mathbf{E}\big[\mathcal{I}_{ \boldsymbol{\varLambda}^{(n)}}(T)\big] + 	\sup_{n\geq 1}\mathbf{E}\big[\boldsymbol{N}^{(n)}(T)\big]  + \sup_{n\geq 1} \mathbf{E}\bigg[\sup_{t\in[0,T]} \Big|\widetilde{\boldsymbol{N}}^{(n)}(t)\Big|^2\bigg]<\infty.
 	\eeqlb 
 \end{proposition}
 \proof Taking expectation on both sides of (\ref{eqn.IntLambda}) with $t=T$, we have 
 \beqlb\label{eqn.5001}
 \mathbf{E}\big[\boldsymbol{N}^{(n)}(T)\big]=\mathbf{E}\big[\mathcal{I}_{ \boldsymbol{\varLambda}^{(n)}}(T)\big] 
 = \mathcal{I}_{\boldsymbol{H}^{(n)}} (T) ,
 \eeqlb
 which is bounded uniformly in $n\geq 1$ because of (\ref{eqn.500}). 
 Moreover, by using the Burkholder-Davis-Gundy inequality to the martingale $\widetilde{\boldsymbol{N}}^{(n)}$, 
 \beqnn
 \mathbf{E}\bigg[\sup_{t\in[0,T]} \big|\widetilde{\boldsymbol{N}}^{(n)}(t)\big|^2\bigg] \leq C\cdot \mathbf{E}\big[ \boldsymbol{N}^{(n)}(T) \big] = C\cdot \mathbf{E}\big[\mathcal{I}_{ \boldsymbol{\varLambda}^{(n)}}(T)\big]
 = C\cdot \mathcal{I}_{\boldsymbol{H}^{(n)}} (T), 
 \eeqnn 
 for some constant $C>0$ independent of $n$ and hence the desired upper bound holds.
 \qed 
 
 \begin{corollary}\label{Coro.601}
   For each $T\geq 0$, we have $\big\| \boldsymbol{N}^{(n)} -\mathcal{I}_{ \boldsymbol{\varLambda}^{(n)}} \big\|_{L^\infty_T} \overset{\rm p}\to 0$  as $n\to\infty$.
 \end{corollary}
 \proof By the definitions of $ \boldsymbol{N}^{(n)}$ and $ \widetilde{\boldsymbol{N}}^{(n)}$; see (\ref{eqn.ScaledProcess}), we have $\boldsymbol{N}^{(n)}(t)-\mathcal{I}_{ \boldsymbol{\varLambda}^{(n)}}(t) = \widetilde{\boldsymbol{N}}^{(n)}(t)/\sqrt{n\theta_n}$ for $t\geq 0$. 
 From this and (\ref{eqn.601}), 
 \beqnn
  \mathbf{E}\bigg[ \sup_{t\in[0,T]} \Big| \boldsymbol{N}^{(n)}(t)-\mathcal{I}_{ \boldsymbol{\varLambda}^{(n)}}(t) \Big|^2 \bigg] 
  =  \frac{1}{n\theta_n} \cdot \mathbf{E}\bigg[ \sup_{t\in[0,T]} \Big|\widetilde{\boldsymbol{N}}^{(n)}(t) \Big|^2 \bigg] \to 0,
 \eeqnn
 as $n\to\infty$, which immediately induces the  desired limit. 
 \qed

 \begin{proposition}\label{Prop.602}
  The martingale sequence $\{ \widetilde{\boldsymbol{N}}^{(n)} \}_{n\geq 1}$ is uniformly tight; see Definition 6.1 in \cite[p.377]{JacodShiryaev2003}.
 \end{proposition} 
 \proof By (\ref{eqn.QV}), the martingale $\widetilde{\boldsymbol{N}}^{(n)}$ has predictable quadratic co-variation 
 \beqlb\label{eqn.615}
 \big\langle\widetilde{N}^{(n)}_i,\widetilde{N}^{(n)}_j\big\rangle= \mathbf{1}_{\{i=j \}}\cdot 	\mathcal{I}_{\varLambda^{(n)}_i},\quad i,j\in\mathtt{D}.
 \eeqlb
 By (\ref{eqn.601}), the sequence $\big\{ \big\langle\widetilde{N}^{(n)}_i,\widetilde{N}^{(n)}_j \big\rangle_t \big\}_{n\geq 1}$ is tight for each $t\geq 0$. 
 Since all jumps of  $\widetilde{\boldsymbol{N}}^{(n)}$ are uniformly bounded by $1/\sqrt{n\theta_n}$, the uniformly tightness of $\{ \widetilde{\boldsymbol{N}}^{(n)} \}_{n\geq 1}$ follows from Proposition~6.13 in \cite[p.379]{JacodShiryaev2003}. 
 \qed

 \medskip
 
  \begin{lemma} \label{Lemma.Tightness}
  The sequence $\{ (\mathcal{I}_{\boldsymbol{\varLambda}^{(n)}}, \boldsymbol{N}^{(n)}, \widetilde{\boldsymbol{N}}^{(n)}) \}_{n\geq 1}$ is relatively compact in   $D_S(\mathbb{R}_+;\mathbb{R}_+^d\times\mathbb{R}_+^d\times\mathbb{R}^d)$.
 \end{lemma}
 \proof The relative compactness of the two sequences $\{ \mathcal{I}_{ \boldsymbol{\varLambda}^{(n)}}  \}_{n\geq 1}$ and $\{\boldsymbol{N}^{(n)} \}_{n\geq 1}$ in $D_S(\mathbb{R}_+;\mathbb{R}^d_+)$  follows directly from Proposition~\ref{Appendix.Prop.S02} along with their monotonicity and (\ref{eqn.601}).
 By Proposition~\ref{Appendix.Prop.S04} and \ref{Prop.602},  the martingale sequence  $\{ \widetilde{\boldsymbol{N}}^{(n)} \}_{n\geq 1}$ is also relatively compact in $D_S(\mathbb{R}_+;\mathbb{R}^d)$. 
 Finally,  the relative compactness of $\{ (\mathcal{I}_{ \boldsymbol{\varLambda}^{(n)}},\boldsymbol{N}^{(n)}, \widetilde{\boldsymbol{N}}^{(n)}) \}_{n\geq 1}$  in $D_S(\mathbb{R}_+;\mathbb{R}^d_+\times\mathbb{R}^d_+\times\mathbb{R}^d)$ follows from Proposition~\ref{Appendix.Prop.S021}.
 \qed 
 
 \begin{lemma}\label{Lemma.CluterPoint}
 Assume that $(\boldsymbol{\varXi},\boldsymbol{\eta},\boldsymbol{M})$ is an accumulation point of $\{ (\mathcal{I}_{\boldsymbol{\varLambda}^{(n)}}, \boldsymbol{N}^{(n)}, \widetilde{\boldsymbol{N}}^{(n)}) \}_{n\geq 1}$ in $D_S(\mathbb{R}_+;\mathbb{R}_+^d\times\mathbb{R}_+^d\times\mathbb{R}^d)$ along a subsequence $\{n_k\}_{k\geq 1}$, we have the following.
 \begin{enumerate}
  \item[(1)] The two processes $\boldsymbol{\varXi}$ and $\boldsymbol{\eta}$ are identical almost surely, i.e. $\boldsymbol{\varXi}\overset{\rm a.s.}=\boldsymbol{\eta}$.
 		
  \item[(2)] There exists a countable subset $I_1\subset \mathbb{R}_+$ such that for any $\boldsymbol{f} \in C^1(\mathbb{R}_+;\mathbb{C}_-^d)$  and  $\boldsymbol{h}\in C^1(\mathbb{R}_+;\mathtt{i}\mathbb{R}^d)$,
  \beqnn
  \big( \boldsymbol{f} * d\mathcal{I}_{\boldsymbol{\varLambda}^{(n_k)}}, \boldsymbol{h}* d\widetilde{\boldsymbol{N}}^{(n_k)}   \big)
  \overset{\rm f.f.d.}\longrightarrow 
  \big( \boldsymbol{f} * d\boldsymbol{\varXi}  + \boldsymbol{h}* d\boldsymbol{M}  \big),
 \eeqnn
 as $k\to\infty$ along $ \mathbb{R}_+\setminus I_1$. 
 Moreover, it holds along $ \mathbb{R}_+$ if $\boldsymbol{f}(0)=\boldsymbol{h}(0) = \mathbf{0}$. 
 
 \item[(3)] If $\boldsymbol{\varUpsilon}$ is continuous on $\mathbb{R}_+$, then $(\boldsymbol{\varXi}, \boldsymbol{M})$ is stochastically continuous. 
 \end{enumerate} 
 \end{lemma}
 \proof Claim (1) follows directly from Corollary~\ref{Coro.601}. For claim (2), by Fubini's theorem we have 
 \beqlb\label{eqn.FubiniTrans}
 \begin{split}
 \boldsymbol{f} * d\mathcal{I}_{\boldsymbol{\varLambda}^{(n_k)}} (T) &= \boldsymbol{f}(0)\cdot \mathcal{I}_{\boldsymbol{\varLambda}^{(n_k)}}(T) +  \boldsymbol{f}'* \mathcal{I}_{\boldsymbol{\varLambda}^{(n_k)}}(T),  \\[1mm]
 \boldsymbol{h}* d\widetilde{\boldsymbol{N}}^{(n_k)} (T)
 &=\boldsymbol{h}(0)\cdot \widetilde{\boldsymbol{N}}^{(n_k)} (T)
 +\boldsymbol{h}' * \widetilde{\boldsymbol{N}}^{(n_k)} (T),
 \\[2mm]
 \boldsymbol{f} * d\boldsymbol{\varXi} (T) &= 
 \boldsymbol{f}(0) \cdot \boldsymbol{\varXi} (T)
 +\boldsymbol{f}' * \boldsymbol{\varXi} (T),
 \\[2mm]
 \boldsymbol{h}* d\boldsymbol{M} (T)
 &= \boldsymbol{h}(0)\cdot \boldsymbol{M} (T)
 + \boldsymbol{h}'* \boldsymbol{M} (T),
 \end{split}
 \eeqlb 
 for any $T\geq 0$.
  By proposition~\ref{Appendix.Prop.S020}(1), there exists a countable set $I_1\in\mathbb{R}_+$, which is empty if $\boldsymbol{f}(0)=\boldsymbol{h}(0) = \mathbf{0}$, such that 
  \beqnn
  (\mathcal{I}_{\boldsymbol{\varLambda}^{(n_k)}}, \widetilde{\boldsymbol{N}}^{(n_k)})\overset{\rm f.f.d.}\longrightarrow  (\boldsymbol{\varXi} ,\boldsymbol{M} ),
  \eeqnn
  as $k\to\infty$ along $\mathbb{R}_+\setminus I_1$. 
 Moreover, by Proposition~\ref{Appendix.Prop.S020}(4) and the continuity of  $(\boldsymbol{f}',\boldsymbol{h}')$ we also have
  \beqnn
  \Big( \boldsymbol{f}'* \mathcal{I}_{\boldsymbol{\varLambda}^{(n_k)}}, \boldsymbol{h}' * \widetilde{\boldsymbol{N}}^{(n_k)}\Big)
  \overset{\rm f.f.d.}\longrightarrow  
  \big(\boldsymbol{f}' * \boldsymbol{\varXi}, \boldsymbol{h}'* \boldsymbol{M} \big)
  \eeqnn
  along $\mathbb{R}_+$.  Consequently, claim (2) follows.
 
 We now start to prove claim (3). For any $s,t\geq 0$, choose two sequences $\{ t_i\}_{i\geq 1}\subset I_1$ and $\{ s_i\}_{i\geq 1}\subset I_1$ such that $s_i\to s+$ and $t_i\to t+$ as $i\to\infty$. 
 By Proposition~\ref{Appendix.Prop.S020}(1), we have  
 \beqnn
 \Big(\mathcal{I}_{\boldsymbol{\varLambda}^{(n_k)}}(s_i), \widetilde{\boldsymbol{N}}^{(n_k)}(s_i), \mathcal{I}_{\boldsymbol{\varLambda}^{(n_k)}}(t_i),\widetilde{\boldsymbol{N}}^{(n_k)}(t_i)\Big) 
 \overset{\rm d}\to 
 \big( \boldsymbol{\varXi}(s_i),\boldsymbol{M}(s_i), \boldsymbol{\varXi}(t_i), \boldsymbol{M}(t_i) \big).
 \eeqnn
 From this, the last equality in (\ref{eqn.5001}) and the fact that $\mathcal{I}_{\boldsymbol{H}^{(n_k)}} \to\boldsymbol{\varUpsilon}$ uniformly on compacts as $k\to\infty$, we have for any $K>0$,
 \beqnn
 \mathbf{E}\Big[ \big|\boldsymbol{\varXi}(t_i)-\boldsymbol{\varXi}(s_i) \big|\wedge K \Big] 
 \ar=\ar \lim_{k\to\infty} 	\mathbf{E}\Big[ \big| \mathcal{I}_{\boldsymbol{\varLambda}^{(n_k)}}(t_i)- \mathcal{I}_{\boldsymbol{\varLambda}^{(n_k)}}(s_i)\big|\wedge K \Big]\cr
 \ar\leq\ar \lim_{k\to\infty} \Big|\mathbf{E}\big[   \mathcal{I}_{\boldsymbol{\varLambda}^{(n_k)}}(t_i)  \big]- \mathbf{E}\big[  \mathcal{I}_{\boldsymbol{\varLambda}^{(n_k)}}(s_i)  \big]\Big|  =\big| \boldsymbol{\varUpsilon}(t_i)- \boldsymbol{\varUpsilon}(s_i)\big|.
 \eeqnn
 By the monotone convergence theorem along with the right-continuity of $(\boldsymbol{\varXi},\boldsymbol{\varUpsilon})$, 
 \beqnn
 \mathbf{E}\Big[ \big|\boldsymbol{\varXi}(t)-\boldsymbol{\varXi}(s)\big| \wedge K \Big] = \lim_{i\to\infty} 	\mathbf{E}\Big[ \big|\boldsymbol{\varXi}(t_i)-\boldsymbol{\varXi}(s_i)\big|\wedge K \Big] 
 \leq \big| \boldsymbol{\varUpsilon}(t)- \boldsymbol{\varUpsilon}(s)\big|,
 \eeqnn
 which goes to $0$ as $|t-s|\to 0$.  This, together with the arbitrariness of $K$, yields the stochastic continuity of $\boldsymbol{\varXi}$.
 Similarly, we also have 
 \beqnn
 \mathbf{E}\Big[ \big|\boldsymbol{M}(t_i)-\boldsymbol{M}(s_i)\big|^2\wedge K \Big] 
 \ar=\ar
 \lim_{k\to\infty}\mathbf{E}\Big[ \big|\widetilde{\boldsymbol{N}}^{(n_k)}(t_i)-\widetilde{\boldsymbol{N}}^{(n_k)}(s_i)\big|^2\wedge K \Big] \cr
 \ar\leq\ar 	\lim_{k\to\infty}\mathbf{E}\Big[ \big|\widetilde{\boldsymbol{N}}^{(n_k)}(t_i)-\widetilde{\boldsymbol{N}}^{(n_k)}(s_i)\big|^2  \Big]. 
 \eeqnn
 Applying the Burkholder-Davis-Gundy inequality to the last expectation, we have 
 \beqnn
 \mathbf{E}\Big[ \big|\boldsymbol{M}(t_i)-\boldsymbol{M}(s_i)\big|^2\wedge K \Big] \leq C\cdot \lim_{k\to\infty} \int_{s_i}^{t_i} \mathbf{E}[\boldsymbol{\varLambda}^{(n_k)}(r)]dr,
 \eeqnn
 for some constant $C>0$. 
 The previous argument tells that 
 \beqnn
 \mathbf{E}\Big[ \big|\boldsymbol{M}(t)-\boldsymbol{M}(s)\big|^2\wedge K \Big]\leq C\cdot \big| \boldsymbol{\varUpsilon}(t)- \boldsymbol{\varUpsilon}(s)\big|,
 \eeqnn
 which also goes $0$ as $|t-s|\to0$ and hence $\boldsymbol{M}$ is also stochastically continuous. 
 \qed
 
%
%
 \begin{lemma}\label{Lemma.ExistUniqueV}
 	Theorem~\ref{MainThm.FourierLaplaceFunctional}(1) holds.
 	\end{lemma}
 \proof For each $j \in \mathtt{D}$ and $n\geq 1$, we set the parameter $( \boldsymbol{\mu}_n ,  \boldsymbol{\phi}_n )$ as in Remark~\ref{Remark.01}(1) with $\boldsymbol{a}=\mathbf{0}$ and $\boldsymbol{A}=\mathbf{Id}$. 
 In this case, we have  $ \boldsymbol{H}^{(n)} = \boldsymbol{R}^{(n)}_{\cdot j} = \big( R^{(n)}_{ij} \big)_{i\in\mathtt{D}}$. 
 Plugging this back into the right-hand side of the second equality in (\ref{eqn.503}) and then using the first equality in (\ref{eqn.504}), we have 
 \beqlb\label{eqn.5004}
  \mathbf{E} \Big[ \exp\Big\{\boldsymbol{f} *\boldsymbol{\varLambda}^{(n)}(T) + \boldsymbol{h}* d \widetilde{\boldsymbol{N}}^{(n)}(T)  \Big\} \Big]
 \ar=\ar \exp \big\{   \boldsymbol{W}^{(n)} * \boldsymbol{R}^{(n)}_{\cdot j}(T) \big\}
 = \exp\big\{V^{(n)}_j(T)\big\}.  
 \eeqlb 
 By Lemma~\ref{Lemma.CluterPoint}(2), we have for any
 $\boldsymbol{f} \in C^1(\mathbb{R}_+;\mathbb{C}_-^d)$  and  $\boldsymbol{h}\in C^1(\mathbb{R}_+;\mathtt{i}\mathbb{R}^d)$, 
 \beqnn
 \mathbf{E} \Big[ \exp\Big\{\boldsymbol{f} *\boldsymbol{\varLambda}^{(n_k)}(T) + \boldsymbol{h}* d \widetilde{\boldsymbol{N}}^{(n_k)}(T)  \Big\} \Big] 
 \to  \mathbf{E} \Big[ \exp\big\{\boldsymbol{f} *d\boldsymbol{\varXi}(T) + \boldsymbol{h}* d \boldsymbol{M}(T)  \big\} \Big],\quad 
 \eeqnn
 as $k\to\infty$ for any $T\in  \mathbb{R}_+ \setminus I_1$, which together with (\ref{eqn.5004}) induces that the limit $V^{(n_k)}_j(T) \to  \mathcal{V}_j(T) \in \mathbb{C}_-$ exists and 
 \beqlb\label{eqn.FourierLap.01}
 \mathbf{E} \Big[ \exp\big\{\boldsymbol{f} *d\boldsymbol{\varXi}(T) + \boldsymbol{h}* d \boldsymbol{M}(T)  \big\} \Big] = \exp\big\{ \mathcal{V}_j(T) \big\}. 
 \eeqlb
 By Lemma~\ref{Lemma.ExitenceRV}(1), the limit function $\boldsymbol{\mathcal{V}}:= \big(\mathcal{V}_j\big)_{j\in\mathtt{D}}  \in L^\infty_{\rm loc}(\mathbb{R}_+;\mathbb{C}_-^d)$ is the unique solution to (\ref{Eqn.AffineVolterra}). 
 By applying the dominated convergence theorem along with the fact that $(\boldsymbol{\varXi},\boldsymbol{M}) \in D(\mathbb{R}_+;\mathbb{R}^d_+\times \mathbb{R}^d)$ and the last two equalities in (\ref{eqn.FubiniTrans}) to the expectation on the left-hand side of (\ref{eqn.FourierLap.01}), one can see that it is c\`adl\`ag and hence  $\boldsymbol{\mathcal{V}} \in D(\mathbb{R}_+;\mathbb{C}^d_-) $. Here we have proved the first claim.

 We now prove the second claim by distinguishing the following two cases.
 \begin{enumerate}
 	\item[$\bullet$] \textit{$\boldsymbol{f}(0)=\boldsymbol{h}(0)=\mathbf{0}$.} We extend $\boldsymbol{f}'$ and $\boldsymbol{h}'$ to the whole real line with $\boldsymbol{f}'(t)=\boldsymbol{f}'(0)$ and $\boldsymbol{h}'(t)=\boldsymbol{h}'(0)$ for $t\leq 0$. 
 	The last two equalities in (\ref{eqn.FubiniTrans}) indicate that
 	\beqnn
 	\boldsymbol{f} *d\boldsymbol{\varXi}(T) \in C(\mathbb{R}_+;\mathbb{C}_-^d)
 	\quad\mbox{and}\quad 
 	\boldsymbol{h}* d \boldsymbol{M}  \in C(\mathbb{R}_+;\mathtt{i}\mathbb{R}^d) . 
 	\eeqnn  
 	The continuity of $\boldsymbol{\mathcal{V}}$ following directly by applying the dominated convergence theorem again to (\ref{eqn.FourierLap.01}).
 	
 	\item[$\bullet$] \textit{$\underline{\boldsymbol{\varPi}}$ is continuous on $\mathbb{R}_+$.} 
 	Applying Lemma~\ref{Lemma.CluterPoint}(3) and the continuity of $(\boldsymbol{f}' * \boldsymbol{\varXi},\boldsymbol{h}' * \boldsymbol{M})$ to the last two equalities in (\ref{eqn.FubiniTrans}),  we have that both $\boldsymbol{f} * d\boldsymbol{\varXi}$ and $\boldsymbol{h}* d\boldsymbol{M}$ are stochastically continuous, which induces that the expectation on the left side of (\ref{eqn.FourierLap.01}) is continuous and so is $ \mathcal{V}_j$. Thus $\boldsymbol{ \mathcal{V}}$ is continuous. 
   \qed 
 	\end{enumerate}

  \begin{lemma}\label{Lemma.FourierLaplaceCluster}
 	Theorem~\ref{MainThm.FourierLaplaceFunctional}(2) holds for any accumulation point of $\big\{ (\mathcal{I}_{\boldsymbol{\varLambda}^{(n)}}, \boldsymbol{N}^{(n)}, \widetilde{\boldsymbol{N}}^{(n)}) \big\}_{n\geq 1}$. 
 \end{lemma}
 \proof 
 Using the settings and results in Lemma~\ref{Lemma.CluterPoint} and \ref{Lemma.FLFLambdaNn},
 we have
 \beqnn
  \exp \big\{   \boldsymbol{W}^{(n_k)} *  \boldsymbol{H}^{(n_k)} (t) \big\} \ar=\ar \mathbf{E}\big[ \exp \big\{ \boldsymbol{f}* d\mathcal{I}_{\boldsymbol{\varLambda}^{(n_k)}} (t) +   \boldsymbol{h}* d\widetilde{\boldsymbol{N}}^{(n_k)} (t) \big\} \big] \\[2mm] 
  \ar\to\ar  
  \mathbf{E}\big[ \exp \big\{  \boldsymbol{f} * d\boldsymbol{\varXi}(t) +   \boldsymbol{h}* d\boldsymbol{M} (t) \big\} \big],
  \eeqnn
  as $k\to\infty$ for any $t\in \mathbb{R}_+\setminus I_1$. 
  We now prove that  there exists a subsequence of $\{  \boldsymbol{W}^{(n_k)} *  \boldsymbol{H}^{(n_k)} \}_{k\geq 1}$ that converges  to $\boldsymbol{\mathcal{W}}  * d \boldsymbol{\varUpsilon}$ almost everywhere.
  Actually,  by the triangle inequality,
  \beqnn
  \big|\boldsymbol{W}^{(n_k)} * \boldsymbol{H}^{(n_k)}  (t)- \boldsymbol{\mathcal{W}}  * d \boldsymbol{\varUpsilon}(t)\big|
  \leq \big|\boldsymbol{W}^{(n_k)}- \boldsymbol{\mathcal{W}}\big| * \boldsymbol{H}^{(n_k)} (t)   +\big| \boldsymbol{\mathcal{W}}* \boldsymbol{H}^{(n_k)} (t)- \boldsymbol{\mathcal{W}}  * d \boldsymbol{\varUpsilon}(t)\big|,
  \eeqnn 
  for any $t\geq 0$.  
  By using  (\ref{Eqn.BoundConvIRn}) and then Lemma~\ref{Lemma.ExitenceRV}(2), he first term on the right-hand side of this inequality can be bounded by
  \beqnn
  \big|\boldsymbol{W}^{(n_k)}- \boldsymbol{\mathcal{W}}\big| *  \boldsymbol{H}^{(n_k)} (t) \leq 
  \big\|\boldsymbol{W}^{(n_k)}- \boldsymbol{\mathcal{W}}\big\|_{L^\infty_t} \cdot \sup_{n\geq 1} \mathcal{I}_{\boldsymbol{H}^{(n_k)}} (t)\leq C\cdot  \big\|\boldsymbol{W}^{(n_k)}- \boldsymbol{\mathcal{W}}\big\|_{L^\infty_t},
  \eeqnn
 which goes to $0$ as $n\to\infty$.
 By using Remark~\ref{Remark.3.5} as well as the fact that $\boldsymbol{\mathcal{W}} \in D(\mathbb{R}_+;\mathbb{C}^d)$ and the assumption that $\mathcal{I}_{\boldsymbol{H}^{(n_k)}}\to \boldsymbol{\varUpsilon} $ for all continuous points of $ \boldsymbol{\varUpsilon}$, we have 
  \beqnn
  \big\| \boldsymbol{\mathcal{W}}* \boldsymbol{H}^{(n_k)}- \boldsymbol{\mathcal{W}}  * d \boldsymbol{\varUpsilon}\big\|_{L^1_T} \to 0,
  \eeqnn
  as $k\to\infty$ for any $T\geq 0$. Hence one can always find a subsequence $\{n_l \}_{l\geq 1}\subset\{ n_k \}_{k\geq 1}$ such that 
  $$  \boldsymbol{\mathcal{W}}* \boldsymbol{H}^{(n_l)} \overset{\rm a.e.}\to \boldsymbol{\mathcal{W}}  * d \boldsymbol{\varUpsilon}$$ as $l\to \infty$. 
  Putting these results together, we have 
  \beqlb\label{eqn.50011}
   \boldsymbol{W}^{(n_l)} * \boldsymbol{H}^{(n_l)}  (t) \to \boldsymbol{\mathcal{W}}  * d \boldsymbol{\varUpsilon}(t) 
   \quad \mbox{and hence}\quad 
   \mathbf{E}\big[ \exp \big\{  \boldsymbol{f} * d\boldsymbol{\varXi}(t) +   \boldsymbol{h}* d\boldsymbol{M} (t) \big\} \big] = \exp\{ \boldsymbol{\mathcal{W}}  * d \boldsymbol{\varUpsilon}(t) \}. 
  \eeqlb
  for almost every $t\geq 0$.
  Moreover, by the fact that $(\boldsymbol{\varXi},\boldsymbol{M}) \in D(\mathbb{R}_+;\mathbb{R}^d_+\times \mathbb{R}^d)$ and the last two equality in (\ref{eqn.FubiniTrans}), both $ \boldsymbol{f} * d\boldsymbol{\varXi}$ and $  \boldsymbol{h}* d\boldsymbol{M}$ are c\`adl\`ag and hence   $\mathbf{E}\big[ \exp \big\{  \boldsymbol{f} * d\boldsymbol{\varXi} +   \boldsymbol{h}* d\boldsymbol{M} \big\} \big] \in D(\mathbb{R}_+;\mathbb{C}_-)$. 
  Moreover, notice that $\boldsymbol{\mathcal{W}}  * d \boldsymbol{\varUpsilon} $ is right-continuous because both $\boldsymbol{\mathcal{W}} $  and $\boldsymbol{\varUpsilon}$  are c\`adl\`ag, then (\ref{eqn.50011}) holds for all $t\geq 0$. 
  \qed

  Armed with the previous auxiliary results, we are now able to provide the detailed proofs for our main theorems one by one. 
  
  \medskip
  
  \textbf{\textit{Proof of Theorem~\ref{MainThm.S-WeakConvergence}.}}  
  By Lemma~\ref{Lemma.Tightness} and \ref{Lemma.CluterPoint}, it remains to prove that all accumulation points are identical in distribution.  
  Assume that $(\boldsymbol{\varXi},\boldsymbol{M})$ and $(\boldsymbol{\varXi}^*,\boldsymbol{M}^*)$ are two limits of $\{ (\mathcal{I}_{ \boldsymbol{\varLambda}^{(n)}},  \widetilde{\boldsymbol{N}}^{(n)}) \}_{n\geq 1}$ along the subsequences $\{ n_k \}_{k\geq 1}$ and $\{ n_k^* \}_{n\geq 1}$ respectively.
  Since $D(\mathbb{R}_+;\mathbb{R}_+^d \times \mathbb{R}^d)$ is a subspace of $L^\infty_{\rm loc}(\mathbb{R}_+;\mathbb{R}_+^d \times\mathbb{R}^d)$ with dual space $L^1_{\rm loc}(\mathbb{R}_+;\mathbb{R}_+^d \times\mathbb{R}^d)$, hence it suffices to prove that 
  \beqlb\label{eqn.5006}
  \big( \boldsymbol{f} * \boldsymbol{\varXi} (T), \boldsymbol{h}* \boldsymbol{M} (T) \big)
  \overset{\rm d}=
  \big( \boldsymbol{f} * \boldsymbol{\varXi}^* (T), \boldsymbol{h}* \boldsymbol{M}^* (T) \big),\quad T\geq 0,
  \eeqlb
  for any $\boldsymbol{f}\in L^1_{\rm loc}(\mathbb{R}_+;\mathbb{C}_-^d)$ and $\boldsymbol{h}\in L^1_{\rm loc}(\mathbb{R}_+;\mathtt{i}\mathbb{R}^d)$. 
  Additionally, notice that $C(\mathbb{R}_+;\mathbb{C}_-^d)\subset L^1_{\rm loc}(\mathbb{R}_+;\mathbb{C}_-^d)$ and $C(\mathbb{R}_+;\mathtt{i}\mathbb{R}^d)\subset L^1_{\rm loc}(\mathbb{R}_+;\mathtt{i}\mathbb{R}^d)$ are dense, 
  thus we just need to prove (\ref{eqn.5006}) with $\boldsymbol{f}\in C(\mathbb{R}_+;\mathbb{C}_-^d)$ and $\boldsymbol{h}\in C(\mathbb{R}_+;\mathtt{i}\mathbb{R}^d)$. 
  By the last two equalities in (\ref{eqn.FubiniTrans}), we  can write (\ref{eqn.5006}) as 
  \beqlb \label{eqn.5007}
  \big( \mathcal{I}_{\boldsymbol{f}} * d\boldsymbol{\varXi} (T), \mathcal{I}_{\boldsymbol{h}}* d\boldsymbol{M} (T) \big)
  \overset{\rm d}=
  \big( \mathcal{I}_{\boldsymbol{f}} * d\boldsymbol{\varXi}^* (T), \mathcal{I}_{\boldsymbol{h}}* d\boldsymbol{M}^* (T) \big),\quad T\geq 0,
  \eeqlb
  which holds if and only if 
   \beqnn
   \mathbf{E}\Big[ \exp \big\{ \mathcal{I}_{\boldsymbol{f}} * d\boldsymbol{\varXi} (T)  + \mathcal{I}_{\boldsymbol{h}}* d\boldsymbol{M} (T) \big\} \Big]
  = \mathbf{E}\Big[ \exp \big\{ \mathcal{I}_{\boldsymbol{f}} * d\boldsymbol{\varXi}^* (T)+ \mathcal{I}_{\boldsymbol{h}}* d\boldsymbol{M}^* (T) \big\} \Big].
  \eeqnn
  It follows directly from Lemma~\ref{Lemma.FourierLaplaceCluster} and the uniqueness of solutions to (\ref{Eqn.AffineVolterra}); see Lemma~\ref{Lemma.ExistUniqueV}. 
  \qed 
  
  \medskip

 \textbf{\textit{Proof of Theorem~\ref{MainThm.FourierLaplaceFunctional}.}}
 The desired two claims have been proved in Lemma~\ref{Lemma.ExistUniqueV} and~\ref{Lemma.FourierLaplaceCluster}. 
 \qed
 
 \medskip
 
 \textbf{\textit{Proof of Theorem~\ref{MainThm.SVERepresentation}.}} By Proposition~\ref{Appendix.Prop.S020}(1), without loss of generality we may assume that 
 \beqlb\label{eqn.6004}
 \Big(\mathcal{I}_{ \boldsymbol{\varLambda}^{(n)}},  \widetilde{\boldsymbol{N}}^{(n)} \Big) \overset{\rm a.s.}\to \big(\boldsymbol{\varXi},\boldsymbol{M} \big) ,
 \eeqlb
 in $D_S(\mathbb{R}_+;\mathbb{R}_+^d\times\mathbb{R}^d)$. 
 For each $T\geq 0$, integrating both sides of (\ref{eqn.IntLambda}) over $[0,T]$ and then using Fubini's theorem we have 
 \beqlb\label{eqn.5017}
 \int_0^T \mathcal{I}_{ \boldsymbol{\varLambda}^{(n)}} (t) dt
 \ar=\ar  \int_0^T\mathcal{I}_{\boldsymbol{H}^{(n)}}(t)dt
 +  \int_0^T d\mathcal{I}_{ \boldsymbol{R}^{(n)}}*\widetilde{\boldsymbol{N}}^{(n)} (t)dt 
 =  \int_0^T\mathcal{I}_{\boldsymbol{H}^{(n)}}(t)dt
 +  \mathcal{I}_{ \boldsymbol{R}^{(n)}}*\widetilde{\boldsymbol{N}}^{(n)} (T) .
 \eeqlb
 By Proposition~\ref{Appendix.Prop.S01}(2), we have almost surely
 \beqlb\label{eqn.5013}
 \int_0^t \mathcal{I}_{ \boldsymbol{\varLambda}^{(n)}} (s) ds\to  \int_0^t \boldsymbol{\varXi} (s) ds
 \quad \mbox{and}\quad 
 \int_0^t \mathcal{I}_{\boldsymbol{H}^{(n)}}(s)ds \to  \int_0^t \boldsymbol{\varUpsilon} (s)ds,\quad t\geq 0,
 \eeqlb
 as $n\to\infty$. Moreover, by the triangle inequality,
 \beqnn
  \big|\mathcal{I}_{ \boldsymbol{R}^{(n)}}*\widetilde{\boldsymbol{N}}^{(n)}(t) - \underline{ \boldsymbol{\varPi}}* \boldsymbol{M} (t)\big|
  \ar\leq\ar  \big|  \mathcal{I}_{ \boldsymbol{R}^{(n)}}-\underline{ \boldsymbol{\varPi}}\big|*\big|\widetilde{\boldsymbol{N}}^{(n)}\big|(t) +\big|  \underline{ \boldsymbol{\varPi}}*\widetilde{\boldsymbol{N}}^{(n)}(t)-\underline{ \boldsymbol{\varPi}}* \boldsymbol{M} (t)\big|.
 \eeqnn
  By Proposition~\ref{Appendix.Prop.S01}(2) again, we have almost surely
  \beqnn
  \big|  \underline{ \boldsymbol{\varPi}}*\widetilde{\boldsymbol{N}}^{(n)}(t)-\underline{ \boldsymbol{\varPi}}* \boldsymbol{M} (t)\big|\to 0 ,\quad t\geq 0.
  \eeqnn
 Moreover, by (\ref{Eqn.BoundConvIRn}) and (\ref{eqn.601}) we also have almost surely,
 \beqnn
 \big|  \mathcal{I}_{ \boldsymbol{R}^{(n)}}-\underline{ \boldsymbol{\varPi}}\big|*\big|\widetilde{\boldsymbol{N}}^{(n)}\big|(t)\leq 
 \big\|  \mathcal{I}_{ \boldsymbol{R}^{(n)}}-\underline{ \boldsymbol{\varPi}}\big\|_{L^1_t} \cdot \Big|\sup_{s\in[0,t]}\widetilde{\boldsymbol{N}}^{(n)}(t)\Big| \to 0,\quad t\geq 0.
 \eeqnn
 Putting these results together, we have almost surely,
 \beqnn
 \mathcal{I}_{ \boldsymbol{R}^{(n)}}*\widetilde{\boldsymbol{N}}^{(n)}(t) \to  \underline{ \boldsymbol{\varPi}}* \boldsymbol{M} (t),\quad t\geq 0.
 \eeqnn
 Taking this and (\ref{eqn.5013}) back into (\ref{eqn.5017}) and then passing both sides  to the corresponding limits, we have almost surely
 \beqlb\label{eqn.5018}
  \int_0^T  \boldsymbol{\varXi} (t) dt = \int_0^T  \boldsymbol{\varUpsilon} (t)dt
 + \underline{ \boldsymbol{\varPi}}* \boldsymbol{M} (T),\quad T\geq 0.
 \eeqlb
 By Fubini's theorem, we have 
 \beqnn
  \underline{ \boldsymbol{\varPi}}* \boldsymbol{M} (T)=  \int_0^T \boldsymbol{\varPi}* \boldsymbol{M} (t)dt . 
 \eeqnn
 Plugging this back into (\ref{eqn.5018}) and the differentiating both sides with respect to $t$ along with the right-continuity of integrands, we have almost surely,
 \beqnn
 \boldsymbol{\varXi} = \boldsymbol{\varUpsilon}
 + \boldsymbol{\varPi}* \boldsymbol{M}. 
 \eeqnn
Here we have finished the whole proof.
\qed 

 \medskip
 
 \textbf{\textit{Proof of Theorem~\ref{MainThm.J1WeakConvergence}.}} 
 The weak convergence (\ref{MainWeakConv}) in $D_{J_1}(\mathbb{R}_+;\mathbb{R}_+^d\times\mathbb{R}_+^d\times\mathbb{R}^d)$ follows directly Theorem~3.37 in \cite[p.354]{JacodShiryaev2003} and Theorem~\ref{MainThm.S-WeakConvergence}(2.i). 
 We now prove claim (2). 
 By (\ref{eqn.QV}) and claim (1), the martingale $\widetilde{\boldsymbol{N}}^{(n)}$ has predictable  quadratic co-variation
 \beqnn
 \Big\langle \widetilde{N}_i^{(n)}, \widetilde{N}_j^{(n)} \Big\rangle = \mathbf{1}_{\{ i=j \}} \cdot \mathcal{I}_{\varLambda_i^{(n)}},\quad i,j\in\mathtt{D}, 
 \eeqnn
 which converges to $\mathbf{1}_{\{ i=j \}} \cdot \varXi_i$ weakly in $D_{J_1}(\mathbb{R}_+;\mathbb{R}_+)$   as $n\to\infty$. 
 Hence the continuous limit martingale $\boldsymbol{M}$ has  the predictable quadratic co-variation 
 \beqnn
 \langle M_i,M_j\rangle=\mathbf{1}_{\{ i=j \}} \cdot \varXi_i,
 \eeqnn
 By Theorem~5.13 in \cite[p.121]{LeGall2016}, there exists a standard $d$-dimensional Brownian motion $\boldsymbol{B}$ such that the continuous martingale $\boldsymbol{M}$ can be represented as
 \beqnn
 \boldsymbol{M}(t)=\big(M_i(t)\big)_{i\in\mathtt{D}}= \big(B_i\circ \varXi_i(t)\big)_{i\in\mathtt{D}} = \big(B_i(\varXi_i(t)) \big)_{i\in\mathtt{D}},\quad t\geq 0.  
 \eeqnn
 Putting this together with (\ref{Eqn.HawkesVoletrra}), we see that $(\boldsymbol{\varXi},\boldsymbol{M},\boldsymbol{B})$ solves (\ref{Eqn.HawkesVoletrra.02}).   
 
 We now prove the weak uniqueness. Assume that $(\boldsymbol{\varXi},\boldsymbol{M})$ and $(\boldsymbol{\varXi}^*,\boldsymbol{M}^*)$ are two  weak solutions of (\ref{Eqn.HawkesVoletrra.02}).
 Similarly as in the proof of Theorem~\ref{MainThm.S-WeakConvergence}, it suffices to prove that 
 \beqnn
 \mathbf{E}\Big[ \exp \big\{ \mathcal{I}_{\boldsymbol{f}} * d\boldsymbol{\varXi} (T)  + \mathcal{I}_{\boldsymbol{h}}* d\boldsymbol{M} (T) \big\} \Big]
 = \mathbf{E}\Big[ \exp \big\{ \mathcal{I}_{\boldsymbol{f}} * d\boldsymbol{\varXi}^* (T)+ \mathcal{I}_{\boldsymbol{h}}* d\boldsymbol{M}^* (T) \big\} \Big].
 \eeqnn
 for any $\boldsymbol{f}\in C(\mathbb{R}_+;\mathbb{C}_-^d)$ and $\boldsymbol{h}\in C(\mathbb{R}_+;\mathtt{i}\mathbb{R}^d)$, which follows directly from Lemma~\ref{Lemma.FLFXiM} and the uniqueness of solutions to (\ref{Eqn.AffineVolterra}); see Lemma~\ref{Lemma.ExistUniqueV}. 
 \qed  
 
  \textbf{\textit{Proof of Corollary~\ref{Corollary.ContinuityXiM}.}} When $\boldsymbol{\varPi}(dt)$ has density $\boldsymbol{\pi}$, we can write (\ref{Eqn.HawkesVoletrra}) as 
  \beqnn
   \boldsymbol{\varXi}= \boldsymbol{\varUpsilon}+ \boldsymbol{\pi}*\boldsymbol{M}.
  \eeqnn
  The continuity of $  \boldsymbol{\varXi}$ follows directly from that of $ \boldsymbol{\varUpsilon}$ and $\boldsymbol{\pi}*\boldsymbol{M}$; see Proposition~\ref{Prop.A1}.   
  For any  $\boldsymbol{f}\in L^\infty_{\rm loc}(\mathbb{R}_+;\mathbb{C}_-^d)$ and $\boldsymbol{h}\in L^\infty_{\rm loc}(\mathbb{R}_+;\mathtt{i}\mathbb{R}^d)$, all claims in Theorem~\ref{MainThm.FourierLaplaceFunctional}  hold because of Corollary~\ref{Corollary.407} and Lemma~\ref{Lemma.FLFXiM}.
  \qed 
  
  \medskip
  
%
%
%
%
  
   \textit{\textbf{Proof of Corollary~\ref{Corollary.2001}.}}  Lemma 2.1 in \cite{Jaber2021} along with our assumptions induces that $\boldsymbol{\varXi} $ is differentiable with intensity satisfying the first equation in \eqref{eqn.SVE04}. 
   It remains to prove the martingale representation of $M$.   
   By Theorem~\ref{MainThm.J1WeakConvergence}(2) and  Theorem~7.1 in \cite[p.84]{IkedaWatanabe1989}, it suffices to identify that $\boldsymbol{\xi}$ is predictable, which unfortunately is not clear. Therefore we need to find a suitable predictable modification of $ \boldsymbol{\xi}$ satisfying the preceding stochastic equation. 
   For convention, we set $\boldsymbol{\xi}(t)
   =\boldsymbol{0}$ for $t<0$.
   Consider the following limit inferior
   \beqnn
    \boldsymbol{\hat\xi}(t):= \liminf_{\delta\to 0+} \boldsymbol{\hat\xi}_\delta(t) \in[0,\infty]
    \quad \mbox{with}\quad  \boldsymbol{\hat\xi}_\delta(t):=\frac{1}{\delta}\int_{t-\delta}^{t}  \boldsymbol{\xi}(s)ds ,\quad t\geq 0.
   \eeqnn
   There exists a null set $\Omega_0\subset\Omega$ such that $\boldsymbol{\hat\xi}(t,\omega)<\infty$ for almost every $t\geq 0$ and all $\omega\in \Omega$, i.e.,  for each $T\geq 0$, by Fatou's lemma and then Fubini's theorem, 
   \beqnn
   \int_0^T\boldsymbol{\hat\xi}(t)dt
   \ar=\ar \int_0^T \liminf_{\delta\to 0+} \frac{1}{\delta}\int_{(t-\delta)^+}^{t}  \boldsymbol{\xi}(s)ds dt\cr
   \ar\leq\ar \liminf_{\delta\to 0+}\int_0^T  \frac{1}{\delta}\int_{t-\delta}^{t}  \boldsymbol{\xi}(s)ds dt\cr
   \ar\leq\ar \liminf_{\delta\to 0+}\int_0^T      \boldsymbol{\xi}(s)ds\cdot \frac{1}{\delta}\int_s^{s+\delta}dt
   = \int_0^T      \boldsymbol{\xi}(s)ds<\infty,\quad a.s.
   \eeqnn
   Additionally, the Lebesgue differentiation theorem induces that $\boldsymbol{\hat\xi}(t,\omega)=\boldsymbol{\xi}(t,\omega)$ for almost every $t\geq 0$ and  
   \beqnn
   \int_0^t \boldsymbol{\hat\xi}(s,\omega)ds=\int_0^t\boldsymbol{\xi}(s,\omega)ds,\quad t\geq 0.
   \eeqnn
   Consequently, the first equation in \eqref{eqn.SVE04} still holds with $\boldsymbol{\xi}$ replaced by $\boldsymbol{\hat\xi}$. 
   Finally, note that the process $\boldsymbol{\hat\xi}_\delta$ is continuous and hence predictable. 
   This induces that the limit process $\boldsymbol{\hat\xi}$ is also predictable. 
   \qed
  
  \medskip
  
  \textit{\textbf{Proof of Theorem~\ref{Thm.ConvergenceDS01}.}} 
  We prove the three claims separately in the following.
  
  {\it Claim (1).} 
  By (\ref{eqn.224}), the function $\boldsymbol{\varPhi}_0$ is $\boldsymbol{K}$-admissible with auxiliary matrices $\boldsymbol{Q},\boldsymbol{U} \in \mathbb{R}^{d\times d}$ and an auxiliary constant $\lambda^+_0$.  
  Then claim (1) follows directly from Lemma~\ref{Lemma.Pi}.
  
  {\it Claim (2).} We just prove the first equality in (\ref{eqn.ResolventEquation}). 
  The second one can be proved in the same way. 
  It suffices to prove that 
  \beqlb\label{eqn.5012}
  \mathcal{L}_{\boldsymbol{\varPi}_0}(\lambda)
  =\mathcal{L}_{\boldsymbol{\varPi}}(\lambda)
  + \mathcal{L}_{\boldsymbol{\varPi}_0}(\lambda) \cdot \boldsymbol{b}^{\boldsymbol{\varPhi}}\cdot \mathcal{L}_{\boldsymbol{\varPi} }(\lambda)
  = \big( \mathbf{Id} + \mathcal{L}_{\boldsymbol{\varPi}_0}(\lambda) \cdot \boldsymbol{b}^{\boldsymbol{\varPhi}} \big)\cdot \mathcal{L}_{\boldsymbol{\varPi}}(\lambda),
  \quad \lambda \geq \lambda^+\vee \lambda^+_0. 
  \eeqlb
  From (\ref{eqn.LimitResolvent}), we have 
  \beqnn
  \mathcal{L}_{\boldsymbol{\varPi}}(\lambda)
  \ar=\ar \boldsymbol{Q}
  \begin{pmatrix}
  	\boldsymbol{\varphi}(\lambda)  & 
  	\mathbf{0} \vspace{7pt} \\
  	\mathbf{0} & \mathbf{Id}
  \end{pmatrix}^{-1}
  \begin{pmatrix}
  	\mathbf{Id} & 
  	\boldsymbol{U}_{\mathtt{I}\mathtt{J}}\big(\mathbf{Id}-\boldsymbol{U}_{\mathtt{J}\mathtt{J}} \big)^{-1} \vspace{7pt} \\
  	\mathbf{0} & \mathbf{0}
  \end{pmatrix}
  \boldsymbol{Q}^{\rm T},\cr
  \mathcal{L}_{\boldsymbol{\varPi}_0}(\lambda)
  \ar=\ar \boldsymbol{Q}
  \begin{pmatrix}
  	\boldsymbol{\varphi}_0(\lambda)  & 
  	\mathbf{0} \vspace{7pt} \\
  	\mathbf{0} & \mathbf{Id}
  \end{pmatrix}^{-1}
  \begin{pmatrix}
  	\mathbf{Id} & 
  	\boldsymbol{U}_{\mathtt{I}\mathtt{J}}\big(\mathbf{Id}-\boldsymbol{U}_{\mathtt{J}\mathtt{J}} \big)^{-1} \vspace{7pt} \\
  	\mathbf{0} & \mathbf{0}
  \end{pmatrix}
  \boldsymbol{Q}^{\rm T},
  \eeqnn
  and hence the matrix $\big( \mathbf{Id} + \mathcal{L}_{\boldsymbol{\varPi}_0}(\lambda) \cdot \boldsymbol{b}^{\boldsymbol{\varPhi}} \big)\cdot \mathcal{L}_{\boldsymbol{\varPi}}(\lambda)$ can be written as 
  \beqnn
  \ar\ar\Bigg(\mathbf{Id} +  \boldsymbol{Q}
  \begin{pmatrix}
  	\boldsymbol{\varphi}_0(\lambda)  & 
  	\mathbf{0} \vspace{3pt} \\
  	\mathbf{0} & \mathbf{Id}
  \end{pmatrix}^{-1}
  \begin{pmatrix}
  	\mathbf{Id} & 
  	\boldsymbol{U}_{\mathtt{I}\mathtt{J}}\big(\mathbf{Id}-\boldsymbol{U}_{\mathtt{J}\mathtt{J}} \big)^{-1} \vspace{3pt} \\
  	\mathbf{0} & \mathbf{0}
  \end{pmatrix}
  \boldsymbol{Q}^{\rm T} \cdot \boldsymbol{b}^{\boldsymbol{\varPhi}} \Bigg)\cr
  \ar\ar \qquad \times \boldsymbol{Q}
  \begin{pmatrix}
  	\boldsymbol{\varphi}(\lambda)  & 
  	\mathbf{0} \vspace{3pt} \\
  	\mathbf{0} & \mathbf{Id}
  \end{pmatrix}^{-1}
  \begin{pmatrix}
  	\mathbf{Id} & 
  	\boldsymbol{U}_{\mathtt{I}\mathtt{J}}\big(\mathbf{Id}-\boldsymbol{U}_{\mathtt{J}\mathtt{J}} \big)^{-1} \vspace{3pt} \\
  	\mathbf{0} & \mathbf{0}
  \end{pmatrix}
  \boldsymbol{Q}^{\rm T}\cr
  \ar=\ar  \boldsymbol{Q}\begin{pmatrix}
  	\boldsymbol{\varphi}_0(\lambda)  & 
  	\mathbf{0} \vspace{3pt} \\
  	\mathbf{0} & \mathbf{Id}
  \end{pmatrix}^{-1}\Bigg(
  \begin{pmatrix}
  	\boldsymbol{\varphi}_0(\lambda)  & 
  	\mathbf{0} \vspace{3pt} \\
  	\mathbf{0} & \mathbf{Id}
  \end{pmatrix}
  +
  \begin{pmatrix}
  	\mathbf{Id} & 
  	\boldsymbol{U}_{\mathtt{I}\mathtt{J}}\big(\mathbf{Id}-\boldsymbol{U}_{\mathtt{J}\mathtt{J}} \big)^{-1} \vspace{3pt} \\
  	\mathbf{0} & \mathbf{0}
  \end{pmatrix}
  \cdot \boldsymbol{Q}^{\rm T}  \boldsymbol{b}^{\boldsymbol{\varPhi}} \boldsymbol{Q}\Bigg)\cr
  \ar\ar \qquad \times 
  \begin{pmatrix}
  	\boldsymbol{\varphi}(\lambda)  & 
  	\mathbf{0} \vspace{3pt} \\
  	\mathbf{0} & \mathbf{Id}
  \end{pmatrix}^{-1}
  \begin{pmatrix}
  	\mathbf{Id} & 
  	\boldsymbol{U}_{\mathtt{I}\mathtt{J}}\big(\mathbf{Id}-\boldsymbol{U}_{\mathtt{J}\mathtt{J}} \big)^{-1} \vspace{3pt} \\
  	\mathbf{0} & \mathbf{0}
  \end{pmatrix}
  \boldsymbol{Q}^{\rm T} . 
  \eeqnn
  By (\ref{eqn.224}) with $\lambda=0$ and (\ref{eqn.225}), we have
  \beqnn
  \begin{pmatrix}
  	\mathbf{Id} & 
  	\boldsymbol{U}_{\mathtt{I}\mathtt{J}}\big(\mathbf{Id}-\boldsymbol{U}_{\mathtt{J}\mathtt{J}} \big)^{-1} \vspace{3pt} \\
  	\mathbf{0} & \mathbf{0}
  \end{pmatrix}
  \cdot \boldsymbol{Q}^{\rm T}  \boldsymbol{b}^{\boldsymbol{\varPhi}} \boldsymbol{Q} = \begin{pmatrix}
  	\boldsymbol{\varphi}_0(0) &  \mathbf{0} \vspace{3pt} \\
  	\mathbf{0} & \mathbf{0}
  \end{pmatrix}
  \eeqnn
  and then
  \beqnn
  \big( \mathbf{Id} + \mathcal{L}_{d\boldsymbol{\varPi}_0}(\lambda) \cdot \boldsymbol{b}^{\boldsymbol{\varPhi}} \big)\cdot \mathcal{L}_{d\boldsymbol{\varPi}}(\lambda)
  \ar=\ar  \boldsymbol{Q}\begin{pmatrix}
  	\boldsymbol{\varphi}_0(\lambda)  & 
  	\mathbf{0} \vspace{3pt} \\
  	\mathbf{0} & \mathbf{Id}
  \end{pmatrix}^{-1}\Bigg(
  \begin{pmatrix}
  	\boldsymbol{\varphi}_0(\lambda)  & 
  	\mathbf{0} \vspace{3pt} \\
  	\mathbf{0} & \mathbf{Id}
  \end{pmatrix}
  +
  \begin{pmatrix}
  	\boldsymbol{\varphi}_0(0) &  \mathbf{0} \vspace{3pt} \\
  	\mathbf{0} & \mathbf{0}
  \end{pmatrix}\Bigg)\cr
  \ar\ar \qquad \times 
  \begin{pmatrix}
  	\boldsymbol{\varphi}(\lambda)  & 
  	\mathbf{0} \vspace{3pt} \\
  	\mathbf{0} & \mathbf{Id}
  \end{pmatrix}^{-1}
  \begin{pmatrix}
  	\mathbf{Id} & 
  	\boldsymbol{U}_{\mathtt{I}\mathtt{J}}\big(\mathbf{Id}-\boldsymbol{U}_{\mathtt{J}\mathtt{J}} \big)^{-1} \vspace{3pt} \\
  	\mathbf{0} & \mathbf{0}
  \end{pmatrix}
  \boldsymbol{Q}^{\rm T}\cr
  \ar=\ar  \boldsymbol{Q}\begin{pmatrix}
  	\boldsymbol{\varphi}_0(\lambda)  & 
  	\mathbf{0} \vspace{3pt} \\
  	\mathbf{0} & \mathbf{Id}
  \end{pmatrix}^{-1} 
  \begin{pmatrix}
  	\mathbf{Id} & 
  	\boldsymbol{U}_{\mathtt{I}\mathtt{J}}\big(\mathbf{Id}-\boldsymbol{U}_{\mathtt{J}\mathtt{J}} \big)^{-1} \vspace{3pt} \\
  	\mathbf{0} & \mathbf{0}
  \end{pmatrix}
  \boldsymbol{Q}^{\rm T},
  \eeqnn
  which equals to $\mathcal{L}_{\boldsymbol{\varPi}_0}(\lambda)$ and hence  (\ref{eqn.5012})  holds.

  {\it Claim (3).}  We first prove that any solution $(\boldsymbol{\varXi},\boldsymbol{M})$ of  (\ref{Eqn.HawkesVoletrra}) solves (\ref{eqn.SVE0101}). 
  Plugging  the sum on the right-hand side of (\ref{Eqn.HawkesVoletrra}) into  $\boldsymbol{\varPi}_0 * \big(\boldsymbol{b}^{\boldsymbol{\varPhi}}\cdot \boldsymbol{\varXi} \big)$ and then using Fubini's theorem, 
  \beqnn
  \boldsymbol{\varPi}_0 * \big(\boldsymbol{b}^{\boldsymbol{\varPhi}}\cdot \boldsymbol{\varXi} \big)
  \ar=\ar \boldsymbol{\varPi}_0 * \big(\boldsymbol{b}^{\boldsymbol{\varPhi}}\cdot \boldsymbol{\varPi}  * \boldsymbol{\varGamma}   + \boldsymbol{b}^{\boldsymbol{\varPhi}}\cdot \boldsymbol{\varPi}  * \boldsymbol{M} \big) .
  \eeqnn
  Taking this back into the right side of (\ref{eqn.SVE0101}), then using Fubini's theorem and (\ref{eqn.ResolventEquation}), we have 
  \beqnn
  \lefteqn{\boldsymbol{\varPi}_0  * \boldsymbol{\varGamma} - \boldsymbol{\varPi}_0 * \big(\boldsymbol{b}^{\boldsymbol{\varPhi}}\cdot \boldsymbol{\varXi} \big) + \boldsymbol{\varPi}_0  * \boldsymbol{M}}\ar\ar\cr
  \ar\ar\cr
  \ar=\ar \boldsymbol{\varPi}_0  * \boldsymbol{\varGamma} -\boldsymbol{\varPi}_0 *  \boldsymbol{b}^{\boldsymbol{\varPhi}}\cdot \boldsymbol{\varPi}  * \boldsymbol{\varGamma}  +  \boldsymbol{\varPi}_0  * \boldsymbol{M} - \boldsymbol{\varPi}_0 *   \boldsymbol{b}^{\boldsymbol{\varPhi}}\cdot \boldsymbol{\varPi}  * \boldsymbol{M}   
  = \boldsymbol{\varPi}  * \boldsymbol{\varGamma} + \boldsymbol{\varPi}  * \boldsymbol{M} = \boldsymbol{\varXi},
  \eeqnn
  which induces that $(\boldsymbol{\varXi},\boldsymbol{M})$ satisfies (\ref{eqn.SVE0101}).
  
  For the converse, assume that $(\boldsymbol{\varXi},\boldsymbol{M})$ solves (\ref{eqn.SVE0101}) and then multiply both sides of (\ref{eqn.SVE0101}) and (\ref{eqn.ResolventEquation}) by the matrix $\boldsymbol{b}^{\boldsymbol{\varPhi}}$,  we have 
  \beqlb\label{eqn.631}
  \boldsymbol{b}^{\boldsymbol{\varPhi}} \cdot \boldsymbol{\varXi} \ar=\ar  \boldsymbol{b}^{\boldsymbol{\varPhi}} \cdot \boldsymbol{\varPi}_0  * \boldsymbol{\varGamma} - \big(\boldsymbol{b}^{\boldsymbol{\varPhi}} \cdot \boldsymbol{\varPi}_0\big) * \big(\boldsymbol{b}^{\boldsymbol{\varPhi}}\cdot \boldsymbol{\varXi} \big) + \big(\boldsymbol{b}^{\boldsymbol{\varPhi}} \cdot \boldsymbol{\varPi}_0 \big) * \boldsymbol{M} 
  \eeqlb
  and 
  \beqlb\label{eqn.632}
  \boldsymbol{b}^{\boldsymbol{\varPhi}} \cdot	\boldsymbol{\varPi}_0 
  = \boldsymbol{b}^{\boldsymbol{\varPhi}} \cdot \boldsymbol{\varPi}
  + \big(\boldsymbol{b}^{\boldsymbol{\varPhi}} \cdot \boldsymbol{\varPi}_0\big) * \big(\boldsymbol{b}^{\boldsymbol{\varPhi}}\cdot \boldsymbol{\varPi} \big)  ,
  \eeqlb
 which shows that $-\boldsymbol{b}^{\boldsymbol{\varPhi}} \cdot \boldsymbol{\varPi}$ is the resolvent associated to $-\boldsymbol{b}^{\boldsymbol{\varPhi}} \cdot \boldsymbol{\varPi}_0$.
 In view of Proposition~\ref{Prop.A3}, we can write (\ref{eqn.631}) into 
  \beqnn
  \boldsymbol{b}^{\boldsymbol{\varPhi}} \cdot \boldsymbol{\varXi} 
  \ar=\ar  \boldsymbol{b}^{\boldsymbol{\varPhi}} \cdot \boldsymbol{\varPi}_0  * \boldsymbol{\varGamma} - \big(\boldsymbol{b}^{\boldsymbol{\varPhi}} \cdot \boldsymbol{\varPi}\big) * \big(\boldsymbol{b}^{\boldsymbol{\varPhi}} \cdot \boldsymbol{\varPi}_0 \big) * \boldsymbol{\varGamma}  + \big(\boldsymbol{b}^{\boldsymbol{\varPhi}} \cdot \boldsymbol{\varPi}_0 \big) * \boldsymbol{M}  -\big(\boldsymbol{b}^{\boldsymbol{\varPhi}} \cdot \boldsymbol{\varPi}\big) * \big(\boldsymbol{b}^{\boldsymbol{\varPhi}} \cdot \boldsymbol{\varPi}_0 \big) * \boldsymbol{M} .
  \eeqnn
  Applying Fubini's theorem and then (\ref{eqn.632}) to the right side of this equation, we have 
  \beqnn
  \boldsymbol{b}^{\boldsymbol{\varPhi}} \cdot \boldsymbol{\varXi} 
  =\big(\boldsymbol{b}^{\boldsymbol{\varPhi}} \cdot \boldsymbol{\varPi} \big) * \boldsymbol{a} + \big(\boldsymbol{b}^{\boldsymbol{\varPhi}} \cdot \boldsymbol{\varPi}\big) * \boldsymbol{M}.
  \eeqnn
  Plugging this into (\ref{eqn.SVE0101}) and then using Fubini's theorem and (\ref{eqn.ResolventEquation}), 
  \beqnn
  \boldsymbol{\varXi} \ar=\ar   \boldsymbol{\varPi}_0  * \boldsymbol{\varGamma} - \boldsymbol{\varPi}_0  *\big(\boldsymbol{b}^{\boldsymbol{\varPhi}} \cdot \boldsymbol{\varPi} \big) * \boldsymbol{\varGamma} -  \boldsymbol{\varPi}_0  *\big(\boldsymbol{b}^{\boldsymbol{\varPhi}} \cdot \boldsymbol{\varPi}\big) * \boldsymbol{M}  + \boldsymbol{\varPi}_0  * \boldsymbol{M} 
  = \boldsymbol{\varPi}  * \boldsymbol{\varGamma} + \boldsymbol{\varPi} * \boldsymbol{M} ,
  \eeqnn
  which tells that $(\boldsymbol{\varXi},\boldsymbol{M})$ is a solution of (\ref{Eqn.HawkesVoletrra}).  
  The proof is end. 
  \qed

 \appendix
 
 \renewcommand{\theequation}{A.\arabic{equation}}
  \section{Basic properties of matrices}
 \label{Sec.Matrices}
 \setcounter{equation}{0}
 
 In this section, we recall some elementary properties of (non-negative) $d\times d$-matrices. The reader may refer to \cite{BermanPlemmons1994,Seneta1981} for more details. 

 \begin{proposition}\label{Appendix.Prop.BlockInversion}
 	For a matrix $\boldsymbol{A} \in  \mathbb{R}^{d\times d}$ that is partitioned into four blocks, i.e.,
 	\beqnn
 	\boldsymbol{A} =  \begin{pmatrix}
 		\boldsymbol{A}_{11} & \boldsymbol{A}_{12}\\
 		\boldsymbol{A}_{21} & \boldsymbol{A}_{11}
 	\end{pmatrix},
 	\eeqnn
  the following hold.

 	\begin{enumerate}
 		\item[(1)]  We have ${\rm det}(	\boldsymbol{A} ) = {\rm det}(	\boldsymbol{A}_{11} \boldsymbol{A}_{22}-\boldsymbol{A}_{12}\boldsymbol{A}_{21} )$.
 		
 		\item[(2)] If both $\boldsymbol{A}_{11}$ and $\boldsymbol{A}_{22}-\boldsymbol{A}_{21}\boldsymbol{A}_{11}^{-1}\boldsymbol{A}_{12}$ are invertible, then $\boldsymbol{A}$ is invertible if and only if the following well-defined matrix is invertible
 		\beqnn
 		\begin{pmatrix}
 			\big(\boldsymbol{A}_{22}-\boldsymbol{A}_{21}\boldsymbol{A}_{11}^{-1}\boldsymbol{A}_{12}\big)^{-1} & -\big(\boldsymbol{A}_{22}-\boldsymbol{A}_{21}\boldsymbol{A}_{11}^{-1}\boldsymbol{A}_{12}\big)^{-1}\boldsymbol{A}_{12}\boldsymbol{A}_{22}^{-1}\\
 			-\boldsymbol{A}_{22}^{-1}\boldsymbol{A}_{21} \big(\boldsymbol{A}_{22}-\boldsymbol{A}_{21}\boldsymbol{A}_{11}^{-1}\boldsymbol{A}_{12}\big)^{-1} & \boldsymbol{A}_{22}^{-1}\big( \mathbf{Id} + \boldsymbol{A}_{21} \big(\boldsymbol{A}_{22}-\boldsymbol{A}_{21}\boldsymbol{A}_{11}^{-1}\boldsymbol{A}_{12}\big)^{-1} \boldsymbol{A}_{12}\boldsymbol{A}_{22}^{-1} \big)
 		\end{pmatrix},
 		\eeqnn
 		which, in this case, is also the inverse of $\boldsymbol{A}$.
  	\end{enumerate}
 	 
 \end{proposition}

 \begin{proposition}\label{Appendix.Prop.LimitMatrix}
  For a sequence of matrices $\{\boldsymbol{A}_n\}_{n\geq 1} \subset \mathbb{R}^{d\times d}$, if the element-wise limit $\lim_{n\to\infty} \boldsymbol{A}_n = \boldsymbol{A} \in  \mathbb{R}^{d\times d}$ exists as $n\to\infty$, the following hold.
  \begin{enumerate}
  	\item[(1)] We have $\rho (\boldsymbol{A}_n) \to \rho (\boldsymbol{A})$  and ${\rm det} (\boldsymbol{A}_n) \to {\rm det} (\boldsymbol{A})$ as $n\to\infty$;
  	
  	\item[(2)] If $\boldsymbol{A}$ is invertible, i.e. ${\rm det} (\boldsymbol{A})\neq 0$, then there exists a constant $n_0\geq 1$ such that ${\rm det} (\boldsymbol{A}_n)\neq 0$ and hence $\boldsymbol{A}_n$ is invertible for any $n\geq n_0$.
  \end{enumerate}

 \end{proposition}

 \begin{proposition}[Perron-Frobenius’ Theorem] \label{Appendix.Prop.PFT}
 	For each $\boldsymbol{A} \in \mathbb{R}^{d\times d}_+$, we have $\rho(\boldsymbol{A})$ is an (largest) eigenvalue of $\boldsymbol{A}$. Moreover, there exists a non-negative eigenvector associated to the eigenvalue $\rho (A)$. 
 \end{proposition}
 
 \begin{proposition}[Real Schur decomposition]\label{Appendix.Prop.RSD}
 	For each $\boldsymbol{A} \in \mathbb{R}_{+,r,\kappa}^{d\times d}$ with $r>0$ and $\kappa\in\mathtt{D}$, there exists a real orthogonal  matrix $\boldsymbol{O} \in \mathbb{R}^{k\times k}$ (i.e. $\boldsymbol{O}^{\rm T}\boldsymbol{O}=\boldsymbol{O}\boldsymbol{O}^{\rm T} = \mathbf{Id}$) such that 
 	\beqnn 
 	\boldsymbol{A}=\boldsymbol{O}^{\rm T} \boldsymbol{B}\boldsymbol{O}
 	\quad \mbox{and}\quad
 	\boldsymbol{B} = \begin{pmatrix}
 		\boldsymbol{B}_{\mathtt{I}\mathtt{I}} & \boldsymbol{B}_{\mathtt{I}\mathtt{J}}\\
 		\boldsymbol{0} & \boldsymbol{B}_{\mathtt{J}\mathtt{J}}
 	\end{pmatrix},
 	\eeqnn 
   where $\boldsymbol{B}_{\mathtt{I}\mathtt{I}} \in \mathbb{R}^{\kappa\times\kappa}$ is an upper triangular matrix with all diagonal entries equal to $r$ and all eigenvalues of  $\boldsymbol{B}_{\mathtt{J}\mathtt{J}} \in \mathbb{R}^{(d-\kappa)\times(d-\kappa)}$ have real parts less than $r$.  
  \end{proposition}

  \renewcommand{\theequation}{B.\arabic{equation}}
  
 \section{Convolution and Volterra equation}
 \label{Appendix-Volterra}
 \setcounter{equation}{0}
  
 We recall some elements of convolutions and Volterra equations; see \cite{GripenbergLondenStaffans1990} for details. The next proposition is a direct corollary of Young's convolution inequality and H\"older's inequality. 
 
 \begin{proposition}\label{Prop.A1}
 	For each $f\in L^p_{\rm loc}(\mathbb{R}_+;\mathbb{C})$ and $g\in L^q_{\rm loc}(\mathbb{R}_+;\mathbb{C})$ with $p,q\in [1,\infty]$, we have the following.
 	\begin{enumerate}
 		\item[(1)] If   $\frac{1}{p}+\frac{1}{q}\leq1$, then $f*g \in C(\mathbb{R}_+;\mathbb{C})$.
 		
 		\item[(2)] If $\frac{1}{p}+\frac{1}{q}>1$, then $	\|f*g\|_{L^r_T} \leq \|f\|_{L^p_T} \cdot \|g\|_{L^q_T}$ for any $T\geq 0$ and $r\in[1,\infty]$ with $\frac{1}{r}\geq  \frac{1}{p}+\frac{1}{q}-1$. 
 		 
 	\end{enumerate}

 \end{proposition}
 
 Fix a $\sigma$-finite measure $\boldsymbol{\upsilon}(dt) \in M_{\rm}(\mathbb{R}_+;\mathbb{R}_+^{d\times d})$ as well as a column vector function  $\boldsymbol{f}: \mathbb{R}_+\mapsto \mathbb{C}^d $ and a row vector function $\boldsymbol{g}: \mathbb{C}^d\mapsto \mathbb{C}^d$, consider the following nonlinear Volterra equation with measure kernel
 \beqlb\label{eqn.Volterra}
 \boldsymbol{x}(t)= \boldsymbol{f}(t) + \int_{[0,t]} \boldsymbol{g}\big(t-s, \boldsymbol{x}(t-s)\big)  \boldsymbol{\upsilon}(ds),\quad t\geq 0.
 \eeqlb 
 For $p\in[1,\infty]$, a pair $(T_\infty,\boldsymbol{x})$ with $T_\infty\in (0,\infty]$ and $\boldsymbol{x}\in L^p_{\rm loc}([0,T_\infty),\mathbb{C}^d)$ is called a \textit{non-continuable solution} of (\ref{eqn.Volterra}) if $\boldsymbol{x}$  satisfies (\ref{eqn.Volterra}) on $[0,T_\infty)$ and $\|\boldsymbol{x}\|_{L^p_{T_\infty}}=\infty$ if $T_\infty<\infty$; we say that $\boldsymbol{x}$ is a {\it global solution} if $T_\infty=\infty$.  
 The next proposition comes from Corollary~4.3 in \cite[p.364]{GripenbergLondenStaffans1990}.

 \begin{proposition}\label{Prop.A2}
 Suppose that $\boldsymbol{\upsilon}(dt)$ is absolutely continuous with respect to Lebesgue measure and $\boldsymbol{g}$ is locally Lipschitz continuous with respect to the second argument. For $\boldsymbol{f}\in L^\infty_{\rm loc}(\mathbb{R}_+,\mathbb{C}^d)$, the nonlinear Volterra equation (\ref{eqn.Volterra}) has a unique  non-continuable solution  $(T_\infty, \boldsymbol{ x}) \in (0,\infty)\times L^\infty_{\rm loc}([0,T_\infty),\mathbb{C}^d)$. 
 Moreover, if $\boldsymbol{f}$ is continuous on $[0,T_\infty)$, then so is $\boldsymbol{ x}$. 
 \end{proposition} 
 
 When (\ref{eqn.Volterra}) is linear, i.e. $\boldsymbol{g}(t,\boldsymbol{z})=\boldsymbol{z}$ for $\boldsymbol{z}\in\mathbb{C}^d$, then existence and uniqueness of solutions have been well-established; see Chapter 2 and 4 in \cite{GripenbergLondenStaffans1990}. 
  
 \begin{proposition}\label{Prop.A3}
  If $\boldsymbol{\upsilon} $ has density  $\boldsymbol{\upsilon}'$, there exists a unique solution $\boldsymbol{R}_{\boldsymbol{\upsilon}'} \in  L^1_{\rm loc}(\mathbb{R}_+;\mathbb{R}_+^{d\times d})$ to 
  \beqnn
  \boldsymbol{R}_{\boldsymbol{\upsilon}'} = \boldsymbol{\upsilon}'+  \boldsymbol{R}_{\boldsymbol{\upsilon}'} * \boldsymbol{\upsilon}'.
  \eeqnn
 Moreover, if $\boldsymbol{g}(t,\boldsymbol{z})=\boldsymbol{z}$ and $\boldsymbol{f}\in L^p_{\rm loc}(\mathbb{R}_+,\mathbb{C}^d)$ [resp. $C(\mathbb{R}_+;\mathbb{C}^d)$], the linear Volterra equation (\ref{eqn.Volterra})  has a unique global solution in $L^p_{\rm loc}(\mathbb{R}_+,\mathbb{C}^d)$ [resp. $C(\mathbb{R}_+;\mathbb{C}^d)$]
 and the solution is given by 
 \beqnn
 \boldsymbol{x}=\boldsymbol{f} + \boldsymbol{R}_{\boldsymbol{\upsilon}'}*\boldsymbol{f} .
 \eeqnn  
 \end{proposition}
   
    \renewcommand{\theequation}{C.\arabic{equation}}
    \section{$S$-topology} 
 \label{Sec.TopologiesD}
 \setcounter{equation}{0}
 
  The $S$-topology was firstly introduced by Jakubowski \cite{Jakubowski1997} on $D([0,T];\mathbb{R}^d)$ for $T\geq 0$ as a sequential topology generated by naturally arising criteria of relative compactness; see \cite{Jakubowski1996,Jakubowski1997} for more detailed explanations and results.
  It can be extended to the space  $D(\mathbb{R}_+;\mathbb{R}^d)$ following the standard argument. For simplicity, we state all results for functions/processes in $D(\mathbb{R}_+;\mathbb{R})$, which can be easily extended into the multi-dimensional case. Their proofs are given here or can be found in \cite{Jakubowski1996}. 
 
  Let $V(\mathbb{R}_+;\mathbb{R})$ be the space of all c\`adl\`ag functions with bounded variation. 
  Every element of $V(\mathbb{R}_+;\mathbb{R})$ determines a $\mathbb{R}$-valued signed measure on $[0,T]$.  
  For a sequence $\{f_n\}_{n\geq 1}\subset D(\mathbb{R}_+;\mathbb{R})$, we write $f_n \overset{\rm S}\to f \in D(\mathbb{R}_+;\mathbb{R})$ as $n\to\infty$ if for each $T\geq 0$ and  $\epsilon >0$, there exist a function $f^\epsilon \in V(\mathbb{R}_+;\mathbb{R})$ a sequence of functions $\{f_n^\epsilon\}_{n\geq 1}\subset V(\mathbb{R}_+;\mathbb{R})$ such that
  \beqlb\label{eqn.D01}
  \big\|f_n-f^\epsilon_n\big\|_{L^\infty_T} \leq \epsilon,\quad  \big\|f-f^\epsilon\big\|_{L^\infty_T} \leq \epsilon
  \quad \mbox{and}\quad
  \lim_{n\to\infty}\int_0^T g(s) df_n^\epsilon  = \int_0^T g(s) df^\epsilon,\quad T\geq 0, 
  \eeqlb
  for any continuous function $g$ on $\mathbb{R}_+$.
 The topology on $ D(\mathbb{R}_+;\mathbb{R})$ induced by the  $\overset{\rm S}\to$ convergence is called \textit{$S$-topology} and the corresponding topological space is denoted as $D_S(\mathbb{R}_+;\mathbb{R})$. 
 As a direct consequence, the next proposition show that the $\overset{\rm S}\to$ convergence  induces the almost everywhere convergence. 
  
  \begin{proposition}\label{Appendix.Prop.S01}
  	If $f_n \overset{\rm S}\to f$ in $D_S(\mathbb{R}_+;\mathbb{R})$ as $n\to\infty$, then the following hold.
  	\begin{enumerate}
  		\item[(1)]  $f_n(t) \to f(t)$ for every $t\geq 0$ outside a countable set. 
  		
  		\item[(2)] $g*f_n\to g*f$ point-wisely for any $g\in L^\infty_{\rm loc}(\mathbb{R}_+;\mathbb{R})$.
  		\end{enumerate}
  \proof The first claim comes from Remark 2.4 in \cite{Jakubowski1997}. For the second one, recall the sequences $\{ f^\epsilon_n \}_{n\geq 1,\epsilon>0}$ and $\{ f^\epsilon  \}_{\epsilon>0}$, we have for any $T\geq 0$, 
  \beqnn
  |g*f_n(t)- g*f(t)| \leq | g*(f_n- f^\epsilon_n)(t)|  + \|g*f^\epsilon_n(t)- g*f^\epsilon*(t)|   + |g*(f^\epsilon- f)(t)|  .   
  \eeqnn
  	By the first two inequalities in (\ref{eqn.D01}), 
  	\beqnn
  	| g*(f_n- f^\epsilon_n)(t)| 
  	+|g*(f^\epsilon- f)(t)|  \leq  \| g\|_{L^1_t} \cdot \big( \| f_n- f^\epsilon_n\|_{L^\infty_t} + \|  f^\epsilon-f\|_{L^\infty_t} \big) \to 2\epsilon \cdot  \| g\|_{L^1_t}.
  	\eeqnn
  	By using  Fubini's theorem and then the limit in  (\ref{eqn.D01}),
  	\beqnn
  	|g*f^\epsilon_n(t)- g*f^\epsilon(t)|  = | \mathcal{I}_g*df^\epsilon_n(t)- \mathcal{I}_g*f^\epsilon(t)| \to 0,
  	\eeqnn
  	as $n\to\infty$ for any $\epsilon >0$ and $t\geq 0$. The second claim follows by putting these estimates together.
  	\qed 
  	 
  \end{proposition}
  
  For each $a< b$ and $t\geq 0$, let  $N_a^b(f,t)$ be the up-crossing number of $[a,b]$ by the function $f \in D(\mathbb{R}_+;\mathbb{R})$ in the time interval $[0,t]$. 
  
  \begin{proposition} \label{Appendix.Prop.S02}
  	For a sequence of stochastic processes $\{ X_n \}_{n\geq 1}$ in $D_S(\mathbb{R}_+;\mathbb{R})$, it is relatively compact if and only if both $\{\|X_n\|_{L^\infty_T} \}_{n\geq 1}$ and $\{N_a^b(X_n,T)\}$ are tight for any $T\geq 0$, which, in particular,  hold if $\{ X_n \}_{n\geq 1} \subset V(\mathbb{R}_+;\mathbb{R})$ and $\sup_{n\geq 1}\mathbf{E}[|X_n(t)|]<\infty$ for any $t\geq 0$. 
  \end{proposition}

  \begin{proposition}\label{Appendix.Prop.S020}
  	For a sequence of processes $\{X_n \}_{n\geq 1}$, if  $X_n \to  X$ weakly in $D_S(\mathbb{R}_+;\mathbb{R})$, 
   then the following hold.
  	\begin{enumerate}
  		\item[(1)] For a countable set $Q\in\mathbb{R}_+$, we have $ X_n \overset{\rm f.f.d.}\to  X$  along $\mathbb{R}_+\setminus Q$.
  		
  		\item[(2)]  If $X_n$ is non-decreasing for each $n\geq 1$, so is $X$. 
  		
  		\item[(3)] There exist a sequence $\{ \hat{X}_n \}_{n\geq 1}$ and $\hat{X}$ defined on a common probability space with
  		\beqnn
  		\hat{X}_n \overset{\rm d}= X_n
  		\quad\mbox{and}\quad
  		\hat{X}  \overset{\rm d}= X , \quad n\geq 1,
  		\eeqnn
  		such that $\hat{X}_n \overset{\rm a.s.}\to \hat{X}$  in $D_S(\mathbb{R}_+;\mathbb{R})$. 
  		
  		\item[(4)] For any $g\in L^\infty_{\rm loc}(\mathbb{R}_+;\mathbb{R})$, we have $g*X^{(n)} \overset{\rm f.f.d.}\longrightarrow g*X$ along $\mathbb{R}_+$.
  		\end{enumerate}
  	 
  	\end{proposition}

  \begin{proposition}\label{Appendix.Prop.S021}
  If two sequences of processes $\{ X_n \}_{n\geq 1}$ and  $\{ Z_n \}_{n\geq 1}$ are relatively compact in $D_S(\mathbb{R}_+;\mathbb{R})$,  
  then  $\{ (X_n,Z_n) \}_{n\geq 1}$ is relatively compact in $D_S(\mathbb{R}_+;\mathbb{R}\times \mathbb{R})$.  
  \end{proposition}

  \begin{proposition}\label{Appendix.Prop.S03}
  Suppose that a sequence of processes $\{ X_n \}_{n\geq 1}$ is relatively compact in $D_S(\mathbb{R}_+;\mathbb{R})$ and for any finite sequence $ t_1,\cdots,t_k $ of a dense subset $A\subset\mathbb{R}_+$, 
  \beqnn
   (X_n(t_1),\cdots, X_n(t_k))\overset{\rm d}\to \nu_{t_1,\cdots,t_k},
  \eeqnn
  as $n\to\infty$ with $\nu_{t_1,\cdots,t_k}$ being a probability law on $\mathbb{R}^k$, then there exists process $X$ in $D(\mathbb{R}_+;\mathbb{R})$ such that 
  \begin{enumerate}
  	\item[(1)]  $(X(t_1),\cdots, X(t_k)) \overset{\rm d}=\nu_{t_1,\cdots,t_k}$ for any finite sequence $ \{t_1,\cdots,t_k\}\subset A$;
  	
  	\item[(2)]  $ X_n \to  X$ weakly in $D_S(\mathbb{R}_+;\mathbb{R})$.
  	\end{enumerate} 
  
  \end{proposition}

  \begin{proposition}  \label{Appendix.Prop.S04}
  A sequence of uniformly tight martingales $\{ Y_n \}_{n\geq 1}$ is relatively compact in  $D_S(\mathbb{R}_+;\mathbb{R})$. 
  \end{proposition}
   \begin{proposition} \label{Appendix.Prop.S07}
   For a sequence of relatively compact martingales $\{ Y_n \}_{n\geq 1}$ is in  $D_S(\mathbb{R}_+;\mathbb{R})$, we have each limit process is a square-integral martingale if 
  \beqnn
  \sup_{n\geq 1}\mathbf{E}\bigg[\sup_{t\in[0,T]}\big|Y_n(t)\big|^2\bigg]<\infty,\quad T\geq 0 . 
  \eeqnn 
  \end{proposition}
  

     \renewcommand{\theequation}{D.\arabic{equation}}
    
 \section{Extended Bernstein functions}
 \label{Appendix--EBF}
 \setcounter{equation}{0}
 
 In this section, we recall and provide some elementary properties of extended Bernstein functions. 
 The reader may refer to \cite{SchillingSongVondracek2012} for a detailed discussion of Bernstein functions.
 A real valued function $f$ on $\mathbb{R}_+$ is a \textit{extended Bernstein function} if it is smooth on $(0,\infty)$ and
 	\beqnn 
 	(-1)^{k-1} \frac{d^k}{d\lambda^k}f(\lambda)\geq 0,\quad k\in \mathbb{Z}_+, \, \lambda >0.
 	\eeqnn
 Additionally, if $f(\lambda)\geq 0$ for all $\lambda>0$ then $f$ is said to be a \textit{Bernstein function}. 
 The spaces of all Bernstein functions and extended Bernstein functions are denoted by $\mathcal{BF}$ and $\mathcal{EBF}$ respectively.  
 For $k,m\in\mathbb{Z}_+$, let $\mathcal{EBF}^{k\times m}$ be the space of all functions $\boldsymbol{f} \in C(\mathbb{R}_+;\mathbb{R}^{k\times m})$ with $f_{ij} \in \mathcal{EBF}$.
 It is obvious that every $f \in \mathcal{EBF}$ is non-decreasing and $f +\beta  \in \mathcal{BF}$ for any constant $\beta\geq| f(0)\wedge 0 |$.

 \begin{proposition}\label{Prop.EBF01}
 	A function $f\in\mathcal{EBF}$  if and only if it admits the representation
 	\beqlb\label{eqn.RepEBF}
 	f(\lambda)= b+ \sigma \cdot \lambda + \int_{(0,\infty)} (1-e^{-\lambda t})\nu(dt),
 	\eeqlb
 	for two constants $b \in\mathbb{R}$, $\sigma\geq 0$ and a $\sigma$-finite measure $\nu(dt)$ on $(0,\infty)$ satisfying
 	\beqlb \label{eqn.A02}
 	\int_{(0,\infty)} (1\wedge t)\nu(dt)<\infty.
 	\eeqlb
  	In particular, the triplet $(b,c,\nu)$ is usually known as {\rm L\'evy triplet} and determines the Bernstein function $f$ uniformly and vice versa.  
  	 Moreover, $f\in \mathcal{BF}$ if and only if $b\geq 0$.
 \end{proposition}
 \proof If $f(0)\geq 0$ then $f\in \mathcal{BF}$ and the representation (\ref{eqn.RepEBF}) follows directly from Theorem~3.2 in \cite{SchillingSongVondracek2012} with $b=f(0)\geq 0$. 
 If $f(0)\leq 0$, the fact that $f$ is non-decreasing yields that  $f(\lambda)-f(0)\geq 0$ and hence $f(\lambda)-f(0)\in \mathcal{BF}$. Using Theorem~3.2 in \cite{SchillingSongVondracek2012} again, the function $f(\lambda)-f(0) $ admits the representation (\ref{eqn.RepEBF}) for some L\'evy triplet $(b,c,\nu)$ and hence the function $f(\lambda)$ also admits  the representation (\ref{eqn.RepEBF}) with L\'evy triplet $(b+f(0),c,\nu)$. 
 \qed

 \begin{proposition} \label{Prop.EBF02}
 	For a sequence $\{f_n \}_{n\geq 1} \subset  \mathcal{EBF}$, if the limit $\lim_{n\to\infty}f_n(\lambda)=f(\lambda)$ exists for any $\lambda \geq 0$, then $f \in \mathcal{EBF}$.
  Additionally,	if $f_n$ and $f$ have L\'evy triplets $ (b_n,c_n,\nu_n) $ and $(b,c,\nu)$ respectively,  then
 	\beqlb\label{eqn.A01}
 	b_n\to b,\quad c_n + \int_0^\infty (t\wedge 1) \nu_n(dt)\to c+ \int_0^\infty (t\wedge 1) \nu(dt) 
 	\quad \mbox{and}\quad \int_0^\infty g(t) \nu_n(dt)\to  \int_0^\infty g(t) \nu(dt) 
 	\eeqlb
  as $n\to\infty$,  for any bounded and continuous $g$ on $\mathbb{R}_+$ satisfying that $g(x)=O(x)$ as $x\to0+$.  
 \end{proposition}
 \proof For each $n\geq 1$ and $\lambda \geq 0$, let  $h_n(\lambda)=f_n(\lambda)- f_n(0)$.
 By Proposition~\ref{Prop.EBF01}, we have $h_n \in \mathcal{BF}$ with L\'evy triplet $(0,c_n,\nu_n)$ and as $n\to\infty$,
 \beqnn
 b_n= f_n(0)\to  f(0)=b
 \quad \mbox{and}\quad
 h_n(\lambda) \to f(\lambda)- f(0),\quad \lambda \geq 0.
 \eeqnn 
 By the first claim in \cite[Corollary~3.9]{SchillingSongVondracek2012}, we have $f(\lambda)- f(0) \in \mathcal{BF}$ and hence $f(\lambda)  \in \mathcal{EBF}$. 
 Notice that and  $f(\lambda)- f(0)$ has $(0,c,\nu)$. 
 The last two limits in (\ref{eqn.A01}) are direct consequences of the second claim in \cite[Corollary~3.9]{SchillingSongVondracek2012}.
 \qed

  \begin{proposition} \label{Prop.EBF03}
 	For any given function $f\in \mathcal{EBF}$, there exists a sequence of scaling parameters $\{\gamma_n\}_{n\geq 1}$ with $\gamma_n=O(n)$ as $n\to\infty$ and a sequence $\{\varphi_n\}_{n\geq 1} \subset L^1(\mathbb{R}_+;\mathbb{R}_+)$ such that 
 	\beqlb\label{eqn.A03}
   \lim_{n\to\infty}	\gamma_n \Big(1-\mathcal{L}_{\varphi_n}(\lambda/n)\Big) =  f(\lambda),\quad \lambda \geq 0. 
 	\eeqlb
 \end{proposition}
 \proof 
 Assume that  the function $f\in \mathcal{EBF}$ has L\'evy triplet $(a,b,\nu)$. 
 Let $\bar\nu(t):=\nu([t,\infty))$ be the tail-measure of $\nu$, which is non-increasing and right-continuous. 
 For each $n\geq 1$, let $\{\bar\nu_n(t):t\geq 0\}_{n\geq 1}$ be a sequence of non-increasing and absolutely continuous approximation for $\bar\nu$ with the Radon-Nikodym derivative $\bar\nu_n$, denoted by $\bar\nu_n'$, satisfying that $\bar\nu_n'(t) \equiv 0$ for any $t\in[0,1/n]$ and 
 \beqlb\label{eqn.001}
 \lim_{n\to\infty} \int_0^\infty h(t) \bar\nu_n'(t)dt = \lim_{n\to\infty} \int_0^\infty h(t) d\bar\nu_n(t) =  \int_0^\infty h(t) d\bar\nu (t)= - \int_0^\infty h(t) \nu (dt),
 \eeqlb
 for any $h\in C_b(\mathbb{R}_+)$ satisfying that $h(x)=O(x)$ as $x\to0+$. 
 For each $n\geq 1$, we consider a function $\varphi_n$ defined on $(0,\infty)$ as follows
 \beqnn
 \varphi_n(t)= \Big( 1-\frac{a}{\gamma_n} \Big) \cdot \Big( (1-\varepsilon_n) \cdot \beta_n e^{-\beta_n t}  + \varepsilon_n \cdot g_n(t) \Big),
 \eeqnn
 for three positive constants $\gamma_n>0$, $\beta_n\geq 0$, $\varepsilon_n\in(0,1)$ and a probability density function $g_n$ on $(0,\infty)$ that will be specified later. 
 
 \textbf{Case~1.} If $ \nu (\mathbb{R}_+)=\infty$, we choose 
 \beqnn
 \gamma_n = 2\cdot \bar\nu_n(0),\quad  
 \varepsilon_n \equiv \frac{1}{2}\cdot \mathbf{1}_{\{b>0\}}+ \mathbf{1}_{\{b=0\}},\quad 
 \beta_n =\frac{ \bar\nu_n(0)}{b\cdot n}\cdot  \mathbf{1}_{\{b>0\}}
 \quad\mbox{and}\quad 
 g_n(t) = -\frac{\bar\nu_n'(t/n)}{n\cdot \bar\nu_n(0)}.
 \eeqnn
 Using integration by parts to   (\ref{eqn.A02}), we have $\gamma_n =O(n)$ as $n\to\infty$.  
 A simple calculation shows that 
 \beqnn
	\gamma_n \Big(1-\mathcal{L}_{\varphi_n}(\lambda/n)\Big)
 \ar=\ar a + \Big( 1-\frac{a}{\gamma_n} \Big) \cdot \Big( \frac{ \bar\nu_n(0)}{ \bar\nu_n(0)+b\lambda}\cdot b\lambda - \int_0^\infty \big(1-e^{- \lambda  t}\big)  \bar\nu_n'(t)dt\Big) , \quad \lambda \geq 0. 
 \eeqnn
 The desired limit (\ref{eqn.A03}) follows from this and (\ref{eqn.001}).
 
 \textbf{Case~2.} If $ \nu (\mathbb{R}_+)<\infty$,   we choose 
 \beqnn
 \gamma_n = n,\quad  
 \varepsilon_n = \frac{\nu (\mathbb{R}_+)}{n} \wedge 1,\quad 
 \beta_n =\frac{1}{b }\cdot  \mathbf{1}_{\{b>0\}} + n \cdot  \mathbf{1}_{\{b>0\}}
 \quad\mbox{and}\quad 
 g_n(t) = -\frac{\bar\nu_n'(t/n)}{n\cdot \bar\nu_n(0)}.
 \eeqnn
 In this case, we have  for each $n\geq \nu (\mathbb{R}_+)$ and $\lambda \geq 0$,
 \beqnn
 	\gamma_n \Big(1-\mathcal{L}_{\varphi_n}(\lambda/n)\Big)
 \ar=\ar a + \Big( 1-\frac{a}{n} \Big) \cdot \Big( \frac{  \lambda}{ \beta_n+\lambda/n} - \int_0^\infty \big(1-e^{- \lambda  t}\big)  \bar\nu_n'(t)dt\Big).
 \eeqnn
 By using (\ref{eqn.001}) again, we can obtain the limit (\ref{eqn.A03}) directly. 
 \qed

 \begin{corollary} \label{Corollary.EBF01}
  For $k,m\in\mathbb{Z}_+$, $\boldsymbol{A}\in \mathbb{R}_+^{k\times m}$ and  $\boldsymbol{f} \in \mathcal{EBF}^{k\times m}$, there exists a sequence of scaling parameters $\{\gamma_n\}_{n\geq 1}$ with $\gamma_n=O(n)$ as $n\to\infty$ and a sequence of functions $\{\boldsymbol{\psi}_n\}_{n\geq 1} \subset L^1(\mathbb{R}_+;\mathbb{R}_+^{k\times m})$
   	\beqlb\label{eqn.A06}
   \lim_{n\to\infty}	\gamma_n \Big(\boldsymbol{A}-\mathcal{L}_{\boldsymbol{\psi}_n}(\lambda/n)\Big)
   =  \boldsymbol{f}(\lambda),\quad \lambda \geq 0. 
   \eeqlb 
 \end{corollary}
 \proof For each $1\leq i\leq k $ and $1\leq j\leq m$, by Proposition~\ref{Prop.EBF03} there exists a sequence of scaling parameters $\{\gamma_n^{ij}\}_{n\geq 1}$ with $\gamma_n^{ij}=O(n)$ and a sequence $\{\varphi_n^{ij}\}_{n\geq 1} \subset L^1(\mathbb{R}_+;\mathbb{R}_+)$ such that 
 \beqlb \label{eqn.A07}
 \lim_{n\to\infty}	\gamma_n^{ij} \Big(A_{ij}-\mathcal{L}_{A_{ij}\cdot\varphi_n^{ij}}(\lambda/n)\Big)
 =\lim_{n\to\infty}	A_{ij}\cdot \gamma_n^{ij} \Big( 1 -\mathcal{L}_{\varphi_n^{ij}}(\lambda/n)\Big) 
  = f_{ij}(\lambda),\quad \lambda \geq 0. 
 \eeqlb
 We choose a sequence $\{\gamma_n\}_{n\geq 1}$ satisfying that $ \gamma_n = O(n)$ as $n\to\infty$ and
 \beqnn 
 \sup_{n\geq 1} \frac{\gamma_n^{ij}}{\gamma_n }<1,\quad  1\leq i\leq k ,\, 1\leq j\leq m . 
 \eeqnn
 Consider a sequence of functions $\{\boldsymbol{\psi}_n\}_{n\geq 1} \subset L^1(\mathbb{R}_+;\mathbb{R}_+^{k\times m})$ defined by
 \beqnn
 \psi_{n,ij}(t):=\Big(1- \frac{\gamma_n^{ij}}{\gamma_n } \Big) \cdot A_{ij} \cdot \beta_n e^{-\beta_n t} + \frac{\gamma_n^{ij}}{\gamma_n }\cdot A_{ij}\cdot \varphi_n^{ij}(t),\quad t\geq 0,\, n\geq 1,\,1\leq i\leq k ,\, 1\leq j\leq m,
 \eeqnn
 where  $\{ \beta_n\}_{n\geq 1}$ a sequence of constants with $\beta_n\to \infty$ as $n\to\infty$. 
 For each $\lambda\geq  0$, a simple calculation together with (\ref{eqn.A07}) shows that as $n\to\infty$,
 \beqnn
 \gamma_n \Big(A_{ij}-\mathcal{L}_{\psi_{n,ij}}(\lambda/n)\Big) = (\gamma_n -\gamma_n^{ij}) \cdot A_{ij} \cdot \Big( 1- \frac{n\beta_n}{n\beta_n + \lambda} \Big) 
 +	\gamma_n^{ij} \Big(A_{ij}-\mathcal{L}_{A_{ij}\cdot\varphi_n^{ij}}(\lambda/n)\Big)
 \to f_{ij}(\lambda).
 \eeqnn  
 \qed

   \begin{proposition} \label{Prop.EBF011}
 	For each $f \in \mathcal{EBF}$ with $f(\infty)>0$,  there exists a unique $\sigma$-finite measure $\mu^f(dt)$ on $\mathbb{R}_+$ such that for any $\lambda \geq 0$ with $f(\lambda)>0$,
 	\beqlb\label{eqn.App.A04}
 	\int_0^\infty e^{-\lambda t}\mu^f(dt) = \frac{1}{f(\lambda)}.
 	\eeqlb 
 \end{proposition}
 \proof If $f(0)\geq 0$,  the desired result follows directly from (5.20) in \cite[p.63]{SchillingSongVondracek2012}. 
 If $f(0)<0$, let $\lambda_0:= \inf\{ \lambda\geq 0: f(\lambda)\geq 0 \}$. Notice that $f(\lambda +\lambda_0)$ is a Bernstein function, by the previous result there exists a unique $\sigma$-finite measure $\mu^f_0(dt)$ on $\mathbb{R}_+$ such that 
 \beqnn
 \int_0^\infty e^{-\lambda t}\mu^f_0(dt) = \frac{1}{f(\lambda+\lambda_0)},\quad \lambda \geq 0
 \eeqnn
 and hence (\ref{eqn.App.A04}) with $\mu^f (dt)= e^{\lambda_0 t} \cdot \mu^f_0(dt)$.
 \qed 
 
 When $f \in \mathcal{BF}$, it is usual to refer to $\mu^f$ as the \text{potential measure}  associated with $f$ in the potential theory. 
 Without confusion, we also use this terminology for convention. 
 
 \begin{proposition}\label{Prop.ResolventEquation}
 	For two functions $f_1,f_2 \in \mathcal{EBF}$ with L\'evy triplets $(b_1,\sigma,\nu)$ and $(b_2,\sigma,\nu)$ respectively, assume that $f_1(\infty),f_2(\infty)>0$,   the corresponding extended potential measures $\mu^{f_1}(dt)$ and $\mu^{f_2}(dt)$ satisfy
 	\beqlb\label{ExResolventEquation}
 	\mu^{f_1} -\mu^{f_2} +(b_1-b_2)\cdot \mu^{f_1}*\mu^{f_2}  =0.
 	\eeqlb 
 \end{proposition}
 \proof By (\ref{eqn.App.A04}), for any $\lambda \geq 0$ with $f_1(\lambda),f_2(\lambda)\geq 0 $ we have 
 \beqnn
 \frac{1}{f_1(\lambda)} -\frac{1}{f_2(\lambda)} +(b_1-b_2)\cdot \frac{1}{f_1(\lambda)}\cdot \frac{1}{f_2(\lambda)}
 =0.
 \eeqnn
 The desired equation (\ref{ExResolventEquation}) follows  directly from the one-to-one correspondence among measures on $\mathbb{R}_+$ and their Laplace transforms.
 \qed

%
%
%
%
%
%
%
%
%
%
%
%
%
%
%
%


 \bibliographystyle{plain}

 \bibliography{Reference}

 \end{document}